\documentclass[11pt,oneside]{amsart}

\usepackage{bm}

\usepackage{fix-cm}
\DeclareMathAlphabet{\mathcal}{OMS}{cmsy}{m}{n}

\usepackage{tikz-cd}

\usepackage{tcolorbox}

\usepackage[arrow,curve,matrix,tips,2cell]{xy}

\usepackage[hmargin=2cm,vmargin=1.5cm]{geometry}

\usepackage{fancyhdr}

\usepackage{anyfontsize}
\usepackage{mathrsfs}

\usepackage{accents}

\usepackage{bbm}

\usepackage{amsmath}
\usepackage{amscd}
\usepackage{amssymb}
\usepackage{latexsym}
\usepackage{url}

\usepackage{graphicx}

\newtheorem{theorem}{Theorem}[section]
\newtheorem*{theorem*}{Theorem}
\newtheorem{lemma}[theorem]{Lemma}
\newtheorem*{lemma*}{Lemma}
\newtheorem{corollary}[theorem]{Corollary}
\newtheorem{proposition}[theorem]{Proposition}

\newtheorem{remark}[theorem]{Remark}
\newtheorem{definition}[theorem]{Definition}
\newtheorem*{definition*}{Definition}

\newtheorem{question}[theorem]{Question}
\newtheorem*{question*}{Question}

\newtheorem{example}[theorem]{Example}
\newtheorem{examples}[theorem]{Examples}

\newtheorem{introtheorem}{Theorem}

\makeatletter
\def\revddots{\mathinner{\mkern1mu\raise\p@
\vbox{\kern7\p@\hbox{.}}\mkern2mu
\raise4\p@\hbox{.}\mkern2mu\raise7\p@\hbox{.}\mkern1mu}}
\makeatother 
\newcommand{\bgl}{\begin{equation}} 
\newcommand{\egl}{\end{equation}}
\newcommand{\bgloz}{\begin{equation*}} 
\newcommand{\egloz}{\end{equation*}}
\newcommand{\bgln}{\begin{eqnarray}} 
\newcommand{\egln}{\end{eqnarray}}
\newcommand{\bglnoz}{\begin{eqnarray*}} 
\newcommand{\eglnoz}{\end{eqnarray*}}
\newcommand{\btheo}{\begin{theorem}}
\newcommand{\etheo}{\end{theorem}}
\newcommand{\btheooz}{\begin{theorem*}}
\newcommand{\etheooz}{\end{theorem*}}

\newcommand{\blemma}{\begin{lemma}}
\newcommand{\elemma}{\end{lemma}}
\newcommand{\blemmaoz}{\begin{lemma*}}
\newcommand{\elemmaoz}{\end{lemma*}}
\newcommand{\bproof}{\begin{proof}}
\newcommand{\eproof}{\end{proof}}
\newcommand{\bbew}{\begin{beweis}}
\newcommand{\ebew}{\end{beweis}}
\newcommand{\bremark}{\begin{remark}\em}
\newcommand{\eremark}{\end{remark}}
\newcommand{\bdefin}{\begin{definition}}
\newcommand{\edefin}{\end{definition}}
\newcommand{\bdefinoz}{\begin{definition*}}
\newcommand{\edefinoz}{\end{definition*}}
\newcommand{\bex}{\begin{example}}
\newcommand{\eex}{\end{example}}
\newcommand{\bexs}{\begin{examples}}
\newcommand{\eexs}{\end{examples}}

\newcommand{\bprop}{\begin{proposition}}
\newcommand{\eprop}{\end{proposition}}
\newcommand{\bcor}{\begin{corollary}}
\newcommand{\ecor}{\end{corollary}}
\newcommand{\bfa}{\begin{cases}} 
\newcommand{\efa}{\end{cases}}

\newcommand{\bquestion}{\begin{question}}
\newcommand{\equestion}{\end{question}}
\newcommand{\bquestionoz}{\begin{question*}}
\newcommand{\equestionoz}{\end{question*}}
\newcommand{\cF}{\mathcal F}
\newcommand{\cG}{\mathcal G}
\newcommand{\cH}{\mathcal H}

\newcommand{\cJ}{\mathcal J}
\newcommand{\cK}{\mathcal K}
\newcommand{\cL}{\mathcal L}

\newcommand{\cO}{\mathcal O}

\newcommand{\cR}{\mathcal R}

\newcommand{\cV}{\mathcal V}
\newcommand{\cW}{\mathcal W}

\def\Nz{\mathbb{N}}

\def\Rz{\mathbb{R}}

\def\Zz{\mathbb{Z}}

\def\1z{\mathbb{1}}
\newcommand{\fC}{\mathfrak C}
\newcommand{\fD}{\mathfrak D}

\newcommand{\fH}{\mathfrak H}

\newcommand{\fS}{\mathfrak S}
\newcommand{\fT}{\mathfrak T}

\newcommand{\mfd}{\mathfrak d}

\newcommand{\mff}{\mathfrak f}

\newcommand{\mft}{\mathfrak t}

\newcommand{\mfv}{\mathfrak v}
\newcommand{\mfw}{\mathfrak w}
\newcommand{\mfy}{\mathfrak y}
\newcommand{\an}[1]{``#1''} 
\newcommand{\ti}{\tilde}

\newcommand{\ma}{\mapsto} 
\newcommand{\onto}{\twoheadrightarrow} 
\newcommand{\into}{\hookrightarrow} 
\newcommand{\Rarr}{\Rightarrow} 
\newcommand{\Larr}{\Leftarrow} 
\newcommand{\LRarr}{\Leftrightarrow} 

\def\SEMI{\mbox{$\times\kern-2pt\vrule height5pt width.6pt \kern3pt $}}

\newcommand{\lcm}{{\rm lcm}} 
\newcommand{\reg}{^\times} 
\newcommand{\lspan}{{\rm span}} 
\newcommand{\clspan}{\overline{\lspan}} 
\newcommand{\defeq}{\mathrel{:=}} 

\newcommand{\dop}{\text{: }} 
\newcommand{\lge}{\left\{} 
\newcommand{\rge}{\right\}} 
\newcommand{\lsp}{\left\langle} 
\newcommand{\rsp}{\right\rangle} 
\newcommand{\gekl}[1]{\lge #1 \rge} 
\newcommand{\spkl}[1]{\lsp #1 \rsp} 
\newcommand{\menge}[2]{\gekl{ #1 \dop #2 }} 
\newcommand{\isom}{\xrightarrow{\raisebox{-1ex}[0ex][0ex]{$\sim$}}} 
\newcommand{\dom}{{\rm dom\,}}
\newcommand{\im}{{\rm im\,}}

\newcommand{\hcJ}{\widehat{\cJ}}

\newcommand{\rmr}{\ensuremath{\mathrm{r}}}
\newcommand{\rms}{\ensuremath{\mathrm{s}}}

\newcommand{\bltimes}{\bar{\ltimes}}
\newcommand{\bcJ}{\bar{\cJ}}
\newcommand{\bIl}{\bar{I}_l}
\newcommand{\bsim}{\bar{\sim}}

\newcommand{\mcm}{{\rm mcm}}

\newcommand{\Omegamax}{\Omega_{\max}}
\newcommand{\bOmega}{\partial \Omega}
\newcommand{\Omegainf}{\Omega_{\infty}}

\newcommand{\bcG}{\partial \cG}

\newcommand{\kerbd}{{\rm Ker}_{\partial}}

\newcommand{\oset}[2]{%
  \mathop{#2}\limits^{\vbox to -1.66ex{%
  \kern -1.4ex\hbox{$#1$}\vss}}}

\newcommand{\tipreceq}{\oset{\sim}{\preceq}}
\newcommand{\tiprec}{\oset{\sim}{\prec}}

\newcommand{\bbd}{\mathbbm{d}}

\usepackage{newtxtext}
\usepackage{newtxmath}

\begin{document}

\title[Left regular representations of Garside categories I]{Left regular representations of Garside categories I. \\ C*-algebras and groupoids}

\thispagestyle{fancy}

\author{Xin Li}

\address{Xin Li, School of Mathematics and Statistics, University of Glasgow, University Place, Glasgow G12 8QQ, United Kingdom}
\email{Xin.Li@glasgow.ac.uk}

\subjclass[2010]{Primary 46L05, 20F36; Secondary 46L55, 37A55}

\thanks{This project has received funding from the European Research Council (ERC) under the European Union's Horizon 2020 research
and innovation programme (grant agreement No. 817597).}

\begin{abstract}
We initiate the study of C*-algebras and groupoids arising from left regular representations of Garside categories, a notion which originated from the study of Braid groups. Every higher rank graph is a Garside category in a natural way. We develop a general classification result for closed invariant subspaces of our groupoids as well as criteria for topological freeness and local contractiveness, properties which are relevant for the structure of the corresponding C*-algebras. Our results provide a conceptual explanation for previous results on gauge-invariant ideals of higher rank graph C*-algebras. As another application, we give a complete analysis of the ideal structures of C*-algebras generated by left regular representations of Artin-Tits monoids.
\end{abstract}

\maketitle


\setlength{\parindent}{0cm} \setlength{\parskip}{0.5cm}

\section{Introduction}

C*-algebras generated by partial isometries form a rich class of examples, including C*-algebras attached to shifts of finite type \cite{Cun77,CK80,Cun81}, graph C*-algebras \cite{Rae}, higher rank graph C*-algebras \cite{KuPa}, C*-algebras attached to self-similiar groups \cite{Nek09} and semigroup C*-algebras \cite{Co,Li12,Li13,CELY}. For instance, it was shown in \cite{EP17} that every UCT Kirchberg algebra arises in this way. The class of UCT Kirchberg algebras plays an important role in the Elliott classification programme for C*-algebras (see \cite{KiPh,Phi,Kir,Ror}). Spielberg observed that all the classes of C*-algebras mentioned above can be viewed as special cases of a general, unifying construction of C*-algebras generated by left regular representations of left cancellative small categories \cite{Sp14,Sp20}. This is a very general construction, as it contains, up to Morita equivalence, all inverse semigroup C*-algebras (see \cite{DGKMW}). These C*-algebras come with a distinguished quotient which is called the boundary quotient. The passage from the C*-algebra to its boundary quotient is analogous to the passage from the Toeplitz-type C*-algebra of a shift of finite type or graph to its Cuntz-Krieger-type C*-algebra. 

A powerful way to study these C*-algebras of small categories is to construct a groupoid model and study properties of the C*-algebra through a detailed analysis of the groupoid \cite{Ren80,Sp14,Sp20}. Actually, there are two candidates for such groupoid models, which both arise from actions of an inverse semigroup on a space of certain filters attached to the small category. The inverse semigroup is given by the left inverse hull, i.e., the smallest inverse semigroup of partial bijections of the small category containing all left multiplication maps by individual elements of the small category. In \cite{Sp20}, a refined (and enlarged) version of the left inverse hull is considered, leading to the second groupoid model. In both cases, the filters which give rise to the unit space of the groupoid models are defined on the semilattice of idempotents of the inverse semigroup and take into account that elements of this semilattice are subsets of the original small category. The language of inverse semigroups provides an interpretation of the distinguished boundary quotient as the tight quotient, which is induced from the subspace of tight filters (see \cite{Ex08,Ex17,EP16}).

It is an interesting observation that in this very general framework, every left cancellative small category generates -- in an entirely natural and intrinsic way -- a dynamical system in terms of an inverse semigroup action or a groupoid. The same statement applies to the even more general setting of $0$-left cancellative semigroups as considered by Exel and Steinberg in \cite{ES18a,ES18b,ES18c,ES19}. Generally speaking, the goal would be to find a dictionary between properties of the small category, properties of the inverse semigroup action or groupoid, and properties of the C*-algebra and its boundary quotient. Indeed, we present criteria in terms of the underlying small category which completely characterize when the boundary groupoids -- which model the boundary quotients -- are Hausdorff, minimal or effective (or topologically free). We also establish a sufficient criterion for the boundary groupoid to be locally contractive. These properties have immediate consequences for the corresponding boundary quotient C*-algebras concerning ideal structure and pure infiniteness. Such criteria have been established in the general context of inverse semigroup actions and tight groupoids attached to inverse semigroups in \cite{EP16}, and it turns out to be fruitful to translate between the work in \cite{EP16} and our setting of small categories. For instance, this leads to generalizations of the results in \cite{OP}, which covers classes of finitely aligned small categories. In the special case of submonoids of groups, we are naturally led to the following characterization of topological freeness of the boundary action:
\begin{introtheorem}
\label{thm:intro_PinG}
Let $P$ be a submonoid of a group $G$ and denote by $G \curvearrowright \partial \Omega$ its boundary action (in the sense of \cite[Definition~5.7.8]{CELY}). Define $G^c \defeq \menge{g \in G}{(pP) \cap (gpP) \neq \emptyset \quad \forall \ p \in P}$. 
\setlength{\parindent}{0.5cm} \setlength{\parskip}{0cm}

Then $G \curvearrowright \partial \Omega$ is topologically free if and only if $G^c$ is the trivial group. In this case, $\partial C^*_{\lambda}(P)$ is simple, and $\partial C^*_{\lambda}(P)$ is purely infinite simple unless $P$ is the trivial monoid.
\end{introtheorem}
\setlength{\parindent}{0cm} \setlength{\parskip}{0cm}
$G^c$ is always a subgroup of $G$. Theorem~\ref{thm:intro_PinG} tells us that this subgroup captures topological freeness of the boundary action in an arguably more efficient way than the \an{core} as in \cite{CL07} (see also \cite[\S~5.7]{CELY}). In this form, with $G^c$ as the key ingredient, our characterization of topological freeness of the boundary action has not appeared before, but, as Marcelo Laca and Camila F. Sehnem kindly informed me, it also follows from \cite[Proposition~6.18]{LS}. We give a self-contained (and short) proof of Theorem~\ref{thm:intro_PinG} in \S~\ref{s:Boundary} (see Theorem~\ref{thm:Gc=1}).
\setlength{\parindent}{0.5cm} \setlength{\parskip}{0cm}

At the same time, our study of boundary groupoids arising from left regular representations of small categories led us to a characterization of topological freeness of tight groupoids attached to general inverse semigroups (see Theorem~\ref{thm:InvSgpBdTopFree}). To the best of the author's knowledge, such a characterization was not known before. 

We also clarify the relationship between the different groupoid models mentioned above and the analogous variations of the boundary groupoids. For the groupoids themselves, while minimality and local contractiveness are rather rare phenomena, we succeed in completely characterizing, in terms of the underlying small category, when the groupoids are Hausdorff or effective (or topologically free). Our criterion for topological freeness is inspired by \cite[Theorem~5.9]{LS}, which treats the special case of submonoids of groups. Furthermore, we establish a characterization when the boundary is the smallest non-empty closed invariant subspace of the character space, and determine in this case when the boundary groupoid is purely infinite (see Proposition~\ref{prop:smallest}).
\setlength{\parindent}{0cm} \setlength{\parskip}{0.5cm}

Having identified a natural and unifying general framework, it is important to find classes of small categories which are general enough so that they cover interesting classes of examples and yet concrete enough so that a detailed analysis is possible. 

The main goal of the present paper is to discuss one such class of small categories called Garside categories, and in this way contribute to our understanding of C*-algebras attached to small categories. The idea behind Garside categories originated from the study of Braid groups and monoids, and of the more general Artin-Tits groups and monoids. Roughly speaking, Garside structures allow us to carry over classical results and methods from Braid groups and monoids to more general groups, monoids or small categories. The concept of Garside categories feature in proofs of the $K(\pi,1)$-conjecture for various classes of groups \cite{Bes,Par,PS}. Recently, a connection has been discovered between Garside categories and Helly graphs, which has several applications, for instance to isomorphism conjectures such as the Farrell-Jones conjecture or the coarse Baum-Connes conjecture \cite{HO}. We refer the reader to \cite{Deh15} for more details on Garside categories.
\setlength{\parindent}{0.5cm} \setlength{\parskip}{0cm}

In our context, Garside structures allow us to establish normal forms for filters which form the unit spaces of our groupoids. This in turn leads to very concrete descriptions of the groupoid models themselves. As a result, we succeed in describing all closed invariant subspaces in terms of the underlying small category. 
\setlength{\parindent}{0cm} \setlength{\parskip}{0cm}

\begin{introtheorem}
Let $\fC$ be a finitely aligned, left cancellative, countable small category and $\fS$ a Garside family in $\fC$ with $\fS \cap \fC^* = \emptyset$ which is $=^*$-transverse and locally bounded. Let $I_l \ltimes \Omega$ be the groupoid model for $C^*_{\lambda}(\fC)$. 
\setlength{\parindent}{0.5cm} \setlength{\parskip}{0cm}

There is a one-to-one correspondence between closed invariant subspaces of $I_l \ltimes \Omega$ and admissible, $H$-invariant, $\max_{\preceq}^{\infty}$-closed pairs $(\fT,\fD)$ with $\fT \subseteq \fS$ and $\fD \subseteq \fC^0$.
\end{introtheorem}
\setlength{\parindent}{0cm} \setlength{\parskip}{0cm}

The reader will find more explanations and details in \S~\ref{s:C*Gar} (see Theorem~\ref{thm:ClInvSubsp}). The point is that our description is purely in terms of the Garside family $\fS$. We also explicitly characterize which of these closed invariant subspaces belong to the boundary. In addition, we establish criteria for topological freeness and local contractiveness. Again, these properties have consequences for ideal structure and pure infiniteness of our C*-algebras. Our analysis is made possible by the key property of Garside categories that every element admits a normal form, generalizing the classical normal form (also called greedy, Garside or Thurston normal form) of elements in Braid and Artin-Tits monoids. Indeed, as explained in \cite{Deh15}, the general notion of Garside categories (as in \cite{Deh15}) has been designed to allow for this kind of normal forms. For the purpose of studying groupoids and C*-algebras, the usefulness of normal forms has been observed already, for instance in the context of semigroup C*-algebras of right-angled or spherical Artin-Tits monoids \cite{CL02,CL07,LOS}, or of Baumslag-Solitar monoids \cite{Sp12}.
\setlength{\parindent}{0cm} \setlength{\parskip}{0.5cm}

As particular examples, we discuss higher rank graphs in \S~\ref{ss:k-graphs}. Actually, the starting point for this paper was the observation that every higher rank graph is a Garside category in a very natural way. Our results lead to a new interpretation of gauge-invariant ideals (see Lemma~\ref{lem:gauge=induced}). Moreover, not only do our results cover the C*-algebras of higher rank graphs, but they also treat Toeplitz algebras. Furthermore, our analysis extends to categories arising from self-similar actions on graphs or higher rank graphs. As another class of concrete examples, we discuss general Artin-Tits monoids. We complete the study of the ideal structure of their semigroup C*-algebras, which has been started in \cite{CL02,CL07,LOS}, by proving the following result:
 
\begin{introtheorem}
\label{thm:C}
Let $P$ be an irreducible Artin-Tits monoid with set of atoms $A$. If $P$ is spherical, then $\kerbd = \cK(\ell^2 P)$ if $\# A = 1$ and $\cK(\ell^2 P)$ is the only non-trivial ideal of $\kerbd$ if $2 \leq \#A < \infty$. In the latter case, $\kerbd / \cK(\ell^2 P)$ is purely infinite simple. If $P$ is not finitely generated and left reversible, then $\kerbd$ is purely infinite simple. If $P$ is finitely generated and not spherical, then $\cK(\ell^2 P)$ is the only non-trivial ideal of $C^*_{\lambda}(P)$, and $C^*_{\lambda}(P) / \cK(\ell^2 P)$ is purely infinite simple. If $P$ is not finitely generated and not left reversible, then $C^*_{\lambda}(P)$ is purely infinite simple.
\end{introtheorem}
\setlength{\parindent}{0cm} \setlength{\parskip}{0cm}

Here $\kerbd$ is the kernel of the canoncial projection $C^*_{\lambda}(P) \onto \partial C^*_{\lambda}(P)$. In the spherical or left reversible case, $\partial C^*_{\lambda}(P)$ coincides with the reduced group C*-algebra of the Artin-Tits group corresponding to $P$. In Theorem~\ref{thm:C}, the finitely generated, spherical case is treated in \cite{LOS}, and the right-angled case is treated in \cite{CL02,CL07}. Our contribution concerns the remaining cases. We can also characterize when $C^*_{\lambda}(P)$ or $\kerbd$ is nuclear (see also \cite[Theorem~4.2]{LL}). Moreover, we point out that K-theory for semigroup C*-algebras of Artin-Tits monoids has been computed in \cite{Li20}, assuming that the corresponding Artin-Tits group satisfies the Baum-Connes conjecture with coefficients.
\setlength{\parindent}{0cm} \setlength{\parskip}{0.5cm}

Higher rank graphs and Artin-Tits monoids are just some examples of Garside categories. The reader will find many more examples in \cite{Deh15}.

Apart from providing a natural class of examples where we can test and develop our understanding of C*-algebras attached to small categories, this paper at the same time sets the stage for a detailed analysis of the groupoids arising from left regular representations of small categories. These groupoids are not only auxiliary structures to translate between small categories and their C*-algebras, but they are also interesting on their own right as they lead to interesting new structures, for instance topological full groups. Our original motivation which led to the present paper was the natural question left open by Matui in \cite[\S~5.3]{Mat16} whether topological full groups of groupoids attached to products of shifts of finite type are of type ${\rm F}_{\infty}$. We answer this question in \cite{Li21b}.

I am indebted to Marcelo Laca and Camila F. Sehnem for pointing out to me that Theorem~\ref{thm:intro_PinG} also follows from their work in \cite{LS}. I would also like to thank Chris Bruce and the anonymous referee for helpful comments which improved the paper.

This paper is a contribution to the special issue of the Glasgow Mathematical Journal on the occasion of the BMC/BAMC 2021, and I would like to thank the organizers for creating such a wonderful and highly successful event.

\section{Preliminaries}

Let us recall some basics regarding left regular representations of left cancellative categories, C*-algebras generated by these representations and groupoid models for these C*-algebras. Note that we view categories -- which will all be assumed to be small in this paper -- as generalizations of monoids (as in \cite{Wit}), so that no sophisticated category theory will be used.

\subsection{Left cancellative small categories, their left regular representations and C*-algebras}
\label{ss:cat-rr-C}

Given a small category with set of morphisms $\fC$, let $\fC^0$ be its set of objects. We will identify $\mfv \in \fC^0$ with the identity morphism at $\mfv$, so that $\fC^0$ is identified with a subset of $\fC$. Often, we will abuse notation and simply call $\fC$ the small category. Let $\mfd: \: \fC \to \fC^0$ and $\mft: \: \fC \to \fC^0$ be the domain and target maps, so that for $c, d \in \fC$, the product $cd$ is defined if and only if $\mfd(c) = \mft(d)$. This means that our convention is the same as the one in \cite{Sp20,Wit}, while it is opposite to the one used in \cite{Deh15} (see \cite[Remark~1.1]{Wit}). For $c \in \fC$ and $S \subseteq \fC$, we set $cS \defeq \menge{cs}{s \in S, \, \mft(s) = \mfd(c)}$. Moreover, $\fC^*$ denotes the set of invertible elements of $\fC$, i.e., elements $c \in \fC$ for which there exists $c^{-1} \in \fC$ with $c^{-1} c = \mfd(c)$ and $c c^{-1} = \mft(c)$. Note that $\fC^*$ is denoted by $\fC\reg$ in \cite{Deh15,Wit}.
\bdefin
A small category $\fC$ is called left cancellative if for all $c, x, y \in \fC$ with $\mfd(c) = \mft(x) = \mft(y)$, $cx = cy$ implies $x = y$.
\edefin
\setlength{\parindent}{0cm} \setlength{\parskip}{0cm}

From now on, all our small categories will be assumed to be left cancellative. Let $\fC$ be such a small category and form the Hilbert space $\ell^2 \fC$, with canonical orthonormal basis given by $\delta_x(y) = 1$ if $x=y$ and $\delta_x(y) = 0$ if $x \neq y$. For each $c \in \fC$, the assignment $\delta_x \mapsto \delta_{cx}$ if $\mft(x) = \mfd(c)$ and $\delta_x \mapsto 0$ if $\mft(x) \neq \mfd(c)$ extends to a bounded linear operator on $\ell^2 \fC$ which we denote by $\lambda_c$. Note that it is at this point, i.e., to ensure boundedness, that we need left cancellation, which actually implies that $\lambda_c$ is a partial isometry. The left regular representation of $\fC$ is given by $\fC \to {\rm PIsom}(\ell^2 \fC), \, c \mapsto \lambda_c$, where ${\rm PIsom}$ stands for the set of partial isometries.

\bdefin
The left reduced C*-algebra of $\fC$ is given by $C^*_{\lambda}(\fC) \defeq C^*(\menge{\lambda_c}{c \in \fC}) \subseteq \cL(\ell^2 \fC)$.
\edefin
\setlength{\parindent}{0cm} \setlength{\parskip}{0.5cm}

\subsection{Inverse semigroup actions and groupoid models}
\label{ss:InvSgp-GPD}

Let us now describe (candidates for) groupoid models for $C^*_{\lambda}(\fC)$. First of all, every $c \in \fC$ induces the partial bijection $\mfd(c) \fC \isom c \fC, \, x \ma cx$. For brevity, we denote this partial bijection by $c$ again. 
\bdefin
The left inverse hull $I_l$ of $\fC$ is the smallest inverse semigroup containing the partial bijections $\menge{c}{c \in \fC}$, i.e., the smallest semigroup of partial bijections of $\fC$ containing the partial bijections $\menge{c}{c \in \fC}$ and closed under inverses.
\edefin
\setlength{\parindent}{0cm} \setlength{\parskip}{0cm}

For more details on inverse semigroups, we refer the reader to \cite[\S~5.5.1]{CELY}. For $s \in I_l$, we denote its domain by $\dom(s)$ and its image by $\im(s)$. Following \cite[\S~5.5.1]{CELY}, in case $I_l$ contains the partial bijection $0$ which is nowhere defined, $\emptyset \isom \emptyset$, we say that $I_l$ contains zero, and we view $I_l$ as an inverse semigroup with zero. A typical non-zero element $s \in I_l$ is of the form $s = d_n^{-1} c_n \dotso d_1^{-1} c_1$ for some $d_i, c_i \in \fC$ with $\mft(c_i) = \mft(d_i)$ and $\mfd(d_i) = \mfd(c_{i+1})$.
\bremark
Elements of $I_l$ are called zigzags in \cite{Sp20}.
\eremark

\bdefin
For $0 \neq s \in I_l$, define $\mfd(s)$ as the unique $\mfv \in \fC^0$ such that $\dom(s) \subseteq \mfv \fC$, and define $\mft(s)$ as the unique $\mfw \in \fC^0$ such that $\im(s) \subseteq \mfw \fC$.
\edefin
Such $\mfv$ and $\mfw$ exist because, if $s = d_n^{-1} c_n \dotso d_1^{-1} c_1$, then $\dom(s) \subseteq \dom(c_1) \subseteq \mfd(c_1) \fC$ and $\im(s) \subseteq \im(d_n^{-1}) \subseteq \mfd(d_n) \fC$.

\bdefin
The semilattice of idempotents of $I_l$ is denoted by $\cJ \defeq \menge{s^{-1}s}{s \in I_l} = \menge{ss^{-1}}{s \in I_l}$.
\edefin
$I_l$ contains $0$ if and only if $\cJ$ contains $\emptyset$. In that case we denote $\emptyset \in \cJ$ by $0$ again. 
\setlength{\parindent}{0cm} \setlength{\parskip}{0.5cm}

Alternatively, we could set $\cJ = \menge{\dom(s)}{s \in I_l} = \menge{\im(s)}{s \in I_l}$. $\cJ$ is the analogue of the set of constructible right ideals in the semigroup context (see \cite{Li12}). Multiplication in $\cJ$ (denoted by $ef$ for $e, f \in \cJ$) corresponds to intersection of subsets of $\fC$, and the partial order \an{$\leq$} on $\cJ$ corresponds to inclusion of subsets.

At this point, we present a variation of $I_l$, following \cite{Sp20}. 
\bdefin
Let $\bar{\cJ}$ denote the set of subsets of $\fC$ of the form $e \setminus \bigcup_{i=1}^n f_n$ for some $e, f_1, \dotsc, f_n \in \cJ$ with $f_1, \dotsc, f_n \leq e$.
\setlength{\parindent}{0.5cm} \setlength{\parskip}{0cm}

Let $\bIl$ be the set of all partial bijections of $\fC$ of the form $s \varepsilon$ for $s \in I_l$ and $\varepsilon \in \bcJ$ with $\varepsilon \leq s^{-1}s$.
\edefin
\setlength{\parindent}{0cm} \setlength{\parskip}{0cm}

It is easy to see that $\bIl$ is again an inverse semigroup, whose semilattice of idempotents is given by $\bcJ$.

\bdefin
The space of characters $\widehat{\cJ}$ is given by the set of non-zero multiplicative maps $\cJ \to \gekl{0,1}$, which send $0 \in \cJ$ to $0 \in \gekl{0,1}$ in case $I_l$ contains $0$. Here multiplication in $\gekl{0,1}$ is the usual one induced by multiplication in $\Rz$. The topology on $\widehat{\cJ}$ is given by point-wise convergence.
\edefin
A basis of compact open sets for the topology of $\widehat{\cJ}$ is given by sets of the form
$$
 \widehat{\cJ}(e;\mff) \defeq \big \lbrace \chi \in \widehat{\cJ}: \: \chi(e) = 1, \, \chi(f) = 0 \ \ \forall \, f \in \mff \big \rbrace,
$$
where $e \in \cJ$ and $\mff \subseteq \cJ$ is a finite subset. By replacing $\mff$ by $\menge{ef}{f \in \mff}$, we can always arrange that $f \leq e$ for all $f \in \mff$. We will also set $\widehat{\cJ}(e) \defeq \lbrace \chi \in \widehat{\cJ}: \: \chi(e) = 1 \rbrace$. Since $\mfv \fC \cap \mfw \fC = \emptyset$ if $\mfv \neq \mfw$, for every $\chi \in \widehat{\cJ}$ there exists a unique $\mfv \in \fC^0$ with $\chi(\mfv \fC^0) = 1$. In other words, we have
$
 \widehat{\cJ} = \coprod_{\mfv \in \fC^0} \widehat{\cJ}(\mfv)
$.
As explained in \cite[\S~5.5.1]{CELY}, there is a one-to-one correspondence between elements in $\widehat{\cJ}$ and filters (on $\cJ$), i.e., non-empty subsets $\cF$ of $\cJ$ with the properties that $0 \notin \cF$ if $I_l$ contains $0$, whenever $e, f \in \cJ$ satisfy $e \leq f$, then $e \in \cF$ implies $f \in \cF$, and whenever $e, f \in \cJ$ lie in $\cF$, then $ef$ must lie in $\cF$ as well. To be concrete, the one-to-one correspondence is implemented by $\widehat{\cJ} \ni \chi \mapsto \chi^{-1}(1) \subseteq \cJ$.
\setlength{\parindent}{0cm} \setlength{\parskip}{0.5cm}

Following \cite[\S~5.6.7]{CELY} and \cite{Sp20}, we now construct a subspace of $\widehat{\cJ}$ which takes into account that elements of $\cJ$ are subsets of $\fC$. First, let $D_{\lambda}(\fC) \defeq \clspan(\menge{1_e}{e \in \cJ}) \subseteq \ell^{\infty}(\fC)$. Here $1_e$ denotes the characteristic function of $e \subseteq \fC$. As explained in \cite[Corollary~5.6.28]{CELY}, the spectrum of $D_{\lambda}(\fC)$ can be identified with the following subspace of $\widehat{\cJ}$:
\bdefin
\label{def:Omega}
Let $\Omega$ be the subspace of $\widehat{\cJ}$ consisting of characters $\chi$ with the property that whenever $e, f_1, \dotsc, f_n \in \cJ$ satisfy $e = \bigcup_{i=1}^n f_i$ as subsets of $\fC$, then $\chi(e) = 1$ implies that $\chi(f_i) = 1$ for some $1 \leq i \leq n$.
\edefin
\setlength{\parindent}{0cm} \setlength{\parskip}{0cm}

\bremark
Following \cite[Corollary~5.6.28]{CELY}, we will view every $\chi \in \Omega$ as a character on $D_{\lambda}(\fC)$, again denoted by $\chi$. Given $\varepsilon = e \setminus \bigcup_{i=1}^n f_n \in \bcJ$, we have $1_{\varepsilon} \in D_{\lambda}(\fC)$, and we set $\chi(\varepsilon) \defeq \chi(1_{\varepsilon})$.
\eremark

\bex
\label{ex:chix}
Given $x \in \fC$, define $\chi_x(e) \defeq 1$ if $x \fC \leq e$ and $\chi_x(e) \defeq 0$ if $x \fC \not\leq e$. It is easy to see that $\chi_x \in \Omega$.
\eex

The following is immediate from the definition of the topology of $\Omega$, using the basis of compact open sets as defined above.
\blemma
\label{lem:chix-dense}
$\menge{\chi_x}{x \in \fC}$ is a dense subset of $\Omega$.
\elemma

The following observation is an immediate consequence of \cite[Corollary~5.6.29]{CELY}.
\blemma
We have $\Omega = \widehat{\cJ}$ if and only if whenever $e, f_1, \dotsc, f_n \in \cJ$ satisfy $e = \bigcup_{i=1}^n f_i$ as subsets of $\fC$, then there exists $1 \leq i \leq n$ with $e = f_i$.
\elemma
\setlength{\parindent}{0cm} \setlength{\parskip}{0.5cm}

Let us now dualize and obtain the following action of $I_l$ on $\widehat{\cJ}$. A given $s \in I_l$ induces the partial homeomorphism $\widehat{\cJ}(s^{-1}s) \isom \widehat{\cJ}(ss^{-1}), \, \chi \ma s.\chi \defeq \chi(s^{-1} \sqcup s)$. These partial homeomorphism give rise to an action $I_l \curvearrowright \widehat{\cJ}$. The same proof as for \cite[Lemma~5.6.40]{CELY} shows that $\Omega$ is $I_l$-invariant, so that we obtain an $I_l$-action $I_l \curvearrowright \Omega$ by restriction. As before, a given $s \in I_l$ acts via the partial homeomorphism $\Omega(s^{-1}s) \isom \Omega(ss^{-1}), \, \chi \ma \chi(s^{-1} \sqcup s)$. Here and in the sequel, given a subspace $X \subseteq \widehat{\cJ}$, we set $X(e) \defeq X \cap \widehat{\cJ}(e)$ and $X(e;\mff) \defeq X \cap \widehat{\cJ}(e;\mff)$. 

We now set out to describe two candidates for a groupoid model for $C^*_{\lambda}(\fC)$. First, we set
$$
 I_l * \Omega \defeq \menge{(s,\chi) \in I_l \times \Omega}{\chi(s^{-1}s) = 1}.
$$
\bdefin
The transformation groupoid $I_l \ltimes \Omega$ is given by $I_l * \Omega / { }_{\sim}$, where we set $(s,\chi) \sim (t,\psi)$ if $\chi = \psi$ and there exists $e \in \cJ$ with $\chi(e) = 1$ and $se = te$. Equivalence classes with respect to $\sim$ are denoted by $[\cdot]$, and for $s \in I_l$ and $U \subseteq \Omega$, we set $[s,U] \defeq \menge{[s,\chi]}{\chi \in U}$. Range and source maps are given by $\rmr([s,\chi]) = s.\chi$ and $\rms([s,\chi]) = \chi$. Multiplication and inversion are defined by $[s,t.\chi][t,\chi] = [st,\chi]$ and $[s,\chi]^{-1} = [s^{-1},s.\chi]$. 
\setlength{\parindent}{0.5cm} \setlength{\parskip}{0cm}

We equip $I_l \ltimes \Omega$ with the unique topology such that for all $s \in I_l$, $[s,\Omega(s^{-1}s)]$ is an open subset of $I_l \ltimes \Omega$ and the source map induces a homeomorphism $[s,\Omega(s^{-1}s)] \isom \Omega(s^{-1}s)$.
\edefin
\setlength{\parindent}{0cm} \setlength{\parskip}{0cm}

As explained in \cite[\S~2.1]{KM}, we call $I_l \ltimes \Omega$ the transformation groupoid and not the groupoid of germs (as in for instance \cite{EP16}) because in other contexts, the groupoid of germs denotes the quotient of a groupoid by the interior of its isotropy subgroupoid (see for instance \cite{Ren08}). 
\setlength{\parindent}{0cm} \setlength{\parskip}{0.5cm}

Now we follow \cite[\S~5]{Sp20} and construct a variation of $I_l \ltimes \Omega$.
\bdefin
We define $I_l \bltimes \Omega \defeq I_l * \Omega / { }_{\bsim}$, where we set $(s,\chi) \bsim (t,\psi)$ if $\chi = \psi$ and there exists $\varepsilon \in \bcJ$ with $\chi(\varepsilon) = 1$ and $s\varepsilon = t\varepsilon$ in $\bIl$. Equivalence classes with respect to $\bsim$ are denoted by $[\cdot]^{\bsim}$. The groupoid structure on $\bIl \ltimes \Omega$ is defined in the same way as for $I_l \ltimes \Omega$.
\setlength{\parindent}{0.5cm} \setlength{\parskip}{0cm}

We equip $I_l \bltimes \Omega$ with the unique topology such that for all $s \in I_l$, $[s,\Omega(s^{-1}s)]^{\bsim}$ is an open subset of $I_l \bltimes \Omega$ and the source map induces a homeomorphism $[s,\Omega(s^{-1}s)]^{\bsim} \isom \Omega(s^{-1}s)$.
\edefin

\bremark
\label{rem:Iluni-acts}
It is straightforward to check that the $I_l$-action on $\Omega$ induces an $\bIl$-action $\bIl \curvearrowright \Omega$ such that the inclusion $I_l \into \bIl$ induces an isomorphism between the transformation groupoid $\bIl \ltimes \Omega$ for $\bIl \curvearrowright \Omega$ and $I_l \bltimes \Omega$ given by $l_l \bltimes \Omega \isom \bIl \ltimes \Omega, \, [s,\chi]^{\bsim} \ma [s,\chi]$.
\eremark
\setlength{\parindent}{0cm} \setlength{\parskip}{0cm}

By construction, we have a canonical projection $I_l \ltimes \Omega \onto I_l \bltimes \Omega$. It is easy to see that this projection induces an isomorphism of the groupoids of germs.
\setlength{\parindent}{0cm} \setlength{\parskip}{0.5cm}

\subsection{Finite alignment}

Let us now introduce a condition which allows us to reduce the discussion from general constructible right ideals to principal right ideals.
\bdefin[{\cite[Definition~3.2]{Sp20}}]
$\fC$ is finitely aligned if for all $a,b \in \fC$, there exists a finite subset $F \subseteq \fC$ such that $a \fC \cap b \fC = \bigcup_{c \in F} c \fC$.
\edefin
\bremark
\label{rem:fin-align--mcm}
The notion of finite alignment is closely related to the notion of minimal common right-multiple (see \cite[Definition~2.38]{Deh15}), which we abbreviate by mcm. Given $a, b, c \in \fC$, $c$ is called an mcm if $c \in a \fC \cap b \fC$ and no proper left divisor $d$ (i.e., an element $d \in \fC$ with $c \in d \fC$) satisfies $d \in a \fC \cap b \fC$. It is immediate from \cite[Lemma~3.3]{Sp20} that $\fC$ is finitely aligned if and only if for all $a, b \in \fC$, the set of mcms $\mcm(a,b)$ is non-empty and finite up to right multiplication by $\fC^*$.
\eremark
\setlength{\parindent}{0cm} \setlength{\parskip}{0cm}

The following observations are immediate from our definitions (see also \cite[\S~3]{Sp20}).
\blemma
\label{lem:FinAlign}
Suppose that $\fC$ is finitely aligned. Then the following hold:
\begin{enumerate}
\item[(i)] For all $e \in \cJ$ there exists a finite subset $F \subseteq \fC$ such that $e = \bigcup_{x \in F} x \fC$, and every $\varepsilon \in \bcJ$ is a finite disjoint union of sets of the form $x \fC \setminus \bigcup_{i=1}^n y_i \fC$ for $x, y_1, \dotsc, y_n \in \fC$.
\item[(ii)] Every $\chi \in \Omega$ is determined by $\cF_{\rm p} \defeq \menge{x \fC \subseteq \fC}{x \in \fC, \, \chi(x \fC) = 1}$, in the sense that for arbitrary $e \in \cJ$, $\chi(e) = 1$ if and only if there exists $x \fC \in \cF_{\rm p}$ with $x \fC \leq e$. Moreover, a basis of compact open sets for $\Omega$ is given by sets of the form $\Omega(x \fC; y_1 \fC, \dotsc, y_n \fC)$. 
\item[(iii)] Every $s \in I_l$ is a finite union of partial bijections of the form $c d^{-1}$, where $d,c \in \fC$ satisfy $\mfd(c) = \mfd(d)$.
\item[(iv)] We have 
\begin{eqnarray*}
 I_l \ltimes \Omega &=& \menge{[cd^{-1},\chi]}{c, d \in \fC, \, \mfd(c) = \mfd(d); \; (cd^{-1},\chi) \in I_l * \Omega},\\
 I_l \bltimes \Omega &=& \menge{[cd^{-1},\chi]^{\bsim}}{c, d \in \fC, \, \mfd(c) = \mfd(d); \; (cd^{-1},\chi) \in I_l * \Omega}.
\end{eqnarray*}
\end{enumerate}
\elemma
In this sense, finite alignment allows us to reduce to principal right ideals.
\setlength{\parindent}{0cm} \setlength{\parskip}{0.5cm}

\subsection{Groupoid models for left regular C*-algebras}
\label{ss:GPD-LeftReg}

Following \cite{Sp20}, we now explain in what sense $I_l \bltimes \Omega$ is a groupoid model for $C^*_{\lambda}(\fC)$. First of all, as explained in \cite[\S~11]{Sp20}, there is a canonical projection $\Lambda: \: C^*_r(I_l \bltimes \Omega) \onto C^*_{\lambda}(\fC)$ given by $\Lambda(1_{[s,\Omega(s^{-1}s)]^{\bsim}})(\delta_x) =  \delta_{s(x)}$ if $x \in \dom(s)$ and $\Lambda(1_{[s,\Omega(s^{-1}s)]^{\bsim}})(\delta_x) = 0$ if $x \notin \dom(s)$. Moreover, it is shown in \cite[\S~11]{Sp20} that $\Lambda$ is an isomorphism if $\fC$ is finitely aligned or $I_l \bltimes \Omega$ is Hausdorff. We present a characterization for the Hausdorff property in Lemma~\ref{lem:HdOmega}. After comparing the groupoids $I_l \ltimes \Omega$ and $I_l \bltimes \Omega$, we obtain similar results for $I_l \ltimes \Omega$. The reader will also find examples for which $\Lambda$ fails to be injective in \cite[\S~11]{Sp20}.

\subsection{The boundary}
\label{ss:boundary}

Finally, we introduce the boundary, following \cite[\S~5.7]{CELY}. 
\bdefin
$\hcJ_{\max}$ denotes the set of characters $\chi \in \hcJ$ for which $\chi^{-1}(1)$ is maximal among all characters $\chi \in \hcJ$.
\edefin
\setlength{\parindent}{0cm} \setlength{\parskip}{0cm}

The same proof as for \cite[Lemma~5.7.7]{CELY} shows that $\hcJ_{\max} \subseteq \Omega$. Hence this justifies the notation $\Omega_{\max} \defeq \hcJ_{\max}$.
The following collects observations about $\Omega_{\max}$, which are proven in the same way as in \cite[\S~5.7]{CELY}.
\blemma
\label{lem:Omegamax}
\begin{enumerate}
\item[(i)] If $I_l$ contains $0$, then $\chi \in \hcJ$ lies in $\Omega_{\max}$ if and only if for all $e \in \cJ$ with $\chi(e) = 0$, there exists $f \in \cJ$ with $\chi(f) = 1$ such that $ef = 0$.
\item[(ii)] For all $0 \neq e \in \cJ$, there exists $\chi \in \Omega_{\max}$ with $\chi(e) = 1$.
\item[(iii)] $\Omegamax$ is $I_l$-invariant.
\end{enumerate}
\elemma

\bdefin
We define the boundary as $\partial \Omega \defeq \overline{\Omegamax} \subseteq \Omega$.
\edefin
By Lemma~\ref{lem:Omegamax}~(iii), $\partial \Omega$ is $I_l$-invariant, so that we may form the boundary groupoids.
\bdefin
We define the boundary groupoids as $I_l \ltimes \bOmega$ and $I_l \bltimes \bOmega$.
\edefin
This also leads to the boundary quotients $C^*_r(I_l \ltimes \bOmega)$ and $C^*_r(I_l \bltimes \bOmega)$.

\bremark
The boundary groupoid $I_l \ltimes \bOmega$ can be identified with the tight groupoid of the left inverse hull $I_l$, in the sense of \cite{Ex08,EP16}. However, an analogous statement does not hold for $I_l \bltimes \bOmega$. Indeed, as noted in \cite[\S~6]{Sp20}, $\widehat{\bcJ}_{\max}$ can be identified with $\Omega$. It follows that $\widehat{\bcJ}_{\max} = \partial \widehat{\bcJ}$, i.e., $\widehat{\bcJ}_{\max}$ itself is already closed. It is also easy to see this directly. This means that the tight groupoid of the inverse semigroup $\bIl$ is given by $I_l \bltimes \Omega$. Thus $I_l \bltimes \bOmega$ does not have an obvious description as a tight groupoid attached to an inverse semigroup.
\eremark
\setlength{\parindent}{0cm} \setlength{\parskip}{0.5cm}

\section{Comparison of groupoid models}

Let us address the natural question when the groupoids $I_l \ltimes \Omega$ and $I_l \bltimes \Omega$ are isomorphic. By construction, there is a canonical projection $I_l \ltimes \Omega \onto I_l \bltimes \Omega$.

First, we collect a few observations which are immediate consequences of our construction.
\setlength{\parindent}{0cm} \setlength{\parskip}{0cm}

\blemma
\label{lem:quotIxO-->IuxO}
\begin{enumerate}
\item[(i)] The canonical projection $I_l \ltimes \Omega \onto I_l \bltimes \Omega$ is an open quotient map. 
\item[(ii)] The canonical projection $I_l \ltimes \Omega \onto I_l \bltimes \Omega$ maps bisections to bisections.
\item[(iii)] The identity map on $\Omega$ induces a bijection between subsets which are invariant for $I_l \ltimes \Omega$ and subsets which are invariant for $I_l \bltimes \Omega$.
\end{enumerate}
\elemma

\blemma
\label{lem:IxO=IuxO}
The canonical projection $I_l \ltimes \Omega \onto I_l \bltimes \Omega$ is an isomorphism if one of the following holds:
\begin{enumerate}
\item[(i)] $\fC$ is finitely aligned.
\item[(ii)] $I_l \ltimes \Omega$ is Hausdorff.
\end{enumerate}
\elemma
\bproof
Take $(s,\chi), (t,\chi) \in I_l * \Omega$ with $(s,\chi) \bsim (t,\chi)$. Then there exists $\varepsilon \in \bcJ$ with $\chi(\varepsilon) = 1$ and $s \varepsilon = t \varepsilon$.
\setlength{\parindent}{0.5cm} \setlength{\parskip}{0cm}

Suppose that (i) holds. By Lemma~\ref{lem:FinAlign}~(i), we may assume that $\varepsilon = x \fC \setminus \bigcup_{i=1}^n y_i \fC$ for some $x, y_1, \dotsc, y_n \in \fC$. Then $s \varepsilon = t \varepsilon$ implies $s(x) = t(x)$, so that, with $e \defeq x \fC$, $s e = t e$. Moreover, $\chi(\varepsilon) = 1$ implies $\chi(e) = 1$ since $\varepsilon \leq e$. This shows that $(s,\chi) \sim (t,\chi)$.

Now assume that (ii) holds. By Lemma~\ref{lem:chix-dense}, we can find $x_i \in \fC$ with $\lim_i \chi_{x_i} = \chi$. As $\chi(\varepsilon) = 1$, we may assume $\chi_{x_i}(\varepsilon) = 1$, i.e., $x_i \in \varepsilon$. Setting $e_i \defeq x_i \fC$, $s \varepsilon = t \varepsilon$ implies $s e_i = t e_i$, and thus $(s,\chi_{x_i}) \sim (t,\chi_{x_i})$. Because $\lim_i (s,\chi_{x_i}) = (s,\chi)$ and $\lim_i (t,\chi_{x_i}) = (t,\chi)$, and since $I_l \ltimes \Omega$ is Hausdorff, we conclude that $(s,\chi) \sim (t,\chi)$.
\eproof
\setlength{\parindent}{0cm} \setlength{\parskip}{0cm}

For a characterization of the Hausdorff property for $I_l \ltimes \Omega$, see Lemma~\ref{lem:HdOmega}.

\bremark
As observed in \S~\ref{ss:InvSgp-GPD}, the canonical projection $I_l \ltimes \Omega \onto I_l \bltimes \Omega$ induces an isomorphism at the level of groupoids of germs. Hence if $I_l \ltimes \Omega$ is effective, the canonical projection $I_l \ltimes \Omega \onto I_l \bltimes \Omega$ must be an isomorphism.
\eremark

The following is an immediate consequence of the results mentioned in \S~\ref{ss:GPD-LeftReg} and Lemma~\ref{lem:IxO=IuxO}
\bcor
If $\fC$ is finitely aligned or $I_l \ltimes \Omega$ is Hausdorff, then $C^*_r(I_l \ltimes \Omega)$ is isomorphic to $C^*_{\lambda}(\fC)$.
\ecor

Let us now compare boundary groupoids.
\setlength{\parindent}{0cm} \setlength{\parskip}{0cm}

\blemma
\label{lem:IxbO=IuxbO}
The canonical projection $I_l \ltimes \bOmega \onto I_l \bltimes \bOmega$ is an isomorphism if one of the following holds:
\begin{enumerate}
\item[(i)] The canonical projection $I_l \ltimes \Omega \onto I_l \bltimes \Omega$ is an isomorphism.
\item[(ii)] $I_l \ltimes \bOmega$ is Hausdorff.
\item[(iii)] $\partial \Omega = \Omegamax$.
\end{enumerate}
\elemma
\bproof
It is easy to see that (i) is a sufficient condition. Now take $(s,\chi), (t,\chi) \in I_l * \Omega$ with $(s,\chi) \bsim (t,\chi)$. Then there exists $\varepsilon \in \bcJ$ with $\chi(\varepsilon) = 1$ and $s \varepsilon = t \varepsilon$, where $\varepsilon = e \setminus \bigcup_{i=1} f_i$ for $e, f_1, \dotsc, f_n \in \cJ$. We first show that if $\chi \in \Omegamax$, then $(s,\chi) \sim (t,\chi)$: Indeed, $\chi(\varepsilon) = 1$ implies that $\chi(f_i) = 0$ for all $1 \leq i \leq n$. By Lemma~\ref{lem:Omegamax}~(i), $\chi(f_i) = 0$ implies that there exists $f'_i \in \cJ$ with $\chi(f'_i) = 1$ and $f_i f'_i = 0$. Set $f' \defeq f'_1 \dotsm f'_n$. Then $\chi(f') = 1$ and $f' f_i = 0$ for all $1 \leq i \leq n$. We conclude that $\chi(e f') = 1$. Moreover, $e f' \subseteq \varepsilon$, so that $s e f' = t e f'$. It follows that $(s,\chi) \sim (t,\chi)$, as desired. This immediately implies that (iii) is a sufficient condition. To treat (ii), assume now that $(s,\chi) \bsim (t,\chi)$ for some $\chi \in \bOmega$. Then there exist $\chi_i \in \Omegamax$ with $\lim_i \chi_i = \chi$. We may assume $\chi_i(\varepsilon) = 1$ since $\chi(\varepsilon) = 1$. It follows that $(s,\chi_i) \bsim (t,\chi_i)$, and, by what we just proved, $(s,\chi_i) \sim (t,\chi_i)$. Since $I_l \ltimes \bOmega$ is Hausdorff, we conclude $\lim_i (s,\chi_i) = (s,\chi) \sim (t,\chi) = \lim_i (t,\chi_i)$, as desired.
\eproof

\bquestion
Do we always have isomorphisms $I_l \ltimes \Omega \onto I_l \bltimes \Omega$ and $I_l \ltimes \bOmega \onto I_l \bltimes \bOmega$? Most likely the answer will be negative, in which case it would be interesting to find concrete examples where the canonical projections fail to be injective.
\equestion
\setlength{\parindent}{0cm} \setlength{\parskip}{0.5cm}

\section{Properties of the groupoids}

We characterize when $I_l \ltimes \Omega$ and $I_l \bltimes \Omega$ are Hausdorff, when $I_l \ltimes \Omega$ is topologically free, and when $I_l \bltimes \Omega$ is effective. These properties have consequences for the reduced C*-algebras of $I_l \ltimes \Omega$ and $I_l \bltimes \Omega$ (see Corollary~\ref{cor:Cons.CG}).

Let us start with the Hausdorff property. The following will be an application of \cite[Theorem~3.15]{EP16}.
\blemma
\label{lem:HdOmega}
\begin{enumerate}
\item[(i)] $I_l \ltimes \Omega$ is Hausdorff if and only if for all $s \in I_l$, there exists a (possibly empty) finite subset $\gekl{e_1, \dotsc, e_n} \subseteq \cJ$ with 
$
 \menge{x \in \dom(s)}{s(x) = x} = \bigcup_{i=1}^n e_i
$.
\item[(ii)] $I_l \bltimes \Omega$ is Hausdorff if and only if for all $s \in I_l$, there exists a (possibly empty) finite subset $\gekl{\varepsilon_1, \dotsc, \varepsilon_n} \subseteq \bcJ$ with 
$
 \menge{x \in \dom(s)}{s(x) = x} = \bigcup_{i=1}^n \varepsilon_i
$.
\end{enumerate}
\elemma
\setlength{\parindent}{0cm} \setlength{\parskip}{0cm}

\bproof
(i) \cite[Theorem~3.15]{EP16} implies that $I_l \ltimes \Omega$ is Hausdorff if and only if for all $s \in I_l$, the subset 
\begin{equation}
\label{e:chie=1,se=e}
 \menge{\chi \in \Omega}{\exists \, e \in \cJ \text{ with } se = e \text{ and } \chi(e) = 1}
\end{equation}
is closed in $\menge{\chi \in \Omega}{\chi(s^{-1}s) = 1}$. The latter statement is equivalent to compactness of the set in \eqref{e:chie=1,se=e} because $\menge{\chi \in \Omega}{\chi(s^{-1}s) = 1}$ is compact. This in turn is true if and only if there exists a finite subset $\gekl{e_1, \dotsc, e_n} \subseteq \cJ$ with $se_i = e_i$ for all $1 \leq i \leq n$ and
\begin{equation}
\label{e:chie=1...=Omega(ei)}
 \menge{\chi \in \Omega}{\exists \, e \in \cJ \text{ with } se = e \text{ and } \chi(e) = 1} = \bigcup_{i=1}^n \Omega(e_i).
\end{equation}
We claim that \eqref{e:chie=1...=Omega(ei)} is equivalent to $\menge{x \in \dom(s)}{s(x) = x} = \bigcup_{i=1}^n e_i$. 

As $s e_i = e_i$, we always have $\menge{x \in \dom(s)}{s(x) = x} \supseteq \bigcup_{i=1}^n e_i$. Assume that $\menge{x \in \dom(s)}{s(x) = x} \subseteq \bigcup_{i=1}^n e_i$. Given $\chi \in \Omega$ together with $e \in \cJ$ such that $se = e$ and $\chi(e) = 1$, we must have $e \subseteq \bigcup_{i=1}^n e_i$. As $\chi$ lies in $\Omega$, $\chi(e) = 1$ implies that there exists $1 \leq i \leq n$ with $\chi(e_i) = 1$. Hence \eqref{e:chie=1...=Omega(ei)} holds. Conversely, suppose that \eqref{e:chie=1...=Omega(ei)} holds. Take $x \in \dom(s)$ with $s(x) = x$. Then $\chi_x$ lies in the set on the left-hand side of \eqref{e:chie=1...=Omega(ei)}, hence there exists $1 \leq i \leq n$ with $\chi_x(e_i) = 1$. The latter implies that $x \in e_i$. This shows $\menge{x \in \dom(s)}{s(x) = x} \subseteq \bigcup_{i=1}^n e_i$, as desired.
\setlength{\parindent}{0.5cm} \setlength{\parskip}{0cm}

(ii) \cite[Theorem~3.15]{EP16} implies that $I_l \bltimes \Omega$ is Hausdorff if and only if for all $t \in \bIl$, the subset 
\begin{equation}
\label{e:chieps=1,teps=eps}
 \menge{\chi \in \Omega}{\exists \, \varepsilon \in \bcJ \text{ with } t \varepsilon = \varepsilon \text{ and } \chi(\varepsilon) = 1}
\end{equation}
is closed in $\menge{\chi \in \Omega}{\chi(t^{-1}t) = 1}$. First, we claim that the latter is equivalent to the statement that for all $s \in I_l$, the subset 
\begin{equation}
\label{e:chieps=1,seps=eps}
 \menge{\chi \in \Omega}{\exists \, \varepsilon \in \bcJ \text{ with } s \varepsilon = \varepsilon \text{ and } \chi(\varepsilon) = 1}
\end{equation}
is closed in $\menge{\chi \in \Omega}{\chi(s^{-1}s) = 1}$. Indeed, a general element $t \in \bIl$ is of the form $s \delta$ for some $\delta \in \bcJ$ with $\delta \leq s^{-1}s$. Now it is straightforward to see that the set in \eqref{e:chieps=1,teps=eps} coincides with the intersection of the set in \eqref{e:chieps=1,seps=eps} and $\Omega(\delta)$. If the set in \eqref{e:chieps=1,seps=eps} is closed in $\menge{\chi \in \Omega}{\chi(s^{-1}s) = 1}$, then its intersection with $\Omega(\delta)$ must be closed in $\menge{\chi \in \Omega}{\chi(s^{-1}s) = 1} \cap \Omega(\delta) = \menge{\chi \in \Omega}{\chi(t^{-1}t) = 1}$. This shows our claim. Now the rest of the proof is similar as for (i).
\eproof
\setlength{\parindent}{0cm} \setlength{\parskip}{0cm}

In combination with Lemma~\ref{lem:FinAlign}, the following is immediate.
\bcor
\label{cor:FinAlign-HdOmega}
Assume that $\fC$ is finitely aligned. Then $I_l \ltimes \Omega \cong I_l \bltimes \Omega$ is Hausdorff if and only if for all $c, d \in \fC$ with $\mfd(c) = \mfd(d)$ and $\mft(c) = \mft(d)$, there exists a finite subset $\gekl{x_1, \dotsc, x_n} \subseteq \fC$ with 
$
 \menge{x \in \fC}{cx = dx} = \bigcup_{i=1}^n x_i \fC
$.
\ecor

\bremark
Lemma~\ref{lem:HdOmega} and Corollary~\ref{cor:FinAlign-HdOmega} explain the results in \cite[\S~7]{Sp20} that $I_l \bltimes \Omega$ is Hausdorff if $\fC$ is finitely aligned and right cancellative, or if $\fC$ embeds into a groupoid. In the first case, the set $\menge{x \in \fC}{cx = dx}$ is either empty or we have $c = d$, which implies that $\menge{x \in \fC}{cx = dx} = \fC$. In the second case, the set $\menge{x \in \dom(s)}{s(x) = x}$ is either empty or we have $s \in \bcJ$, in which case $\menge{x \in \dom(s)}{s(x) = x}$ coincides with $s^{-1}s$, where we view the latter as a subset of $\fC$.
\eremark
\setlength{\parindent}{0cm} \setlength{\parskip}{0.5cm}

Let us now consider topological freeness and effectiveness. Recall that an {\'e}tale groupoid $\cG$ is called effective if the interior of its isotropy subgroupoid coincides with the unit space, i.e., ${\rm Iso}(\cG)^{\circ} = \cG^{(0)}$. Following \cite[Definition~2.20]{KM}, we call an {\'e}tale groupoid $\cG$ topologically free if for every open bisection $\gamma$ with $\gamma \subseteq \cG \setminus \cG^{(0)}$, 
$
 \menge{x \in \cG^{(0)}}{\cG_x^x \cap \gamma \neq \emptyset}
$
has empty interior, or equivalently, 
$
 \menge{x \in \rms(\gamma)}{\gamma x \notin \cG_x^x}
$
is dense in $\rms(\gamma)$. By \cite[Lemma~2.23]{KM}, $\cG$ is topologically free if $\cG$ is effective, and the converse holds if $\cG$ is Hausdorff. Topological freeness for groupoids is of interest because it implies the intersection properties for essential groupoid C*-algebras (see \cite[\S~7.5]{KM} for more information).

Now we set $\fC^{*,0} \defeq \menge{u \in \fC^*}{\mft(u) = \mfd(u)}$, and set $\fC^{*,0} \ltimes \Omega \defeq \menge{[u,\chi] \in I_l \ltimes \Omega}{u \in \fC^{*,0}}$.
\btheo
\label{thm:IxOmega-topfree}
The following are equivalent:
\setlength{\parindent}{0cm} \setlength{\parskip}{0cm}

\begin{enumerate}
\item[(i)] $I_l \ltimes \Omega$ is topologically free;
\item[(ii)] $\fC^{*,0} \ltimes \Omega$ is topologically free;
\item[(iii)] For all $\mfv \in \fC^0$, $u \in \mfv \fC^* \mfv$, $f_1, \dotsc, f_n \in \cJ$ with $f_i \lneq \mfv \fC$ for all $1 \leq i \leq n$, $u z \in z \fC^*$ for all $z \in \mfv \fC \setminus \bigcup_{i=1}^n f_i$ implies that there exists $x \in \mfv \fC \setminus \bigcup_{i=1}^n f_i$ with $ux = x$.
\end{enumerate}
\etheo
\setlength{\parindent}{0cm} \setlength{\parskip}{0cm}

\bproof
(i) $\Rarr$ (ii): $\fC^{*,0} \ltimes \Omega$ is an open subgroupoid of $I_l \ltimes \Omega$. Thus an open bisection $\gamma$ of $\fC^{*,0} \ltimes \Omega$ with $\gamma \subseteq (\fC^{*,0} \ltimes \Omega) \setminus \Omega$ is also an open bisection of $I_l \ltimes \Omega$ contained in $(I_l \ltimes \Omega) \setminus \Omega$. Moreover, $\gamma x \notin (I_l \ltimes \Omega)_x^x$ implies that $\gamma x \notin (\fC^{*,0} \ltimes \Omega)_x^x$. This shows that 
$$
 \menge{x \in \rms(\gamma)}{\gamma x \notin (I_l \ltimes \Omega)_x^x} \subseteq \menge{x \in \rms(\gamma)}{\gamma x \notin (\fC^{*,0} \ltimes \Omega)_x^x}.
$$
Hence $\fC^{*,0} \ltimes \Omega$ is topologically free if $I_l \ltimes \Omega$ is topologically free.
\setlength{\parindent}{0.5cm} \setlength{\parskip}{0cm}

(ii) $\Rarr$ (i): Assume that $I_l \ltimes \Omega$ is not topologically free. Then we can find $s \in I_l$ and an open set $U \subseteq \Omega(s^{-1}s)$ with $[s,U] \subseteq (I_l \ltimes \Omega) \setminus \Omega$ and $[s,U] \subseteq {\rm Iso}(I_l \ltimes \Omega)$. As $\menge{\chi_x}{x \in \fC}$ is dense in $\Omega$, there exists $x \in \fC$ with $\chi_x \in U$. $s.\chi_x = \chi_x$ implies that $s(x) = xu$ for some $u \in \fC^{*,0}$. As $[s,\chi_x] \neq \chi_x$, we conclude that $u \notin \fC^0$. Set $V \defeq \Omega(x \fC) \cap U$. $V$ is not empty, so that $x^{-1}.V \neq \emptyset$. It is easy to see that $[x,\Omega(\mfd(x))]^{-1} [s,V] [x,\Omega(\mfd(x))] = [u, x^{-1}.V]$. Moreover, $[x,\Omega(\mfd(x))]^{-1} [s,V] [x,\Omega(\mfd(x))]$ is contained in ${\rm Iso}(I_l \ltimes \Omega) \setminus \Omega$ because $[s,V] \subseteq {\rm Iso}(I_l \ltimes \Omega) \setminus \Omega$. This means that $\fC^{*,0} \ltimes \Omega$ is not topologically free.

(ii) $\Rarr$ (iii): Assume $u z \in z \fC^*$ for all $z \in \mfv \fC \setminus \bigcup_{i=1}^n f_i$. Set $U \defeq \Omega(\mfv \fC; f_1, \dotsc, f_n)$. Then $[u,U] \subseteq {\rm Iso}(\fC^{*,0} \ltimes \Omega)$. As $\fC^{*,0} \ltimes \Omega$ is topologically free, there exists $\chi \in U$ with $[u,\chi] = \chi$, i.e., there exists $e \in \cJ$ with $\chi(e) = 1$ and $ue = e$. $\chi(\mfv \fC \setminus \bigcup_{i=1}^n f_i) = 1$ implies that $e \not\subseteq \bigcup_{i=1}^n f_i$. Hence we can choose $x \in e \setminus \bigcup_{i=1}^n f_i$, and we have $ux = x$.

(iii) $\Rarr$ (ii): First we claim that (iii) is equivalent to the following stronger statement:
\begin{enumerate}
\item[(iii')] For all $\mfv \in \fC^0$, $u \in \mfv \fC^* \mfv$, $e, f_1, \dotsc, f_n \in \cJ$ with $e, f_1, \dotsc, f_n \leq \mfv \fC$ and $\bigcup_{i=1}^n f_i \subsetneq e$, $u z \in z \fC^*$ for all $z \in e \setminus \bigcup_{i=1}^n f_i$ implies that there exists $x \in e \setminus \bigcup_{i=1}^n f_i$ with $ux = x$.
\end{enumerate}
Indeed, to prove (iii) $\Rarr$ (iii'), take $y \in e \setminus \bigcup_{i=1}^n f_i$ and set $\mfv \defeq \mfd(y)$. By assumption, $uy \in y \fC^*$, and hence we have $uy = y \ti{u}$ for some $\ti{u} \in \mfv \fC^* \mfv$. Set $f'_i \defeq y \fC \cap f_i$. Then $\bigcup_{i=1}^n f'_i \subsetneq y \fC$ implies that $\bigcup_{i=1}^n y^{-1} f'_i \subsetneq \mfv \fC$. For every $\ti{x} \in \mfv \fC \setminus \bigcup_{i=1}^n y^{-1} f'_i$, we have by assumption $y \ti{u} \ti{x} = u y \ti{x} \in y \ti{x} \fC^*$ and thus $\ti{u} \ti{x} \in \ti{x} \fC^*$. Hence (iii) implies that there exists $x \in \mfv \fC \setminus \bigcup_{i=1}^n y^{-1} f'_i$ with $\ti{u} x = x$. Then $y x \in e \setminus \bigcup_{i=1}^n f_i$ and $u y x = y \ti{u} x = y x$, as desired.

Now assume that (iii') holds. Let $u \in \fC^{*,0}$, $U = \Omega(e; f_1, \dotsc, f_n)$, and assume that $[u,U] \subseteq {\rm Iso}(\fC^{*,0} \ltimes \Omega)$. Then we must have $u z \in z \fC^*$ for all $z \in e \setminus \bigcup_{i=1}^n f_i$. Hence (iii') implies that there exists $x \in e \setminus \bigcup_{i=1}^n f_i$ with $ux = x$. Then $\chi_x \in U$ because $x \in e \setminus \bigcup_{i=1}^n f_i$. Moreover, $ux = x$ implies that $[u,\chi_x] = \chi_x \in \Omega$. Hence $\fC^{*,0} \ltimes \Omega$ is topologically free.
\eproof
\setlength{\parindent}{0cm} \setlength{\parskip}{0cm}

We now consider $I_l \bltimes \Omega$. As before, we set $\fC^{*,0} \bltimes \Omega \defeq \menge{[u,\chi]^{\bsim} \in I_l \bltimes \Omega}{u \in \fC^{*,0}}$.
\btheo
\label{thm:IunixOmega-eff}
The following are equivalent:
\begin{enumerate}
\item[(i)] $I_l \bltimes \Omega$ is effective;
\item[(ii)] $\fC^{*,0} \bltimes \Omega$ is effective;
\item[(iii)] For all $\mfv \in \fC^0$, $u \in \mfv \fC^* \mfv$, $f_1, \dotsc, f_n \in \cJ$ with $f_i \lneq \mfv \fC$ for all $1 \leq i \leq n$, $u z \in z \fC^*$ for all $z \in \mfv \fC \setminus \bigcup_{i=1}^n f_i$ implies that $ux = x$ for all $x \in \mfv \fC \setminus \bigcup_{i=1}^n f_i$.
\end{enumerate}
\etheo
\bproof
(i) $\Rarr$ (ii) is clear because $\fC^{*,0} \bltimes \Omega$ is an open subgroupoid of $I_l \bltimes \Omega$.
\setlength{\parindent}{0.5cm} \setlength{\parskip}{0cm}

(ii) $\Rarr$ (i): Suppose that $I_l \bltimes \Omega$ is not effective. Then there exist $s \in I_l$, $U \defeq \Omega(e; f_1, \dotsc, f_n)$ and $\chi \in U$ with $[s,U] \subseteq {\rm Iso}(I_l \bltimes \Omega)$ and $[s,\chi] \neq \chi$. Set $\varepsilon \defeq e \setminus \bigcup_{i=1}^n f_i$. $[s,\chi] \neq \chi$ implies that $s \varepsilon \neq \varepsilon$, i.e., there exists $x \in \varepsilon$ with $s(x) \neq x$. We have $\chi_x \in U$, and $s(x) \neq x$ implies $[s,\chi_x] \neq \chi_x$. However, $s.\chi_x = \chi_x$, and thus $s(x) = xu$ for some $u \in \fC^{*,0}$ with $\mft(u) = \mfd(u) = \mfd(x)$. We deduce $x \neq xu$, i.e., $u \neq \mfd(x)$. Set $V \defeq \Omega(x) \cap U$. Then $\chi_x \in V$, so $V$ is not empty. Moreover, $[u,c^{-1}.V] = [x,\Omega(\mfd(x))]^{-1} [s,V] [x,\Omega(\mfd(x))]$ is contained in ${\rm Iso}(\fC^{*,0} \bltimes \Omega)$. We have $[u,\chi_{\mfd(x)}] = [x,\Omega(\mfd(x))]^{-1} [s,\chi_x] [x,\Omega(\mfd(x))] \in [u,c^{-1}.V]$ and $[u,\chi_{\mfd(x)}] \neq \chi_{\mfd(x)}$ because $u \mfd(x) \neq \mfd(x)$. It follows that $\fC^{*,0} \bltimes \Omega$ is not effective.

To prove (ii) $\LRarr$ (iii), we first show that (iii) is equivalent to the following stronger statement:
\begin{enumerate}
\item[(iii')] For all $\mfv \in \fC^0$, $u \in \mfv \fC^* \mfv$, $e, f_1, \dotsc, f_n \in \cJ$ with $e, f_1, \dotsc, f_n \leq \mfv \fC$ and $\bigcup_{i=1}^n f_i \subsetneq e$, $u z \in z \fC^*$ for all $z \in e \setminus \bigcup_{i=1}^n f_i$ implies that $ux = x$ for all $x \in e \setminus \bigcup_{i=1}^n f_i$.
\end{enumerate}
Indeed, to prove (iii) $\Rarr$ (iii'), take $x \in e \setminus \bigcup_{i=1}^n f_i$ and set $\mfv \defeq \mfd(x)$. By assumption, $ux \in x \fC^*$, and hence we have $ux = x \ti{u}$ for some $\ti{u} \in \mfv \fC^* \mfv$. Set $f'_i \defeq x \fC \cap f_i$. Then $\bigcup_{i=1}^n f'_i \subsetneq x \fC$ implies that $\bigcup_{i=1}^n x^{-1} f'_i \subsetneq \mfv \fC$. For every $\ti{x} \in \mfv \fC \setminus \bigcup_{i=1}^n x^{-1} f'_i$, we have by assumption $x \ti{u} \ti{x} = u x \ti{x} \in x \ti{x} \fC^*$ and thus $\ti{u} \ti{x} \in \ti{x} \fC^*$. Hence (iii) implies $\ti{u} = \ti{u} \mfv = \mfv$ and thus $u x = x \ti{u} = x$. As $x$ was an arbitrary element of $e \setminus \bigcup_{i=1}^n f_i$, we are done.

Now let us prove (ii) $\LRarr$ (iii). $\fC^{*,0} \bltimes \Omega$ is effective if and only if for all $u \in \fC^{*,0}$ and $\varepsilon = e \setminus \bigcup_{i=1}^n f_i \in \bcJ$, $[u,\Omega(\varepsilon)] \subseteq {\rm Iso}(\fC^{*,0} \bltimes \Omega)$ implies $[u,\Omega(\varepsilon)] = \Omega(\varepsilon)$. $[u,\Omega(\varepsilon)] \subseteq {\rm Iso}(\fC^{*,0} \bltimes \Omega)$ holds if and only if $uz \in z \fC^*$ for all $z \in \varepsilon$, whereas $[u,\Omega(\varepsilon)] = \Omega(\varepsilon)$ holds if and only if $u \varepsilon = \varepsilon$, i.e., $ux = x$ for all $x \in \varepsilon$. We conclude that (ii) and (iii') are equivalent.
\eproof
\setlength{\parindent}{0cm} \setlength{\parskip}{0cm}

The following are immediate consequences.

\bcor
If $I_l \bltimes \Omega$ is effective, then $I_l \ltimes \Omega$ is topologically free.
\ecor

\bcor
\label{cor:FinAlign-tf,eff}
Assume that $\fC$ is finitely aligned. 
\setlength{\parindent}{0.5cm} \setlength{\parskip}{0cm}

\begin{enumerate}
\item[(i)] $I_l \ltimes \Omega \cong I_l \bltimes \Omega$ is topologically free if and only if for all $\mfv \in \fC^0$, $u \in \mfv \fC^* \mfv$, $c_1, \dotsc, c_n \in \mfv \fC \setminus \mfv \fC^*$, $u z \in z \fC^*$ for all $z \in \mfv \fC \setminus \bigcup_{i=1}^n c_i \fC$ implies that there exists $x \in \mfv \fC \setminus \bigcup_{i=1}^n c_i \fC$ with $ux = x$.
\item[(ii)] $I_l \ltimes \Omega \cong I_l \bltimes \Omega$ is effective if and only if for all $\mfv \in \fC^0$, $u \in \mfv \fC^* \mfv$, $c_1, \dotsc, c_n \in \mfv \fC \setminus \mfv \fC^*$, $u z \in z \fC^*$ for all $z \in \mfv \fC \setminus \bigcup_{i=1}^n c_i \fC$ implies that $ux = x$ for all $x \in \mfv \fC \setminus \bigcup_{i=1}^n c_i \fC$.
\end{enumerate}
\ecor
\setlength{\parindent}{0cm} \setlength{\parskip}{0cm}

We also note the following special case, where our conditions simplify.
\bcor
\label{cor:chiv-open}
Assume that for all $\mfv \in \fC^0$, there exist $f_1, \dotsc, f_n \in \cJ$ with $\mfv \fC \setminus \bigcup_{i=1}^n f_i = \mfv \fC^*$. Then the following are equivalent:
\begin{enumerate}
\item[(i)] $I_l \bltimes \Omega$ is effective.
\item[(ii)] $I_l \ltimes \Omega$ is topologically free.
\item[(iii)] $\fC^{*,0} = \fC^0$.
\end{enumerate}
\ecor
\bproof
(i) $\Rarr$ (ii) has been noted above. Let us prove (ii) $\Rarr$ (iii). We have for all $z \in \mfv \fC \setminus \bigcup_{i=1}^n f_i = \mfv \fC^*$ that $uz = z (z^{-1} u z) \in z \fC^*$. Hence Theorem~\ref{thm:IxOmega-topfree}~(iii) implies that there exists $x \in \mfv \fC \setminus \bigcup_{i=1}^n f_i = \mfv \fC^*$ with $ux = x$. Hence $u = ux x^{-1} = x x^{-1} = \mfv$. (iii) $\Rarr$ (i) is immediate from Theorem~\ref{thm:IunixOmega-eff}.
\eproof

\bremark
Theorems~\ref{thm:IxOmega-topfree} and \ref{thm:IunixOmega-eff} generalize \cite[Theorem~5.9]{LS}.
\eremark

In combination with \cite[Theorem~7.29]{KM} (and the explanations following Theorem~7.29 in \cite{KM}), the following are consequences of our results above.
\bcor
\label{cor:Cons.CG}
If the conditions in Lemma~\ref{lem:HdOmega}~(i) and Theorem~\ref{thm:IxOmega-topfree}~(iii) are satisfied, then $C^*_r(I_l \ltimes \Omega)$ has the intersection property. 
\setlength{\parindent}{0.5cm} \setlength{\parskip}{0cm}

If the conditions in Lemma~\ref{lem:HdOmega}~(ii) and Theorem~\ref{thm:IunixOmega-eff}~(iii) are satisfied, then $C^*_r(I_l \bltimes \Omega)$ has the intersection property.

Suppose that $\fC$ is finitely aligned. If the condition in Corollary~\ref{cor:FinAlign-HdOmega} and one of the conditions in Corollary~\ref{cor:FinAlign-tf,eff} are satisfied, then $C^*_r(I_l \ltimes \Omega) \cong C^*_r(I_l \bltimes \Omega)$ has the intersection property.
\ecor
\setlength{\parindent}{0cm} \setlength{\parskip}{0cm}

\bremark
It is also possible to give a characterization for minimality of $I_l \ltimes \Omega$ and $I_l \bltimes \Omega$ by formulating a characterization when $\Omega = \bOmega$ along the lines of \cite[Lemma~5.7.19]{CELY} and then applying our characterization for minimality of $I_l \ltimes \bOmega$ and $I_l \bltimes \bOmega$ (see Lemma~\ref{lem:bOmega-min}).
\eremark

\bremark
It would also be possible to formulate sufficient criteria for local contractiveness of $I_l \ltimes \Omega$ and $I_l \bltimes \Omega$. However, this happens only in rather special situations (see Proposition~\ref{prop:XminusXLocallyContractive} and Corollary~\ref{cor:InvSub=Ideals,StronglyPurelyInf}, for example). For instance, in the setting of Corollary~\ref{cor:chiv-open}, $I_l \ltimes \Omega$ and $I_l \bltimes \Omega$ and are never locally contractive because the assumptions in Corollary~\ref{cor:chiv-open} imply that $\gekl{\chi_{\mfv}}$ is open for all $\mfv \in \fC^0$.
\eremark
\setlength{\parindent}{0cm} \setlength{\parskip}{0.5cm}

\section{Properties of the boundary groupoid}
\label{s:Boundary}

We characterize when $I_l \ltimes \bOmega$ and $I_l \bltimes \bOmega$ are Hausdorff or minimal, when $I_l \ltimes \bOmega$ is topologically free, when $I_l \bltimes \bOmega$ is effective, and we give a sufficient condition for local contractiveness of $I_l \ltimes \bOmega$ and $I_l \bltimes \bOmega$. These properties have consequences for the boundary quotients (see Corollary~\ref{cor:Cons.CbG}).

Note that if $I_l$ does not contain zero, then $\# \fC^0 = 1$ and $\bOmega$ degenerates to a point. Because of this, it suffices in the following to focus on the case when $I_l$ contains zero.

We first consider the Hausdorff property. The following is an application of \cite[Theorem~3.16]{EP16} because $I_l \ltimes \bOmega$ is the tight groupoid of the inverse semigroup $I_l$.
\blemma
\label{lem:HdIxbOmega}
$I_l \ltimes \bOmega$ is Hausdorff if and only if for all $s \in I_l$ there exist $e_1, \dotsc, e_n \in \cJ$ with $s e_i = e_i$ such that for all $0 \neq e \in \cJ$ with $se = e$, there exists $1 \leq i \leq n$ with $e e_i \neq 0$.
\elemma
\setlength{\parindent}{0cm} \setlength{\parskip}{0cm}

Now we characterize when $I_l \bltimes \bOmega$ is Hausdorff.
\blemma
\label{lem:HdIbxbOmega}
$I_l \bltimes \bOmega$ is Hausdorff if and only if for all $s \in I_l$ there exist $\varepsilon_1, \dotsc, \varepsilon_n \in \bcJ$ with $s \varepsilon_i = \varepsilon_i$ such that for all $0 \neq e \in \cJ$ with $s e = e$, there exists $1 \leq i \leq n$ such that $e \varepsilon_i \neq 0$.
\elemma
\bproof
We make use of the identification $I_l \bltimes \bOmega \cong \bIl \ltimes \bOmega$ (see Remark~\ref{rem:Iluni-acts}). \cite[Theorem~3.15]{EP16}, applied to $\bIl \curvearrowright \bOmega$, implies that $I_l \bltimes \bOmega$ is Hausdorff if and only if for all $s \in \bIl$ there exist $\varepsilon_1, \dotsc, \varepsilon_n \in \bcJ$ with $s \varepsilon_i = \varepsilon_i$ such that for all $\chi \in \bOmega$, $\varepsilon \in \bcJ$ with $\chi(\varepsilon) = 1$ and $s \varepsilon = \varepsilon$, there exists $1 \leq i \leq n$ such that $\chi(\varepsilon_i) = 1$. We may assume that $s \in I_l$ in this statement because every $\bar{s} \in \bIl$ is of the form $s \delta$ for some $s \in I_l$ and $\delta \in \bcJ$, and we can form products of $\varepsilon$ and $\varepsilon_i$ with $\delta$. Next, we claim that the statement is equivalent to the following: For all $s \in I_l$ there exist $\varepsilon_1, \dotsc, \varepsilon_n \in \bcJ$ with $s \varepsilon_i = \varepsilon_i$ such that for all $\chi \in \Omegamax$, $\varepsilon \in \bcJ$ with $\chi(\varepsilon) = 1$ and $s \varepsilon = \varepsilon$, there exists $1 \leq i \leq n$ such that $\chi(\varepsilon_i) = 1$. Indeed, given $\chi \in \bOmega$, we can always find $\eta_{\lambda} \in \Omegamax$ with $\chi = \lim_{\lambda} \eta_{\lambda}$. We may then assume that $\eta_{\lambda}(\varepsilon) = 1$ for all $\lambda$, and then deduce that for all $\lambda$, there exists $1 \leq i \leq n$ with $\eta_{\lambda}(\varepsilon_i) = 1$. By passing to a subnet if necessary, we arrange that there exists $1 \leq i \leq n$ with $\eta_{\lambda}(\varepsilon_i) = 1$ for all $\lambda$, and thus $\chi(\varepsilon_i) = 1$. Now we claim that our new statement is equivalent to the following: For all $s \in I_l$ there exist $\varepsilon_1, \dotsc, \varepsilon_n \in \bcJ$ with $s \varepsilon_i = \varepsilon_i$ such that for all $\chi \in \Omegamax$, $0 \neq e \in \cJ$ with $\chi(e) = 1$ and $s e = e$, there exists $1 \leq i \leq n$ such that $\chi(\varepsilon_i) = 1$. Indeed, given $\chi \in \Omegamax$ and $\varepsilon \in \bcJ$ with $\chi(\varepsilon) = 1$, Lemma~\ref{lem:Omegamax} implies that there exists $e \in \cJ$ with $\chi(e) = 1$ and $e \leq \varepsilon$. Finally, we claim that our statement is equivalent to the desired one: For all $s \in I_l$ there exist $\varepsilon_1, \dotsc, \varepsilon_n \in \bcJ$ with $s \varepsilon_i = \varepsilon_i$ such that for all $0 \neq e \in \cJ$ with $s e = e$, there exists $1 \leq i \leq n$ such that $e \varepsilon_i \neq 0$. To see \an{$\Rarr$}, if there exists $0 \neq e \in \cJ$ with $e \varepsilon_i = 0$ for all $i$, then Lemma~\ref{lem:Omegamax} yields a character $\chi \in \Omegamax$ with $\chi(e) = 1$, and we obtain $\chi(\varepsilon_i) = 0$ for all $i$. For \an{$\Larr$}, assume that there exist $\chi \in \Omegamax$, $0 \neq e \in \cJ$ with $\chi(e) = 1$ and $s e = e$ such that $\chi(\varepsilon_i) = 0$ for all $i$. Write $\varepsilon_i = e_i \setminus \bigcup_{f \in \mff_i} f$. $\chi(\varepsilon_i) = 0$ implies that $\chi(e_i) = 0$ or $\chi(f_i) = 1$ for some $f_i \in \mff_i$. In the first case, Lemma~\ref{lem:Omegamax} yields $e'_i \in \cJ$ with $e_i e'_i = 0$ and $\chi(e'_i) = 1$. In that case set $g_i \defeq e'_i$. In the second case, set $g_i \defeq f_i$. In any case, we obtain $\chi(g_i) = 1$ and $g_i \varepsilon_i = 0$. Now set $e' \defeq e \prod_i g_i$. It follows that $\chi(e') = 1$ (and thus $e' \neq 0$), $s e' = e'$ and $e' \varepsilon_i = 0$ for all $i$.
\eproof

Our characterization simplifies in the finitely aligned case.
\bcor
\label{cor:FinAligned-HdbOmega}
Suppose that $\fC$ is finitely aligned. Then $I_l \ltimes \bOmega \cong I_l \bltimes \bOmega$ is Hausdorff if and only if for all $c, d \in \fC$ with $\mft(d) = \mft(c)$, there exist $x_1, \dotsc, x_n \in \fC$ with $c x_i = d x_i$ for all $1 \leq i \leq n$ such that for all $x \in \fC$ with $cx = dx$, there exists $1 \leq i \leq n$ with $x \fC \cap x_i \fC \neq \emptyset$.
\ecor
\setlength{\parindent}{0cm} \setlength{\parskip}{0.5cm}

Next, we consider minimality.
\blemma
\label{lem:bOmega-min}
The following are equivalent:
\setlength{\parindent}{0cm} \setlength{\parskip}{0cm}

\begin{enumerate}
\item[(i)] $I_l \ltimes \bOmega$ is minimal. 
\item[(ii)] $I_l \bltimes \bOmega$ is minimal. 
\item[(iii)] For all non-zero $e, f \in \cJ$ there exist $s_1, \dotsc, s_n \in I_l$ such that for all $e' \in \cJ$ with $e' \leq e$, there exists $1 \leq i \leq n$ with $e' (s_i f s_i^{-1}) \neq 0$.
\end{enumerate}
\elemma
\bproof
(i) $\LRarr$ (ii) follows from Lemma~\ref{lem:quotIxO-->IuxO}. (i) $\LRarr$ (iii) follows from \cite[Theorem~5.5]{EP16}.
\eproof
\setlength{\parindent}{0cm} \setlength{\parskip}{0cm}

We record the following characterization of minimality in the finitely aligned case.
\bcor
\label{cor:FinAlign_bOmega-min}
Suppose that $\fC$ is finitely aligned. Then $I_l \ltimes \bOmega \cong I_l \bltimes \bOmega$ is minimal if and only if for all $\mfv, \mfw \in \fC^0$ there exist $x_1, \dotsc, x_n \in \mfv \fC$ with $\mfw \fC \mfd(x_i) \neq \emptyset$ for all $1 \leq i \leq n$, such that for all $x \in \mfv \fC$ there exists $1 \leq i \leq n$ with $x \fC \cap x_i \fC \neq \emptyset$.
\ecor
This characterization also appears in \cite[Theorem~6.6]{OP} (the countability assumption on $\fC$ in \cite{OP} is not necessary).
\setlength{\parindent}{0cm} \setlength{\parskip}{0.5cm}

Let us furthermore present a sufficient condition for local contractiveness. 
\blemma
\label{lem:bOmega-LocContr}
$I_l \ltimes \bOmega$ is locally contractive if and only if $I_l \bltimes \bOmega$ is locally contractive.
\setlength{\parindent}{0.5cm} \setlength{\parskip}{0cm}

$I_l \ltimes \bOmega$ is locally contractive if for all $0 \neq e \in \cJ$ there exists $s \in I_l$ and $f_0, \dotsc, f_n \in \cJ \setminus \gekl{0}$ such that $f_i \leq e s^{-1} s$ for all $0 \leq i \leq n$, for all $1 \leq i \leq n$ and $f' \leq s f_i s^{-1}$ there exists $0 \leq j \leq n$ with $f' f_j \neq 0$, and $f_0 s f_i = 0$ for all $0 \leq i \leq n$. 
\elemma 
\setlength{\parindent}{0cm} \setlength{\parskip}{0cm}

\bproof
The first statement follows from Lemma~\ref{lem:quotIxO-->IuxO}. The second statement is an application of \cite[Theorem~6.5]{EP16}.
\eproof

As a consequence, we obtain the following sufficient condition for local contractiveness in the finitely aligned case.
\bcor
\label{cor:FinAlign_bOmega-LocContr}
Suppose that $\fC$ is finitely aligned. Then $I_l \ltimes \bOmega \cong I_l \bltimes \bOmega$ is locally contractive if for all $x \in \fC$ there exist $c, d \in \fC$ with $\mfd(c) = \mfd(d)$ and $y_0, \dotsc, y_n \in \fC$ such that $d y_i \fC \leq x \fC$ for all $0 \leq i \leq n$, for all $1 \leq i \leq n$ and $z \in \fC$ with $z \fC \subseteq c y_i \fC$ there exists $0 \leq j \leq n$ with $z \fC \cap d y_j \fC \neq \emptyset$, and $d y_0 \fC \cap c y_i \fC = \emptyset$ for all $0 \leq i \leq n$. 
\ecor
\setlength{\parindent}{0cm} \setlength{\parskip}{0.5cm}

Finally, we characterize topological freeness or effectiveness of boundary groupoids. First we present a general characterization for topological freeness of tight groupoids attached to inverse semigroups. To the best of the author's knowledge, such a characterization has not appeared before. We work in the setting of \cite{EP16}. Let $S$ be an inverse semigroup with zero and $E$ its semilattice of idempotents. As in \S~\ref{ss:InvSgp-GPD}, we write $\widehat{E}$ for the space of characters of $E$. As in \S~\ref{ss:boundary}, we write $\widehat{E}_{\max}$ for the maximal filters on $E$ and $\partial \widehat{E} \defeq \overline{\widehat{E}_{\max}} \subseteq \widehat{E}$. Note that $\partial \widehat{E}$ is denoted by $\widehat{E}_{\rm tight}$ in \cite{EP16}. The action $S \curvearrowright \widehat{E}$ is defined as in \S~\ref{ss:InvSgp-GPD} and restricts to an action $S \curvearrowright \partial \widehat{E}$ (see also \cite{EP16}, for instance). As in \S~\ref{ss:InvSgp-GPD}, we define $S * \widehat{E} \defeq \lbrace (s,\chi) \in S \times \widehat{E}: \: \chi(s^{-1}s) = 1 \rbrace$ and $S \ltimes \widehat{E} \defeq (S * \widehat{E}) / { }_{\sim}$, where we set $(s,\chi) \sim (t,\psi)$ if $\chi = \psi$ and there exists $e \in E$ with $\chi(e) = 1$ and $se = te$. As above, equivalence classes with respect to $\sim$ are denoted by $[\cdot]$. The groupoid structure is defined as in \S~\ref{ss:InvSgp-GPD}.

\bdefin
Set $S^c \defeq \menge{s \in S}{e (s e s^{-1}) \neq 0 \ \forall \ 0 \neq e \leq s^{-1} s}$.
\edefin
\blemma
\label{lem:S^c}
\begin{enumerate}
\item[(i)] $S^c$ is closed under inverses, i.e., $s \in S^c$ implies $s^{-1} \in S^c$.
\item[(ii)] For all $s, t \in S^c$, $st$ also lies in $S^c$.
\item[(iii)] Whenever $s \in S^c$ and $t \in S$, we have $t^{-1} s t \in S^c$.
\end{enumerate}
\elemma
\setlength{\parindent}{0cm} \setlength{\parskip}{0cm}

\bproof
(i) is straightforward to prove. To prove (ii), take $0 \neq e \leq (st)^{-1}(st)$. Then $e \leq t^{-1}t$. Hence $e (tet^{-1}) \neq 0$. Moreover, $e (tet^{-1}) \leq s^{-1}s$. Thus $0 \neq e (tet^{-1}) s e (tet^{-1}) s^{-1} = e (tet^{-1}) (s e s^{-1}) s tet^{-1} s^{-1} \leq e s tet^{-1} s^{-1}$. For (iii), take $0 \neq e \leq (t^{-1} s t)^{-1} (t^{-1} s t)$. Then $e \leq t^{-1} t$, so that $0 \neq t e t^{-1} \leq s^{-1} s$. Since $s \in S^c$, we deduce that $tet^{-1} s tet^{-1} s^{-1} \neq 0$. Hence it follows that $e (t^{-1}stet^{-1}s^{-1}t) \neq 0$, as desired.
\eproof

In the following, we write $S^c \ltimes \partial \widehat{E} \defeq \lbrace [s,\chi] \in S \ltimes \partial \widehat{E}: \: s \in S^c \rbrace$. We start with a preparatory observation.

\blemma
\label{lem:S^cxhatE}
We have $S^c \ltimes \partial \widehat{E} \subseteq {\rm Iso}(S \ltimes \partial \widehat{E})$.
\elemma
\bproof
Take $s \in S$ and $\chi \in \partial \widehat{E}$ with $\chi(s^{-1}s) = 1$ and $s.\chi \neq \chi$. Since $\widehat{E}_{\max}$ is dense in $\partial \widehat{E}$, we may assume that $\chi \in \widehat{E}_{\max}$. $s.\chi \neq \chi$ implies that there exists $e \in E$ with $\chi(e) = 1$ and $s.\chi(e) = 0$, i.e., $\chi(s^{-1} e s) = 0$. Since $\chi \in \widehat{E}_{\max}$, the analogue of Lemma~\ref{lem:Omegamax} implies that there exists $f \in E$ with $\chi(f) = 1$ and $f (s^{-1} e s) = 0$. Hence $s f s^{-1} e s s^{-1} = 0$. Moreover, $\chi(s^{-1} s) = 1$ implies that $\chi(f s^{-1} s) = 1$ and thus $\chi(e f s^{-1} s) = 1$, so that $e f s^{-1} s \neq 0$. Clearly, we have $e f s^{-1} s \leq s^{-1} s$. Furthermore, $(e f s^{-1} s) s (e f s^{-1} s) s^{-1} = f s^{-1} s (e f s^{-1} s s^{-1}) s e s^{-1} = 0$. We conclude that $s \notin S^c$, as desired.
\eproof

\btheo
\label{thm:InvSgpBdTopFree}
The following are equivalent:
\begin{enumerate}
\item[(i)] $S \ltimes \partial \widehat{E}$ is topologically free. 
\item[(ii)] $S^c \ltimes \partial \widehat{E}$ is topologically free.
\item[(iii)] For all $s \in S^c$, $e, f_1, \dotsc, f_n \in E$ with $f_i \leq e \leq s^{-1}s$ for all $1 \leq i \leq n$ such that there exists $0 \neq f \leq e$ with $f f_i = 0$ for all $1 \leq i \leq n$, there exists $0 \neq f' \leq e$ with $f' f_i = 0$ for all $1 \leq i \leq n$ and $s f' = f'$.
\end{enumerate}
\etheo
\bproof
(i) $\Rarr$ (ii) follows as in the proof of Theorem~\ref{thm:IxOmega-topfree} because $S^c \ltimes \partial \widehat{E}$ is an open subgroupoid of $S \ltimes \partial \widehat{E}$.
\setlength{\parindent}{0.5cm} \setlength{\parskip}{0cm}

For (ii) $\Rarr$ (i), assume that $S \ltimes \partial \widehat{E}$ is not topologically free. Then there exists $s \in S$ and an open set $U \subseteq \partial \widehat{E}$ with $[s,U] \subseteq {\rm Iso}(S \ltimes \partial \widehat{E}) \setminus \partial \widehat{E}$. Take $\psi \in \widehat{E}_{\max} \cap U$ and $t \in S$ with $\psi(tt^{-1}) = 1$. Assume that $t s t^{-1} \notin S^c$. Then there exists $0 \neq f \leq (t^{-1} s^{-1} t)^{-1} (t^{-1} s^{-1} t)$ with $f (t^{-1} s^{-1} t f t^{-1} s t) = 0$. By the analogue of Lemma~\ref{lem:Omegamax}, there exists $\eta_t \in \widehat{E}_{\max}$ with $\eta_t(t f t^{-1}) = 1$. Thus $\eta_t (t^{-1} t) = 1$. Moreover, $f (t^{-1} s^{-1} t f t^{-1} s t) = 0$ implies $(t f t^{-1}) (s^{-1} t f t^{-1} s) = 0$. Hence if $\eta_t(s^{-1}s) = 1$, then $s.\eta_t \neq \eta_t$. Applying this reasoning to all $t \in S$ with $\psi(tt^{-1}) = 1$, we obtain a set $\gekl{\eta_t}_t \subseteq \widehat{E}_{\max}$ with $\eta_t(tt^{-1}) = 1$ for all such $t$. It follows by maximality that $\psi$ lies in the closure of $\gekl{\eta_t}_t$. As $\psi \in U$, this implies that $\eta_t \in U$ for some $t$. In particular, $\eta_t(s^{-1}s) = 1$, which implies $s.\eta_t \neq \eta_t$. This however contradicts the assumption that $[s,U] \subseteq {\rm Iso}(S \ltimes \partial \widehat{E})$. We conclude that there exists $t \in S$ with $\psi(tt^{-1}) = 1$ and $t^{-1}st \in S^c$. The latter implies $tt^{-1}stt^{-1} \in S^c$ by Lemma~\ref{lem:S^c}. Set $V \defeq U \cap \partial \widehat{E}(tt^{-1})$. Then $V$ is not empty because $\psi \in V$. Moreover, we claim $[s,V] \subseteq [tt^{-1}stt^{-1},V]$. Indeed, given $\zeta \in V$, we have $\zeta(tt^{-1}) = 1$ as well as $\zeta = s.\zeta$, so that $\zeta(s^{-1} tt^{-1} s) = s.\zeta(tt^{-1}) = 1$. Hence $\eta(s^{-1} t t^{-1} s t t^{-1}) = 1$. In addition, $tt^{-1}stt^{-1} = s (s^{-1} tt^{-1}stt^{-1})$. This shows that $(s,\zeta) \sim (tt^{-1}stt^{-1},\zeta)$, as desired. We conclude that $[s,V] \subseteq [tt^{-1}stt^{-1},V] \subseteq S^c \ltimes \partial \widehat{E}$. Hence $[s,V] \subseteq [s,U] \subseteq {\rm Iso}(S \ltimes \partial \widehat{E}) \setminus \partial \widehat{E}$ implies $[s,V] \subseteq {\rm Iso}(S^c \ltimes \partial \widehat{E}) \setminus \partial \widehat{E}$. This shows that $S^c \ltimes \partial \widehat{E}$ is not topologically free.

(ii) $\LRarr$ (iii): Lemma~\ref{lem:S^cxhatE} implies $S^c \ltimes \partial \widehat{E} \subseteq {\rm Iso}(S \ltimes \partial \widehat{E})$. Thus $S^c \ltimes \partial \widehat{E}$ is topologically free if and only if every non-empty bisection of $S^c \ltimes \partial \widehat{E}$ has non-empty intersection with $\partial \widehat{E}$. Every non-empty bisection contains a basic open set of the form $[s,\partial \widehat{E}(e; f_1, \dotsc, f_n)]$. $\partial \widehat{E}(e; f_1, \dotsc, f_n)$ is not empty if and only if there exists $\chi \in \widehat{E}_{\max}$ with $\chi(e) = 1$ and $\chi(f_i) = 0$ for all $1 \leq i \leq n$. By the analogue of Lemma~\ref{lem:Omegamax}, the latter holds if and only if there exists $f \in E$ with $f \leq e$ and $f f_i = 0$ for all $1 \leq i \leq n$ such that $\chi(f) = 1$. Hence we may assume that $s, e, f_1, \dotsc, f_n$ are exactly as in (iii). Now $[s,\partial \widehat{E}(e; f_1, \dotsc, f_n)] \cap \partial \widehat{E} \neq \emptyset$ if and only if $[s,\partial \widehat{E}(e; f_1, \dotsc, f_n)] \cap \widehat{E}_{\max} \neq \emptyset$. We claim that the last statement is equivalent to (iii). Indeed, if there exists $\chi \in \widehat{E}_{\max}$ with $\chi(e) = 1$, $\chi(f_1) = \dotso = \chi(f_n) = 0$ and $[s,\chi] = \chi$, then there exists $f \in E$ with $f \leq e$, $f f_i = 0$ for all $1 \leq i \leq n$ and $\chi(f) = 1$. $[s,\chi] = \chi$ implies that there exists $\ti{f} \in E$ with $\chi(\ti{f}) = 1$ and $s \ti{f} = \ti{f}$. Now $f' \defeq f \ti{f}$ has the desired properties. Conversely, if (iii) holds, then by the analogue of Lemma~\ref{lem:Omegamax}, there exists $\chi \in \widehat{E}_{\max}$ with $\chi(f') = 1$. It follows that $\chi \in \partial \widehat{E}(e; f_1, \dotsc, f_n)$, and $s f' = f'$ implies $[s,\chi] = \chi$.
\eproof

To complete the picture, we state the following characterization of effectiveness of $S^c \ltimes \partial \widehat{E}$. It follows from Lemma~\ref{lem:S^cxhatE} and also appears (implicitly) in \cite[\S~4]{EP16}. 
\blemma
$S^c \ltimes \partial \widehat{E}$ is effective if and only if for all $s \in S^c$, there exist $e_1, \dotsc, e_n \in E$ with $e_i \leq s^{-1} s$ and $s e_i = e_i$ for all $1 \leq i \leq n$, such that for all $f \leq s^{-1} s$, there exists $1 \leq i \leq n$ with $f e_i \neq 0$.
\elemma
\setlength{\parindent}{0cm} \setlength{\parskip}{0cm}

The following summarizes our findings, and combines them with the results in \cite[\S~4]{EP16}.
\bcor
\label{cor:SbE-eff,tf}
Consider the following statements:
\begin{enumerate}
\item[(i)] $S \ltimes \partial \widehat{E}$ is effective.
\item[(ii)] $S^c \ltimes \partial \widehat{E}$ is effective.
\item[(iii)] For all $s \in S^c$, there exist $e_1, \dotsc, e_n \in E$ with $e_i \leq s^{-1} s$ and $s e_i = e_i$ for all $1 \leq i \leq n$, such that for all $f \leq s^{-1} s$, there exists $1 \leq i \leq n$ with $f e_i \neq 0$.
\item[(iv)] $S \ltimes \partial \widehat{E}$ is topologically free. 
\item[(v)] $S^c \ltimes \partial \widehat{E}$ is topologically free.
\item[(vi)] For all $s \in S^c$, $e, f_1, \dotsc, f_n \in E$ with $f_i \leq e \leq s^{-1}s$ for all $1 \leq i \leq n$ such that there exists $0 \neq f \leq e$ with $f f_i = 0$ for all $1 \leq i \leq n$, there exists $0 \neq f' \leq e$ with $f' f_i = 0$ for all $1 \leq i \leq n$ and $s f' = f'$.
\end{enumerate}
Then (i) $\Rarr$ (ii) $\LRarr$ (iii) $\Rarr$ (iv) $\LRarr$ (v) $\LRarr$ (vi). If $S \ltimes \partial \widehat{E}$ is Hausdorff, then all these statements are equivalent. If $\partial \widehat{E} = \widehat{E}_{\max}$, then (i) $\LRarr$ (ii) $\LRarr$ (iii).
\ecor
\bproof
All this follows from what we showed above, except for the very last statement, which follows from \cite[Theorem~4.10]{EP16}.
\eproof

Corollary~\ref{cor:SbE-eff,tf} applied to $S = I_l$ yields a characterization when $I_l \ltimes \bOmega$ is topologically free and a necessary condition for effectiveness of $I_l \ltimes \bOmega$. Now we turn to $I_l \bltimes \bOmega$.

\bdefin
We set $\bIl^f \defeq \menge{s \in \bIl}{\exists \, f \in \cJ \text{ with } 0 \neq f \leq s^{-1} s}$.
\edefin
\blemma
\label{lem:Ilfuni}
Given $s \in \bIl$, $s$ lies in $\bIl^f$ if and only if $s^{-1}$ lies in $\bIl^f$ if and only if there exists $\chi \in \Omegamax$ with $\chi(s^{-1}s) = 1$.
\elemma
\bproof
The first equivalence is easy to see. If $s^{-1} s = e \setminus \bigcup_{i=1}^n f_i$, then existence of $\chi \in \Omegamax$ with $\chi(s^{-1}s) = 1$ implies that $\chi(e) = 1$ and $\chi(f_1) = \dotso = \chi(f_n) = 0$. Hence there exists $f \in \cJ$ with $f \leq e$, $ff_1 = \dotso = ff_n = 0$ and $\chi(f) = 1$ by Lemma~\ref{lem:Omegamax}. Thus $0 \neq f \leq s^{-1} s$. Conversely, if there exists $f \in \cJ$ with $0 \neq f \leq s^{-1} s$, then Lemma~\ref{lem:Omegamax} implies that there exists $\chi \in \Omegamax$ with $\chi(f) = 1$ and thus $\chi(s^{-1}s) = 1$.
\eproof

Now we define the analogue of $S^c$ or $I_l^c$.
\bdefin
We define $\bIl^c \defeq \big \lbrace s \in \bIl^f : \: e (s e s^{-1}) \neq 0 \ \forall \, e \in \cJ \text{ with } 0 \neq e \leq s^{-1} s \big \rbrace$.
\edefin

The following is the analogue of Lemma~\ref{lem:S^c}. The proof is similar.
\blemma
\label{lem:Ilcuni}
\begin{enumerate}
\item[(i)] $\bIl^c$ is closed under inverses, i.e., $s \in \bIl^c$ implies $s^{-1} \in \bIl^c$.
\item[(ii)] For all $s, t \in \bIl^c$ with $st \in \bIl^f$, $st$ also lies in $\bIl^c$.
\item[(iii)] For all $s \in \bIl^c$ and $t \in \bIl$ with $t^{-1} s t \in \bIl^f$, we have $t^{-1} s t \in \bIl^c$.
\end{enumerate}
\elemma

As explained in Remark~\ref{rem:Iluni-acts}, we have an identification $I_l \bltimes \bOmega \cong \bIl \ltimes \bOmega$. In the following, we work with $\bIl \ltimes \bOmega$. We set $\bIl^c \ltimes \bOmega \defeq \menge{[s,\chi] \in \bIl \ltimes \bOmega}{s \in \bIl^c}$.
\btheo
\label{thm:IbxbOmega-eff}
We have ${\rm Iso}(\bIl \ltimes \bOmega)^{\circ} = \bIl^c \ltimes \bOmega$.

The following are equivalent:
\begin{enumerate}
\item[(i)] $I_l \bltimes \bOmega \cong \bIl \ltimes \bOmega$ is effective.
\item[(ii)] $\bIl^c \ltimes \bOmega = \bOmega$.
\item[(iii)] For all $s \in \bIl^c$ there exist $\varepsilon_1, \dotsc, \varepsilon_n \in \bcJ$ with $s \varepsilon_i = \varepsilon_i$ for all $1 \leq i \leq n$ such that for all $g \in \cJ$ with $0 \neq g \leq s^{-1}s$, there exists $1 \leq i \leq n$ with $g \varepsilon_i \neq 0$.
\end{enumerate}
\etheo
\bproof
Let us prove ${\rm Iso}(\bIl \ltimes \bOmega)^{\circ} = \bIl^c \ltimes \bOmega$. We first show \an{$\supseteq$}: Take $s \in \bIl^c$ and $\chi \in \Omegamax$ with $\chi(s^{-1}s) = 1$. As we have seen in the proof of Lemma~\ref{lem:Ilfuni}, there exists $f \in \cJ$ with $f \leq s^{-1} s$ and $\chi(f) = 1$. We claim that for all $e \leq f$, $\chi(e) = 1$ implies that $\chi(s e s^{-1}) = 1$. Indeed, if $\chi(s e s^{-1}) = 0$, then Lemma~\ref{lem:Omegamax} implies that there exists $e' \in \cJ$ with $e' (s e s^{-1}) = 0$ and $\chi(e') = 1$. Thus $ee' (s ee' s^{-1}) = 0$, while $ee' \neq 0$ since $\chi(ee') = 1$. This contradicts $s \in \bIl^c$. Now given $g \in \cJ$, $s.\chi(g) = 1$ $\LRarr$ $s.\chi(sfs^{-1}g) = 1$ $\LRarr$ $\chi(f s^{-1} g s) = 1$ $\Rarr$ $\chi(s f s^{-1} g) = 1$ $\Rarr$ $\chi(g) = 1$. Maximality implies $s.\chi = \chi$. Hence $s.\chi = \chi$ for all $\chi \in \bOmega(s^{-1}s)$. Now we show \an{$\subseteq$}: Take $s, t \in \bIl$ and $U = \bOmega(tt^{-1})$ with $U \subseteq \bOmega(s^{-1}s)$ and $[s,U] \subseteq {\rm Iso}(\bIl \ltimes \bOmega)$. Without loss of generality, we may assume $s, t \in \bIl^f$ because of Lemma~\ref{lem:Ilfuni}. Take $\chi \in U$. Then $t^{-1}.\chi((t^{-1}st)^{-1}(t^{-1}st)) = 1$. Lemma~\ref{lem:Ilfuni} implies that $t^{-1}st \in \bIl^f$. If $t^{-1}st \notin \bIl^c$, then there exists $0 \neq e \in \cJ$ with $e \leq (t^{-1}st)(t^{-1}st)^{-1}$ such that $e (t^{-1}s^{-1}t e t^{-1}st) = 0$. Hence $(t e t^{-1})(s^{-1}t e t^{-1}s) = 0$. $e \leq (t^{-1}st)(t^{-1}st)^{-1}$ implies $e \leq t^{-1} t$, hence $t e t^{-1} \in \cJ$. Take $\psi \in \Omegamax$ with $\psi(t e t^{-1}) = 1$. This is possible by Lemma~\ref{lem:Omegamax}. Then $\psi(tt^{-1}) = 1$ and thus $\psi \in U$. In particular, $s.\psi = \psi$. However, $\psi(s^{-1} t e t^{-1} s) = s.\psi(t e t^{-1}) = \psi(t e t^{-1}) = 1$, which contradicts $(t e t^{-1})(s^{-1}t e t^{-1}s) = 0$. We conclude that $t^{-1}st \in \bIl^c$. Since there exists $f \in \cJ$ with $0 \neq f \leq (t^{-1}st)(t^{-1}st)^{-1}$, we have $0 \neq t f t^{-1} \leq tt^{-1}s^{-1}tt^{-1}st^{-1}$, which implies $tt^{-1}stt^{-1} \in \bIl^f$. Since $t^{-1}st \in \bIl^c$, Lemma~\ref{lem:Ilcuni} implies that $tt^{-1}stt^{-1} \in \bIl^c$. Finally, $[s,\chi] = [tt^{-1}stt^{-1},\chi]$ in $\bIl \ltimes \bOmega$ because $\chi(s^{-1}tt^{-1}stt^{-1}) = 1$ and $tt^{-1}stt^{-1} = s (s^{-1}tt^{-1}stt^{-1})$. We conclude that $[s,\chi] = [tt^{-1}stt^{-1},\chi] \in \bIl^c \ltimes \bOmega$. This shows \an{$\subseteq$}.
\setlength{\parindent}{0.5cm} \setlength{\parskip}{0cm}

(i) $\LRarr$ (ii) follows from what we just proved. To prove (ii) $\LRarr$ (iii), observe that $\bIl^c \ltimes \bOmega = \bOmega$ if and only if for all $s \in \bIl^c$ and $\chi \in \bOmega$ with $\chi(s^{-1}s) = 1$, there exists $\varepsilon \in \bcJ$ with $\varepsilon \leq s^{-1}s$, $\chi(\varepsilon) = 1$ and $s \varepsilon = \varepsilon$. Fix $s \in \bIl^c$. By compactness of $\bOmega(s^{-1}s)$, we deduce that there are $\varepsilon_1, \dotsc, \varepsilon_n \in \bcJ$ with $s \varepsilon_i = \varepsilon_i$ for all $1 \leq i \leq n$ such that $\bOmega(s^{-1}s) = \bigcup_{i=1}^n \bOmega(\varepsilon_i)$. We claim that this last equality is equivalent to the statement that for all $g \in \cJ$ with $0 \neq g \leq s^{-1}s$, there exists $1 \leq i \leq n$ with $g \varepsilon_i \neq 0$. Indeed, given $g \in \cJ$ with $0 \neq g \leq s^{-1}s$, Lemma~\ref{lem:Omegamax} provides $\chi \in \Omegamax$ with $\chi(g) = 1$. This implies $\chi(s^{-1}s) = 1$ and hence $\chi(\varepsilon_i) = 1$ for some $1 \leq i \leq n$. We deduce that $g \varepsilon_i \neq 0$. Conversely, assume that there exists $\chi \in \bOmega$ with $\chi(s^{-1}s) = 1$ and $\chi(\varepsilon_i) = 0$ for all $1 \leq i \leq n$. By density, we may assume that $\chi \in \Omegamax$. The proof of Lemma~\ref{lem:Ilfuni} shows that there exists $f \in \cJ$ with $f \leq s^{-1}s$ such that $\chi(f) = 1$. Write $\varepsilon_i = e_i \setminus \bigcup_j f_{ij}$ for some $e_i, f_{ij} \in \cJ$. Either $\chi(e_i) = 0$ or $\chi(e_i) = 1 = \chi(f_{ij})$ for some $j$. In the first case, Lemma~\ref{lem:Omegamax} implies that there exists $g_i \in \cJ$ with $\chi(g_i) = 1$ and $g_i e_i = 0$, which implies $g_i \varepsilon_i = 0$. In the second case, set $g_i \defeq f_{ij}$. Then we also obtain $\chi(g_i) = 1$ and $g_i \varepsilon_i = 0$. Now define $g \defeq f g_1 \dotsm g_n$. Then $\chi(g) = 1$ implies $g \neq 0$. By construction, we have $g \leq s^{-1}s$. In addition, we have $g \varepsilon_i = 0$ for all $1 \leq i \leq n$.
\eproof
\setlength{\parindent}{0cm} \setlength{\parskip}{0cm}

Let us now specialize to the finitely aligned case.
\bcor
\label{cor:FinAlign_bOmega-tf,eff}
Suppose that $\fC$ is finitely aligned.
\begin{enumerate}
\item[(i)] $I_l \ltimes \bOmega \cong I_l \bltimes \bOmega$ is topologically free if and only if for all $c, d \in \fC$ with $\mfd(c) = \mfd(c)$ and $\mft(c) = \mft(d)$ with the property that for all $x \in \fC$ with $cx \in d \fC$, we have $cx\fC \cap dx\fC \neq \emptyset$, there exist $y, z_1, \dotsc, z_n \in \fC$ with $z_i \in y \fC$ for all $1 \leq i \leq n$, $cy \in d \fC$, and with the property that there exists $z \in y \fC$ with $z \fC \cap z_i \fC = \emptyset$ for all $1 \leq i \leq n$, such that there exists $z' \in \fC$ with $z' \fC \cap z_i \fC = \emptyset$ for all $1 \leq i \leq n$ and $cz' = dz'$.
\item[(ii)] $I_l \ltimes \bOmega \cong I_l \bltimes \bOmega$ is effective if and only if for all $c, d \in \fC$ with $\mfd(c) = \mfd(c)$ and $\mft(c) = \mft(d)$, $\delta = a \fC \setminus \bigcup_h a_h \fC$ for $a \in \fC$ and $\gekl{a_h} \subseteq \fC$ finite with the property that $\delta \subseteq c^{-1}(d \fC \cap c \fC)$ and for all $x \in \fC$ with $x \fC \subseteq \delta$, we have $cx\fC \cap dx\fC \neq \emptyset$, there exist $\varepsilon_i = b_i \fC \setminus \bigcup_j b_{ij} \fC$ (for some $b_i \in \fC$ and finite subset $\gekl{b_{ij}} \subseteq \fC$), $1 \leq i \leq n$, with the property that $\varepsilon_i \subseteq \delta$ and $cy = dy$ for all $y \in \varepsilon_i$ for all $1 \leq i \leq n$, such that for all $z \in \delta$, we have $z \fC \cap \varepsilon_i \neq \emptyset$ for some $1 \leq i \leq n$.
\end{enumerate}
\ecor
\setlength{\parindent}{0cm} \setlength{\parskip}{0cm}

Note that for certain finitely aligned left cancellative small categories, effectiveness of the boundary groupoid has been characterized in \cite[Theorem~6.4]{OP}.
\setlength{\parindent}{0cm} \setlength{\parskip}{0cm}

\bcor
\label{cor:Cons.CbG}
If the conditions in Lemma~\ref{lem:HdIxbOmega} and Corollary~\ref{cor:SbE-eff,tf}~(iv) for $S = I_l$ are satisfied, then $C^*_r(I_l \ltimes \bOmega)$ has the intersection property. If, in addition, the condition in Lemma~\ref{lem:bOmega-min} is satisfied, then $C^*_r(I_l \ltimes \bOmega)$ is simple. And if, in addition, the condition in Lemma~\ref{lem:bOmega-min} and the condition in Lemma~\ref{lem:bOmega-LocContr} are satisfied, then $C^*_r(I_l \ltimes \bOmega)$ is purely infinite simple.
\setlength{\parindent}{0.5cm} \setlength{\parskip}{0cm}

If the conditions in Lemma~\ref{lem:HdIbxbOmega} and Theorem~\ref{thm:IbxbOmega-eff}~(iii) are satisfied, then $C^*_r(I_l \bltimes \bOmega)$ has the intersection property. If, in addition, the condition in Lemma~\ref{lem:bOmega-min} is satisfied, then $C^*_r(I_l \bltimes \bOmega)$ is simple. And if, in addition, the condition in Lemma~\ref{lem:bOmega-min} and the condition in Lemma~\ref{lem:bOmega-LocContr} are satisfied, then $C^*_r(I_l \bltimes \bOmega)$ is purely infinite simple.

Suppose that $\fC$ is finitely aligned. If the condition in Corollary~\ref{cor:FinAligned-HdbOmega} and one of the conditions in Corollary~\ref{cor:FinAlign_bOmega-tf,eff} are satisfied, then $C^*_r(I_l \ltimes \Omega) \cong C^*_r(I_l \bltimes \Omega)$ has the intersection property. If, in addition, the condition in Corollary~\ref{cor:FinAlign_bOmega-min} is satisfied, then $C^*_r(I_l \ltimes \bOmega) \cong C^*_r(I_l \bltimes \bOmega)$ is simple. And if, in addition, the condition in Corollary~\ref{cor:FinAlign_bOmega-min} and the condition in Corollary~\ref{cor:FinAlign_bOmega-LocContr} are satisfied, then $C^*_r(I_l \ltimes \bOmega) \cong C^*_r(I_l \bltimes \bOmega)$ is purely infinite simple.
\ecor
\setlength{\parindent}{0cm} \setlength{\parskip}{0.5cm}

We also present the following observation, which is inspired by \cite[Lemma~5.7.10 and Theorem~5.7.16]{CELY}.
\bprop
\label{prop:smallest}
$\bOmega$ is the smallest non-empty closed invariant subspace of $\hcJ$ if and only if $\bOmega$ is the smallest non-empty closed invariant subspace of $\Omega$ if and only if for all $\mfv, \mfw \in \fC^0$, we have $\mfw \fC \mfv \neq \emptyset$. In that case, $I_l \ltimes \bOmega$ is purely infinite in the sense of \cite[Definition~4.9]{Mat15} if and only if for all $\mfv \in \fC^0$, there exist $a, b \in \mfv \fC$ with $a \fC \cap b \fC = \emptyset$.
\eprop
\setlength{\parindent}{0cm} \setlength{\parskip}{0cm}

\bproof
First suppose that $\bOmega$ is the smallest non-empty closed invariant subspace of $\Omega$. Then for every $\mfv \in \fC^0$, the orbit closure of $\chi_{\mfv}$ contains all of $\Omegamax$. Hence, given $\mfw \in \fC^0$, there exists $s \in I_l$ such that $s.\chi_{\mfv}(\mfw \fC) = 1$. It follows that there exists $x \in \mfw \fC$ such that $s^{-1}(x)$ is defined and $\chi_{\mfv}(s^{-1}(x)) = 1$, which implies that $s^{-1}(x) = u$ for some $u \in \mfv \fC^*$. Hence $x = s(u)$ and thus $x \in \mfw \fC \mfd(u)$, which implies that $x u^{-1} \in \mfw \fC \mfv$. 
\setlength{\parindent}{0.5cm} \setlength{\parskip}{0cm}

Now suppose that for all $\mfv, \mfw \in \fC^0$, we have $\mfw \fC \mfv \neq \emptyset$. Take $\chi \in \hcJ_{\max}$ and $\eta \in \hcJ$ arbitrary. Let $\mfv \in \fC^0$ be such that $\eta(\mfv \fC) = 1$. For every $e \in \cJ$ with $\chi(e) = 1$, take $a \in e$ and $b \in \mfd(a) \fC \mfv$. Then $\mfd(ab) = \mfv$, so that $ab.\eta$ is defined, and we have $ab.\eta(ab) = 1$, which implies $ab.\chi(e) = 1$. Set $\chi_e \defeq ab.\eta$. By maximality, $\chi$ must lie in $\overline{\menge{\chi_e}{e \in \cJ, \, \chi(e) = 1}}$.

This concludes the proof of the first part. For the second claim, if there exists $\mfv \in \fC^0$ such that $a \fC \cap b \fC \neq \emptyset$ for all $a, b \in \mfv \fC$, then $\bOmega(\mfv)$ degenerates to a point and $I_l \ltimes \bOmega$ cannot be purely infinite. Conversely, given $\mfv \in \fC^0$ and a basic open set $U \subseteq \bOmega(\mfv)$ of the form $U = \bOmega(x;\mfy)$ for some $x \in \mfv \fC$, $\mfy \subseteq \mfv \fC$, we can find $\chi \in U \cap \Omegamax$. Lemma~\ref{lem:Omegamax} implies that there exists $z \in \mfv \fC$ such that $z \fC \leq x \fC$ and $z \fC \cap y \fC = \emptyset$ for all $y \in \mfy$. It follows that every $\eta \in \bOmega$ with $\eta(z \fC) = 1$ satisfies $\eta \in U$. Now take $a, b \in \mfd(z) \fC$ with $a \fC \cap b \fC = \emptyset$ and $a' \in \mfd(a) \fC \mfv$, $b' \in \mfd(b) \fC \mfv$. Then $aa' \fC \cap bb' \fC = \emptyset$ and $aa', bb' \in \mfd(z) \fC \mfv$. We conclude that $zaa'.U \subseteq U$, $zbb'.U \subseteq U$ and $zaa'.U \cap zbb'.U = \emptyset$.
\eproof
\setlength{\parindent}{0cm} \setlength{\parskip}{0.5cm}

Now let us prove Theorem~\ref{thm:intro_PinG}. Let $P$ be a submonoid of a group $G$. Denote by $1$ the identity element of $P$ and $G$. The remaining part of this section deals with the special case where $\fC = P$.
\bdefin
We set $G^c \defeq \menge{g \in G}{(pP) \cap (gpP) \neq \emptyset \quad \forall \ p \in P}$.
\edefin

\btheo
\label{thm:Gc=1}
The boundary groupoid $I_l \ltimes \bOmega$ for $P$ is effective if and only if $G^c = \gekl{1}$.
\etheo
\setlength{\parindent}{0cm} \setlength{\parskip}{0cm}

\bproof
We will use the same notation as in \cite[\S~5]{CELY}. 
\setlength{\parindent}{0.5cm} \setlength{\parskip}{0cm}

First, we show that for all $g \in G^c$ and $\chi \in U_{g^{-1}} \cap \bOmega$, we have $g.\chi =  \chi$. Without loss of generality, we may assume that $\chi \in \Omegamax$. $\chi \in U_{g^{-1}}$ implies that there exists $f \in \cJ$ with $g^{-1}(f) \in \cJ$ and $\chi(g^{-1}(f)) = 1$. Now take $e \in \cJ$ arbitrary. $g.\chi(e) = 1$ if and only if $g.\chi(ef) = 1$ if and only if $\chi(g^{-1}(ef)) = 1$. Note that $g^{-1}(ef)$ must lie in $\cJ$ as well. Assume that $\chi(e) = 0$. By Lemma~\ref{lem:Omegamax}, there exists $e' \in \cJ$ with $\chi(e') = 1$ and $e'e = 0$. Take $x \in e' g^{-1}(ef)$. Since $g \in G^c$, we know that $(xP) \cap (gxP) \neq \emptyset$. However, $xP \subseteq e'$ and $gxP \subseteq ef \subseteq e$. Hence we obtain a contradiction to $e'e = 0$. This means that for all $g \in G^c$, $\gekl{g} \times U_{g^{-1}} \cap \bOmega \subseteq {\rm Iso}(I_l \ltimes \bOmega)^{\circ}$. So if $I_l \ltimes \bOmega$ is effective, then we must have $G^c = \gekl{1}$.

Conversely, assume that $G^c = \gekl{1}$. If $I_l \ltimes \bOmega$ is not effective, then there exists $1 \neq g \in G$ and a non-empty open set $U \subseteq U_{g^{-1}} \cap \bOmega$ with $\gekl{g} \times U \subseteq {\rm Iso}(I_l \ltimes \bOmega)$. Take $\chi \in U \cap \Omegamax$ and $e \in \cJ$ with $\chi(e) = 1$. Pick an element $x \in e$. $g \neq 1$ implies that $x^{-1} g x \neq 1$, so that $x^{-1} g x \notin G^c$. Hence there exists $p \in P$ with $(pP) \cap (x^{-1} g x p P) = \emptyset$. It follows that $(g^{-1}xpP) \cap (xpP) = \emptyset$. By Lemma~\ref{lem:Omegamax}, there exists $\chi_e \in \Omegamax$ with $\chi_e(xpP) = 1$. By maximality, $\chi$ lies in the closure of $\menge{\chi_e}{e \in \cJ, \, \chi(e) = 1}$. Hence we can find $\chi_e$ in $U$. Now we have $1 = \chi_e(xpP) = g.\chi_e(xpP) = \chi_e(g^{-1}xpP)$. This contradicts $(g^{-1}xpP) \cap (xpP) = \emptyset$. We conclude that $I_l \ltimes \bOmega$ must be effective, as desired.
\eproof
\setlength{\parindent}{0cm} \setlength{\parskip}{0cm}

\bcor[see {\cite[\S~5.7]{CELY}}]
\label{cor:Gc=1}
Let $P$ be a submonoid of a group $G$. If $G^c = \gekl{1}$, then $\partial C^*_{\lambda}(P)$ is simple, and $\partial C^*_{\lambda}(P)$ is purely infinite simple unless $P = \gekl{1}$.
\ecor
\setlength{\parindent}{0cm} \setlength{\parskip}{0.5cm}

\bremark
As Marcelo Laca and Camila F. Sehnem kindly informed me, Theorem~\ref{thm:Gc=1} and Corollary~\ref{cor:Gc=1} also follow from \cite[Proposition~6.18 and Corollary~6.19]{LS}.
\eremark

\section{C*-algebras of Garside categories}
\label{s:C*Gar}

\subsection{Preliminaries on Garside families in small categories}

In the following, we give a brief introduction to Garside categories and collect a few facts about Garside families in small categories which will be needed later on. Our exposition follows \cite{Deh15}, where the reader will find more details (but note that our convention is opposite to the one in \cite{Deh15}, as explained at the beginning of \S~\ref{ss:cat-rr-C}).

Let $\fC$ be a left cancellative small category. 
\bdefin
Given $a, b \in \fC$, we write $a \preceq b$ if $a$ is a left divisor of $b$, i.e., $b \in a \fC$. We write $a \prec b$ if $b \fC \subsetneq a \fC$. We write $a \tipreceq b$ if $a$ is a right divisor of $b$, i.e., $b \in \fC a$. We write $a \tiprec b$ if $\fC b \subsetneq \fC a$. 
\setlength{\parindent}{0.5cm} \setlength{\parskip}{0cm}

We write $a =^* b$ if $a \in b \fC^*$ (which is equivalent to $a \fC = b \fC$). 
\edefin
\setlength{\parindent}{0cm} \setlength{\parskip}{0cm}

\bdefin
$\fC$ is called left Noetherian if there exists no infinite sequence $\dotso \prec a_3 \prec a_2 \prec a_1$. $\fC$ is called right Noetherian if there exists no infinite sequence $\dotso \tiprec a_3 \tiprec a_2 \tiprec a_1$. $\fC$ is called Noetherian if it is both left and right Noetherian. 
\edefin

The following notion already came up in Remark~\ref{rem:fin-align--mcm}.
\bdefin
Given $a, b, c \in \fC$, $c$ is called an mcm of $a$ and $b$ if $c$ is minimal with respect to $\preceq$ among $\menge{d \in \fC}{a \preceq d \text{ and } b \preceq d}$, i.e., $a \preceq c$, $b \preceq c$, and for all $d \in \fC$ with $a \preceq d$ and $b \preceq d$, if $d \preceq c$ then $d =^* c$.
\setlength{\parindent}{0.5cm} \setlength{\parskip}{0cm}

We write $\mcm(a,b)$ for the set of all mcms of $a, b \in \fC$.
\edefin
\setlength{\parindent}{0cm} \setlength{\parskip}{0cm}

It would be more precise to use the term \an{right mcm}. However, there will be no danger of confusion because left mcms will not appear in this paper.

\bdefin
A finite or infinite sequence $s_1, s_2, \dotsc$ in $\fC$ is called a path if $\mfd(s_k) = \mft(s_{k+1})$ for all $k$. We denote this path by $s_1 s_2 \dotsm$.
\edefin

\bdefin
Given $\fS \subseteq \fC$, we set $\fS^{\sharp} \defeq \fS \fC^* \cup \fC^*$.
\edefin
\bdefin
A subset $\fS \subseteq \fC$ is closed under right comultiples if for all $r, s \in \fS$ and $a \in \fC$ with $r \preceq a$, $s \preceq a$, there exists $t \in \fS$ with $r \preceq t$, $s \preceq t$ and $t \preceq a$.
\edefin

\bdefin
Suppose $\fS \subseteq \fC$ is closed under right comultiples, $\fS \cup \fC^*$ generates $\fC$ and $\fS^{\sharp}$ is closed under right divisors. 
\setlength{\parindent}{0.5cm} \setlength{\parskip}{0cm}

A path $s_1 \dotsm s_l \in \fS^{\sharp}$ is called $\fS$-normal if for all $1 \leq k \leq l-1$ and $r \in \fS$, if $r \preceq s_k s_{k+1}$ then $r \preceq s_k$. We also call a $\fS$-normal path a $\fS$-normal word.

For $a \in \fC$, a $\fS$-normal decomposition (or $\fS$-normal form) of $a$ is given by a $\fS$-normal path $s_1 \dotsm s_l$ in $\fS^{\sharp}$ with $a = s_1 \dotsm s_l$.

$\fS$ is called a Garside family if every element in $\fC$ admits a $\fS$-normal decomposition.
\edefin
\setlength{\parindent}{0cm} \setlength{\parskip}{0cm}

Our assumption on $\fS$ is justified by \cite[Chapter~IV, Proposition~1.23]{Deh15}. We used a simplified version of normal decomposition (compare \cite[Chapter~III]{Deh15} for the general definition), which is allowed because of \cite[Chapter~IV, Proposition~1.20]{Deh15}. In the following, whenever $\fS$ is understood, we will drop the prefix \an{$\fS$-} (for instance, we will write \an{normal} instead of \an{$\fS$-normal}).

\bremark
If $s_1 \dotsm s_l$ is normal, then for all $1 \leq j \leq k \leq l$ and $r \in \fS$, if $r \preceq s_j \dotsm s_k$ then $r \preceq s_j$ (see \cite[Chapter~III, Proposition~1.12]{Deh15}).
\eremark

\bremark
If $\fS$ is a Garside family, then $\fC^* \fS \subseteq \fS^{\sharp}$ by \cite[Chapter~III, Proposition~1.39]{Deh15}.
\eremark

By choosing one representative in each $=^*$-class, we may (and will) always arrange that for all $s_1, s_2 \in \fS$, $s_1 =^* s_2$ implies $s_1 = s_2$ (see \cite[Chapter~III, Corollary~1.34]{Deh15}, i.e., $\fS$ is $=^*$-transverse.

\bprop[see {\cite[Chapter~III, Corollary~1.27]{Deh15}}]
\label{prop:Snormal-unique}
If $\fS$ is a Garside family which is $=^*$-transverse, then every $a \in \fC \setminus \fC^*$ admits a unique normal decomposition $s_1 \dotsm s_l$ with $s_k \in \fS \setminus \fC^*$ for all $1 \leq k \leq l-1$ and $s_l \in \fS^{\sharp} \setminus \fC^*$. 
\eprop

\bdefin
Suppose that $\fS$ is a Garside family which is $=^*$-transverse. Given $a \in \fS$, we define $\Vert a \Vert \defeq 0$ if $a \in \fC^*$ and $\Vert a \Vert \defeq l$ if $s_1 \dotsm s_l$ is the unique normal decomposition of $a$ as in Proposition~\ref{prop:Snormal-unique}.
\edefin

There are many criteria which ensure that a subset $\fS$ of $\fC$ is a Garside family. We mention the following example.
\bprop[see {\cite[Chapter~IV, Corollary~2.26]{Deh15}}]
Suppose that $\fC$ is left cancellative and Noetherian. Then $\fS \subseteq \fC$ is a Garside family if and only if $\fS \cup \fC^*$ generates $\fC$ and $\fS^{\sharp}$ is closed under mcms and right divisors.
\eprop

\bdefin
Given $\fS \subseteq \fC$, $s \in \fS$ is called an $\fS$-head of $a \in \fC$ if $s$ is a maximal left divisor in $\fS$ of $a$ with respect to $\preceq$, i.e., $s \preceq a$, and every $r \in \fS$ with $r \preceq a$ satisfies $r \preceq s$. 
\edefin

If $\fS$ is a Garside family, then by \cite[Chapter~IV, Proposition~1.24]{Deh15}, every non-invertible element $a$ of $\fC$ admits an $\fS$-head, which is unique if $\fS$ is $=^*$-transverse. In that case, the $\fS$-head will be denoted by $H(a)$.

The following are immediate consequences of \cite[Chapter~III, Proposition~1.49]{Deh15}.
\bprop
\label{prop:LeftMult-H}
Suppose that $\fS$ is a Garside family which is $=^*$-transverse with $\fS \cap \fC^* = \emptyset$. 
\setlength{\parindent}{0.5cm} \setlength{\parskip}{0cm}

Given a path $a_1, \dotsc, a_n$ in $\fS$, we have
$
 H(a_1 \dotsm a_n) =  H(a_1 \dotso H(a_{n-2} H(a_{n-1} a_n)) \dotso )
$.
If $r_1 r_2 \dotsm$ is a normal path in $\fS$ and $a_1 \dotsm a_n$ is a path in $\fS$, then the normal form of $a_1 \dotsm a_n r_1 r_2 \dotsm$ starts with $H(a_1 \dotsm a_n r_1)$.
\eprop 
\setlength{\parindent}{0cm} \setlength{\parskip}{0.5cm}

\subsection{Classification of closed invariant subspaces}

From now on, let $\fC$ be a countable left cancellative category. Suppose $\fS \subseteq \fC$ is a subset which generates $\fC$. Given an infinite path $w = s_1 s_2 \dotsm$ in $\fS$, we write $w_n \defeq s_1 \dotsm s_n$, $w_{= n} \defeq s_n$ and $w_{>n} \defeq s_{n+1} s_{n+2} \dotsm$. Moreover, we set $\Omegainf \defeq \Omega \setminus \menge{\chi_x}{x \in \fC}$.

\blemma
\label{lem:infpaths}
\begin{enumerate}
\item[(i)] Define a function $\chi_w: \: \cJ \to \gekl{0,1}$ by setting, for all $e \in \cJ$, $\chi_w(e) \defeq 1$ if $w_n \in e$ for some $n$ and $\chi_w(e) = 0$ otherwise. Then $\chi_w \in \Omega$.
\end{enumerate}
\setlength{\parindent}{0cm} \setlength{\parskip}{0cm}

Now assume that $\fC$ is finitely aligned.
\begin{enumerate}
\item[(ii)] Every $\chi \in \Omegainf$ is of the form $\chi_w$ for some infinite path $w$ in $\fS$.
\item[(iii)] Given $c, d \in \fC$ with $\mft(d) = \mft(c)$ and an infinite path $w$ in $\fS$, we have that $c d^{-1}. \chi_w$ is defined if and only if there exists $n$ such that $d \preceq c w_n$. In that case, if we have $d x_n = c w_n$, then $c d^{-1}. \chi_w = \chi_{c x_n w_{>n}}$.
\end{enumerate}
\elemma
\bproof
(i) It is immediate that $\chi_w$ is indeed a character satisfying the condition in Definition~\ref{def:Omega}.
\setlength{\parindent}{0.5cm} \setlength{\parskip}{0cm}

(ii) Given $\chi \in \Omega \setminus \menge{\chi_x}{x \in \fC}$, set $\cF_{\rm p} \defeq \menge{x \in \fC}{\chi(x \fC) = 1}$. Write $\cF_{\rm p} = \gekl{x_1, x_2, \dotsc}$. Now define $w^{(1)} \defeq x_1$. For $n > 1$, since $\chi(w^{(n-1)} \fC \cap x_n \fC) = 1$, there exists $w^{(n)} \in \mcm(w^{(n-1)},x_n)$ such that $\chi(w^{(n)} \fC) = 1$. Here we use that $\fC$ is finitely aligned and $\chi \in \Omega$. Since $\chi \notin \menge{\chi_x}{x \in \fC}$, we may assume that $w^{(n)} \notin \fC^*$ for all $n$. Thus we can write $w^{(n)} = s_1^{(n)} \dotsm s_{l_n}^{(n)}$. Now set  
$$
 w \defeq s_1^{(1)} \dotsm s_{l_1}^{(1)} s_1^{(2)} \dotsm s_{l_2}^{(2)} s_1^{(3)} \dotsm s_{l_3}^{(3)} \dotsm
$$
We claim that $\chi = \chi_w$. Indeed, for $x \in \fC$, $\chi_w(x \fC) = 1$ if and only if $x \preceq w^{(n)}$ for some $n$ if and only if $\chi(x \fC) = 1$. The last equivalence follows from $\chi(w^{(n)}) = 1$ and $x_n \preceq w^{(n)}$ for all $n$. Now our claim follows from Lemma~\ref{lem:FinAlign}~(ii).

(iii) $c d^{-1}. \chi_w$ is defined if and only if $\chi_{w}(d) = 1$, and hence the first claim follows. If $d x_n = c w_n$, then $\chi_w = \chi_{d x_n w_{>n}}$. Thus $c d^{-1}. \chi_w = c d^{-1} d. \chi_{x_n w_{>n}} = \chi_{c x_n w_{>n}}$.
\eproof
\setlength{\parindent}{0cm} \setlength{\parskip}{0cm}

\bdefin
We call a subset $\fS \subseteq \fC$ locally finite if $\mfv \fS$ is finite for all $\mfv \in \fC^0$.
\setlength{\parindent}{0.5cm} \setlength{\parskip}{0cm}

We call a subset $\fS \subseteq \fC$ locally bounded if for every $\mfv \in \fC^0$, there exists no infinite sequence $s_1, s_2, \dotsc$ in $\mfv \fS$ with $s_1 \prec s_2 \prec \dotso$.
\edefin
\setlength{\parindent}{0cm} \setlength{\parskip}{0cm}

In the following, a finite or infinite word $x$ in $\fS$ is called normal if $x = s_1 s_2 \dotsm$ for a normal path $s_1 s_2 \dotsm$. In case $x$ is an infinite normal path, we set $\Vert x \Vert \defeq \infty$. Given two words $x = s_1 s_2 \dotsm$ and $y = t_1 t_2 \dotsm$ in $\fS$, equality of words $x=y$ means $s_1 = t_1, s_2 = t_2, \dotsc$.
\blemma
\label{lem:NormalWordsExUni}
Suppose that $\fC$ is finitely aligned and that $\fS$ is a Garside family in $\fC$ with $\fS \cap \fC^* = \empty$ which is $=^*$-transverse and locally bounded. Then every $\chi \in \Omega \setminus \menge{\chi_{\mfv}}{\mfv \in \fC^0}$ is of the form $\chi_x$ for some normal word $x$. Moreover, given two normal words $x$ and $y$, we have $\chi_x = \chi_y$ if and only if $x = y$.
\elemma
\bproof
By Lemma~\ref{lem:infpaths}~(ii), there exists $w = r_1 r_2 \dotsm$ (where $r_1, r_2, \dotsc \in \fS$) with $\chi = \chi_w$. Set $s_1^{(n)} \defeq H(w_n)$. As $s_1^{(n)} \preceq w_{n+1}$ and $s_1^{(n+1)}$ is the maximal left divisor of $w_{n+1}$, we must have $s_1^{(n)} \preceq s_1^{(n+1)}$. Because $\fS$ is $=^*$-transverse and locally bounded, it follows that the sequence $s_1^{(1)}, s_2^{(1)}, \dotsc$ must be eventually constant, say eventually equal to $s_1$. We introduce the notation $H(w) \defeq s_1$. Now define $s_2 \defeq H(s_1^{-1} w), \dotsc, s_n \defeq H(s_{n-1}^{-1} \dotsm s_2^{-1} s_1^{-1} w), \dotsc$. Set $x = s_1 s_2 \dotsc$. By construction, $s_1 s_2 \dotsm$ is normal. We claim that $\chi = \chi_x$. Indeed, proceed inductively on $n$ to show that for all $n$, there exists $N(n)$ such that $x_n \preceq w_{N(n)}$: This is true by construction for $n=1$. Now suppose that $x_n \preceq w_{N(n)}$. Then $\chi_{x_n^{-1} w_{N(n)} w_{>N(n)}}(s_{n+1}) = 1$ implies that $s_{n+1} \preceq x_n^{-1} w_{N(n+1)}$ for some sufficiently big $N(n+1)$. It follows that $x_{n+1} \preceq w_{N(n+1)}$, as desired. Thus, given $z \in \fC$, if $\chi_x(z \fC) = 1$, then $z \preceq x_n$ for some $n$, and hence $z \preceq x_n \preceq w_{N(n)}$, which implies $\chi_w(z \fC) = 1$. Let us show that, conversely, given $z \in \fC$, if $\chi_w(z \fC) = 1$, then $\chi_x(z \fC) = 1$. By construction and because of $x_n \preceq w_{N(n)}$, the normal form of $w_{N(n)}$ starts with $s_1 s_2 \dotsm s_n$. It then follows from \cite[Chapter~III, Proposition~1.14]{Deh15} that $w_n \prec s_1 \dotsm s_n = x_n$. Hence, given $z \in \fC$, if $\chi_w(z \fC) = 1$, then $z \preceq w_n$ for some $n$, so that $z \preceq w_n \preceq x_n$, and thus $\chi_x(z \fC) = 1$.
\setlength{\parindent}{0.5cm} \setlength{\parskip}{0cm}

Suppose that $x = s_1 s_2 \dotsm$ for a normal path $s_1 s_2 \dotsm$ and $y = t_1 t_2 \dotsm$ for a normal path $t_1 t_2 \dotsm$. If $\chi_x = \chi_y$, then $\chi_y(s_1 \fC) = 1$, hence $s_1 \preceq y_n$ for some $n$. But since $t_1$ is the maximal left divisor of $y_n$, this implies $s_1 \preceq t_1$. By symmetry, we also obtain $t_1 \preceq s_1$. It follows that $s_1 =^* t_1$. As $\fS$ is $=^*$-transverse, we conclude that $s_1 = t_1$. Now proceed inductively, applying the previous argument to $s_1^{-1} x = s_2 s_3 \dotsm$ and $t_1^{-1} y = t_2 t_3 \dotsm$, using that $\chi_{s_1^{-1}x} = s_1^{-1}.\chi_x = t_1^{-1}.\chi_y = \chi_{t_1^{-1}y}$.
\eproof
\setlength{\parindent}{0cm} \setlength{\parskip}{0cm}

From now on, in the remaining part of \S~\ref{s:C*Gar}, we will assume the following:
\setlength{\parindent}{0cm} \setlength{\parskip}{0.5cm}

\textbf{Standing assumptions:} $\fC$ is finitely aligned and $\fS$ is a Garside family in $\fC$ with $\fS \cap \fC^* = \emptyset$ which is $=^*$-transverse and locally bounded. 

Let $\cW$ be the set of (non-empty) normal words in $\fS$. Then Lemma~\ref{lem:NormalWordsExUni} implies that there is a one-to-one correspondence
$
 \cW \amalg \fC^0 \isom \Omega, \, w \ma \chi_w, \mfv \ma \chi_{\mfv}
$.
In the following, to simplify notation, given $x \in \fC$, we will denote $x \fC$ by $x$.
\setlength{\parindent}{0cm} \setlength{\parskip}{0cm}

\blemma
Given a sequence $w^{(i)} \in \cW$ and $w \in \cW$, we have $\lim_i \chi_{w^{(i)}} = \chi_w$ if and only if for all $n$, $w_n$ is maximal with respect to $\preceq$ among $\lbrace v \in \fC : \: \Vert v \Vert \leq n, \, v \preceq w_n^{(i)} \text{ for almost all } i \rbrace$.
\setlength{\parindent}{0.5cm} \setlength{\parskip}{0cm}

If $\fS$ is locally finite, then $\gekl{\chi_x}$ is open for every $x \in \fC$, $\Omega_{\infty}$ is closed, and given a sequence $w^{(i)} \in \cW$ and $w \in \cW$, we have $\lim_i \chi_{w^{(i)}} = \chi_w$ if and only if for all $n$, $w_n = w_n^{(i)}$ for almost all $i$.
\elemma
\setlength{\parindent}{0cm} \setlength{\parskip}{0cm}

\bproof
The first part follows from the following observations: Given $v \in \fC$ with $\Vert v \Vert \leq n$, we have that $\lim_i \chi_{w^{(i)}}(v) = 1$ if and only if $v \preceq w_n^{(i)}$ for almost all $i$, whereas $\chi_w(v) = 1$ if and only if $v \preceq w_n$.
\setlength{\parindent}{0.5cm} \setlength{\parskip}{0cm}

For the second part, if $\fS$ is locally finite, we have $\gekl{\chi_x} = \Omega(x; \mfd(x) \fS)$ is open for all $x \in \fC$ and thus $\Omegainf = \Omega \setminus \bigcup_{x \in \fC} \gekl{\chi_x}$ is closed. Moreover, we claim that $w_n$ is maximal among $\lbrace v \in \fC : \: \Vert v \Vert \leq n, \, v \preceq w_n^{(i)} \text{ for almost all } i \rbrace $ with respect to $\preceq$ if and only if $w_n = w_n^{(i)}$ for almost all $i$. Indeed, by deleting the first few elements of the sequence, we may assume that $w_n \preceq w_n^{(i)}$ for all $i$. If we do not have $w_n = w_n^{(i)}$ for almost all $i$, then by passing to a subsequence, we may arrange $w_n \prec w_n^{(i)}$ for all $i$. Since $\fS$ is locally finite, $\fS^n$ is also locally finite. Hence, by further passing to a subsequence, we may arrange that $w_n^{(i)}$ is constant, say equal to $v$. It follows that $w_n \prec v$, and thus $\chi_{w^{(i)}}$ does not converge to $\chi_w$ by the first part. This is a contradiction.
\eproof
\setlength{\parindent}{0cm} \setlength{\parskip}{0cm}

Given a sequence $s^{(i)}$ in $\fS$ and $s \in \fS \cup \fC^0$, we write $\lim_i s^{(i)} = s$ if $s$ is maximal with respect to $\preceq$ among $\menge{r \in \fS \cup \fC^0}{r \preceq s^{(i)} \text{ for almost all } i}$.
\setlength{\parindent}{0cm} \setlength{\parskip}{0.5cm}

In the following, we denote $I_l \ltimes \Omega$ by $\cG$.
\bprop
\label{prop:invclosure}
Suppose that $\cV \subseteq \cW \amalg \fC^0$. 
\setlength{\parindent}{0.5cm} \setlength{\parskip}{0cm}

Given a normal word in $\fS$, $w = s_1 s_2 \dotsm$, we have $\chi_w \in \overline{\cG.\menge{\chi_v}{v \in \cV}}$ if and only if for all $j$, there exists a sequence $v^{(i)}$ in $\cV$ such that for all $i$, there exist $a_i \in \fC$ and $m_i \in \Nz$ with $\Vert v^{(i)} \Vert < m_i$ or $a_i \in \fC \mfd(v^{(i)})$ if $v^{(i)} \in \fC$ such that, if we set $s_j^{(i)} \defeq H(a_i v^{(i)}_{= m_i})$ in the first case or $s_j^{(i)} \defeq H(a_i)$ in the second case, then $\lim_i s_j^{(i)} = s_j$.

For $\mfw \in \fC^0$, we have $\chi_{\mfw} \in \overline{\cG.\menge{\chi_v}{v \in \cV}}$ if and only if $\mfw \in \cV$ or there exists a sequence $v^{(i)}$ in $\cV$ such that for all $i$, there exist $a_i \in \fC$ and $m_i \in \Nz$ with $\Vert v^{(i)} \Vert < m_i$ or $a_i \in \fC \mfd(v^{(i)})$ if $v^{(i)} \in \fC$ such that, if we set $s^{(i)} \defeq H(a_i v^{(i)}_{= m_i})$ in the first case or $s^{(i)} \defeq H(a_i)$ in the second case, then $\lim_i s^{(i)} = \mfw$.
\eprop
\setlength{\parindent}{0cm} \setlength{\parskip}{0cm}

\bproof
We prove the first claim, the argument for the second claim is analogous. 
\setlength{\parindent}{0.5cm} \setlength{\parskip}{0cm}

For \an{$\Larr$}, write $\chi_{\ti{v}^{(i)}} \defeq a_i (r_1^{(i)} \dotsm r_{m_{i-1}}^{(i)})^{-1}. \chi_{v^{(i)}}$ in the first case and $\chi_{\ti{v}^{(i)}} \defeq a_i (v^{(i)})^{-1}. \chi_{v^{(i)}}$ in the second case. Then $\chi_{\ti{v}^{(i)}} \in \cG.\menge{\chi_v}{v \in \cV}$. Moreover, the normal form of $\ti{v}^{(i)}$ starts with $s_j^{(i)}$. By compactness of $\Omega(\mft(s_j))$, we may assume without loss of generality that $\chi_{\ti{v}^{(i)}}$ converges to $\chi_x$. Then $\chi_x \in \overline{\cG.\menge{\chi_v}{v \in \cV}}$. The assumption $\lim_i s_j^{(i)} = s_j$ implies that the normal form of $x$ starts with $s_j$. Now set $\chi_{w^{(j)}} \defeq (s_1 \dotsm s_{j-1}).\chi_x$. Then the normal form of $w^{(j)}$ starts with $s_1 \dotsm s_j$. Since the normal form of $w$ also starts with $s_1 \dotsm s_j$, we conclude that $\lim_j \chi_{w^{(j)}} = \chi_w$, as desired.

Now we show \an{$\Rarr$}. Without loss of generality, we may assume that $j=1$. Assume that we can find $c_i, d_i \in \fC$ with $\mft(d_i) = \mft(c_i)$ and $v^{(i)} \in \cV$ such that $\lim_i d_i^{-1} c_i.\chi_{v^{(i)}} = \chi_w$. By Lemma~\ref{lem:infpaths}~(iii), we can write $d_i^{-1} c_i.\chi_{v^{(i)}} = \chi_{a_i v_{>N_i}^{(i)}}$ for some $a_i \in \fC$ or $d_i^{-1} c_i.\chi_{v^{(i)}} = \chi_{a_i}$ for some $a_i \in \fC \mfd(v^{(i)})$ (which implies $v^{(i)} \in \fC$). The normal decomposition of $a_i v_{>N_i}^{(i)}$ starts with $s_1^{(i)} \defeq H(a_i v_{>N_i}^{(i)})$ in the first case, and the normal decomposition of $a_i$ starts with $s_1^{(i)} \defeq H(a_i)$ in the second case. Then $\lim_i d_i^{-1} c_i.\chi_{v^{(i)}} = \chi_w$ implies that $\lim_i s_1^{(i)} = s_1$, as desired.
\eproof
\setlength{\parindent}{0cm} \setlength{\parskip}{0cm}

\bdefin
Let $\fT \subseteq \fS$ and $\fD \subseteq \fC^0$. 
\begin{enumerate}
\item[(i)] $(\fT,\fD)$ is called admissible if for all $t \in \fT$, there exists $t' \in \fT$ such that $t t'$ is normal or $\mfd(t) \in \fD$.
\item[(ii)] $(\fT,\fD)$ is called $H$-invariant if for all $a \in \fC \setminus \fC^*$ and $x \in \fT \cup \fD$ with $\mfd(a) = \mft(x)$, $H(ax)$ lies in $\fT$.
\item[(iii)] $(\fT,\fD)$ is called $\max_{\preceq}^{\infty}$-closed if for every sequence $t_i$ in $\fT$, if $\lim_i t_i$ exists in $\fS$, then $\lim_i t_i \in \fT \cup \fD$.
\end{enumerate}
\edefin

\bdefin
Given $X \subseteq \Omega$, let $\cV \subseteq \cW \amalg \fD^0$ be such that $X = \menge{\chi_v}{v \in \cV}$. 
\setlength{\parindent}{0.5cm} \setlength{\parskip}{0cm}

Define $\fT(X) \defeq \menge{t \in \fS}{t = v_{=i} \text{ for some } v \in \cV \cap \cW \text{ and } i \in \Nz}$ and $\fD(X) \defeq \cV \cap \fC^0 = \menge{\mfv \in \fC^0}{\chi_{\mfv} \in X}$.
\edefin
\setlength{\parindent}{0cm} \setlength{\parskip}{0cm}

\blemma
\label{lem:admissible}
$(\fT,\fD)$ is admissible if and only if there exists $X \subseteq \Omega$ such that $\fT = \fT(X)$ and $\fD = \fD(X)$.
\setlength{\parindent}{0.5cm} \setlength{\parskip}{0cm}

$(\fT(X),\fD(X))$ is $H$-invariant and $\max_{\preceq}^{\infty}$-closed if and only if $X$ is $\cG$-invariant and closed.
\elemma
\setlength{\parindent}{0cm} \setlength{\parskip}{0cm}

\bproof
For the first claim, to see \an{$\Larr$}, suppose that $t = v_{=i}$. Then $t v_{= i+1}$ is normal if $\Vert v \Vert \geq i+1$, and $\mfd(s) \in \fD$ if $\Vert v \Vert = i$. For \an{$\Rarr$}, given $t \in \fT$, we can inductively construct an infinite normal word in $\fT$ starting with $t$ or a finite normal word with $\mfd$ in $\fD$.
\setlength{\parindent}{0.5cm} \setlength{\parskip}{0cm}

For the second claim, \an{$\Rarr$} follows from Proposition~\ref{prop:invclosure}. For \an{$\Larr$}, if $X = \menge{\chi_v}{v \in \cV}$ is $\cG$-invariant, then $(\fT,\fD)$ is $H$-invariant because of Proposition~\ref{prop:LeftMult-H}, and if $X$ is closed, then by compactness of $\Omega(\mfv)$ for all $\mfv \in \fC^0$ and Proposition~\ref{prop:invclosure}, $(\fT,\fD)$ is $\max_{\preceq}^{\infty}$-closed, where $\fT = \menge{t \in \fS}{t = v_{=i} \text{ for some } v \in \cV \cap \cW \text{ and } i \in \Nz}$ and $\fD = \cV \cap \fC^0 = \menge{\mfv \in \fC^0}{\chi_{\mfv} \in X}$.
\eproof
\setlength{\parindent}{0cm} \setlength{\parskip}{0cm}

Given $(\fT,\fD)$, there is a smallest $H$-invariant and $\max_{\preceq}^{\infty}$-closed pair $(\overline{\fT},\overline{\fD})$ containing $(\fT,\fD)$, which can be constructed by adjoining elements $H(ax)$ (for $a \in \fC \setminus \fC^*$, $x \in \fT \cup \fD$ with $\mfd(a) = \mft(x)$) and $\lim_i t_i$ (for sequences $t_i$ in $\fT$) step by step and taking the union at the end. Similarly, given $(\fT,\fD)$, there is a biggest admissible pair $(\check{\fT},\check{\fD})$ contained in $(\fT,\fD)$, which can be constructed by deleting elements $t$ for which there does not exist $t' \in \fT$ such that $t t'$ is normal and for which $\mfd(t) \notin \fD$ step by step and taking the intersection at the end.
\bcor
If $(\fT,\fD)$ is admissible, then $(\overline{\fT},\overline{\fD})$ is admissible. In addition, $(\overline{\fT},\overline{\fD})$ is obtained by first adjoining elements $H(ax)$ (for $a \in \fC \setminus \fC^*$, $x \in \fT \cup \fD$ with $\mfd(a) = \mft(x)$) and then adjoining elements of the form $\lim_i t_i$, i.e., this process does not have to be repeated.
\setlength{\parindent}{0.5cm} \setlength{\parskip}{0cm}

Suppose that $(\fT,\fD)$ is $H$-invariant and $\max_{\preceq}^{\infty}$-closed. Then $(\check{\fT},\check{\fD})$ $H$-invariant and $\max_{\preceq}^{\infty}$-closed.
\ecor
\setlength{\parindent}{0cm} \setlength{\parskip}{0cm}

\bproof
Let us prove the first claim. By Lemma~\ref{lem:admissible}, there exists $X \subseteq \Omega$ such that $\fT = \fT(X)$ and $\fD = \fD(X)$. It now follows from Proposition~\ref{prop:invclosure} that $(\overline{\fT}, \overline{\fD}) = (\fT(\overline{\cG.X}), \fD(\overline{\cG.X}))$, and that $(\overline{\fT},\overline{\fD})$ is obtained by first adjoining elements $H(ax)$ (for $a \in \fC \setminus \fC^*$, $x \in \fT \cup \fD$ with $\mfd(a) = \mft(x)$) and then adjoining elements of the form $\lim_i t_i$ (i.e., this process does not have to be repeated). 
\setlength{\parindent}{0.5cm} \setlength{\parskip}{0cm}

Now we prove the second claim. It follows from the first claim that $(\overline{\check{\fT}},\overline{\check{\fD}})$ is admissible. Moreover, since $(\fT,\fD)$ is $H$-invariant and $\max_{\preceq}^{\infty}$-closed, we must have $\overline{\check{\fT}} \subseteq \fT$ and $\overline{\check{\fD}}) \subseteq \fD$. Hence, by maximality of $(\check{\fT},\check{\fD})$, we conclude that $(\check{\fT}, \check{\fD}) = (\overline{\check{\fT}}, \overline{\check{\fD}})$.
\eproof
\setlength{\parindent}{0cm} \setlength{\parskip}{0cm}

\bdefin
Let $\fT \subseteq \fS$ and $\fD \subseteq \fC^0$. We set 
$
 X(\fT,\fD) \defeq \menge{\chi_v}{v_{=i} \in \fT \ \forall \ i \in \Nz} \cup \menge{\chi_{\mfv}}{\mfv \in \fD}
$.
\edefin

\btheo
\label{thm:ClInvSubsp}
The maps $X \ma (\fT(X),\fD(X))$ and $X(\fT,\fD) \mapsfrom (\fT,\fD)$ establish an inclusion-preserving one-to-one correspondence between $\cG$-invariant, closed subspaces of $\Omega$ and admissible, $H$-invariant, $\max_{\preceq}^{\infty}$-closed pairs $(\fT,\fD)$ with $\fT \subseteq \fS$ and $\fD \subseteq \fC^0$.
\etheo
\setlength{\parindent}{0cm} \setlength{\parskip}{0cm}

Here we write $(\fT_1,\fD_1) \subseteq (\fT_2,\fD_2)$ if $\fT_1 \subseteq \fT_2$ and $\fD_1 \subseteq \fD_2$.
\bproof
Lemma~\ref{lem:admissible} implies that these maps are well-defined. Moreover, Lemma~\ref{lem:admissible} implies that $\fT = \fT(X(\fT,\fD))$ and $\fD = \fD(X(\fT,\fD))$. Finally, it remains to show that $X(\fT(X),\fD(X)) = X$. \an{$\supseteq$} is clear. For \an{$\subseteq$}, take $\chi_v \in X(\fT(X),\fD(X))$. If $v \in \fC^0$, then $\chi_v \in X$. If $v \in \cW$, then, for all $i \in \Nz$ with $i \leq \Vert v \Vert$, there exists $w \in \cW$ with $\chi_w \in X$ and $w_{=i} = v_{=i}$. It then follows from Proposition~\ref{prop:invclosure} that $\chi_v \in X$ because $X$ is $\cG$-invariant and closed. It is clear that the maps preserve inclusions.
\eproof

\bcor
\label{cor:ClInvSubsp_locfin}
If $\fS$ is locally finite, then the maps $X \ma (\fT(X),\fD(X))$ and $X(\fT,\fD) \mapsfrom (\fT,\fD)$ establish an inclusion-preserving one-to-one correspondence between $\cG$-invariant, closed subspaces of $\Omega$ and admissible, $H$-invariant pairs $(\fT,\fD)$ with $\fT \subseteq \fS$ and $\fD \subseteq \fC^0$.
\ecor
\bproof
This follows from Theorem~\ref{thm:ClInvSubsp} because every pair $(\fT,\fD)$ is automatically $\max_{\preceq}^{\infty}$-closed as $\fS$ is locally finite.
\eproof

Next we characterize $\cG$-invariant, closed subsets which are contained in the boundary.
\bdefin
Let $\fD_{\max}$ be the subset of all $\mfv \in \fC^0$ with $\mfv \fC = \mfv \fC^*$. Define
\begin{align*}
 \fT_{\rm Max} \defeq &\menge{t \in \fS}{\forall \ F \subseteq \mfd(t) \fS \text{ with } \# \, F < \infty, \, t \ti{t} \in \fS \ \forall \ \ti{t} \in F \ \exists \ x \in \fC \text{ with } x \cap \ti{t} = \emptyset \ \forall \ \ti{t} \in F}\\
 \cup &\menge{t \in \fS}{\exists \text{ finite normal path } v \text{ and } i \leq \Vert v \Vert \text{ with } v_{=i} = t \text{ and } \mfd(v) \in \fD_{\max}}.
\end{align*}
\edefin

\bprop
We have $\fD_{\max} = \fD(\Omegamax)$ and $\fT(\Omegamax) \subseteq \fT_{\rm Max} \subseteq \fT(\bOmega)$.
\eprop
\setlength{\parindent}{0cm} \setlength{\parskip}{0cm}

\bproof
The first claim is clear. Now take $t \in \fT(\Omegamax)$. Then there exists $\chi_w \in \Omegamax$ with $w_{=i} = t$. If $w$ is a finite normal path, then $\mfd(w) \in \fD_{\max}$. Now suppose that $w$ is an infinite normal path. As $\Omegamax$ is $\cG$-invariant, we may assume $i=1$. Take a finite subset $F \subseteq \mfd(t) \fS$ with $t \ti{t} \in \fS$ for all $ti{t} \in F$. Then $\chi_w(t \ti{t}) = 0$ for all $\ti{t} \in F$. By Lemma~\ref{lem:Omegamax}, there exists $x \in \fC$ with $\chi_w(t x) = 1$ and $t x \cap t \ti{t} = \emptyset$ for all $\ti{t} \in F$. Hence $x \cap \ti{t} = \emptyset$ for all $\ti{t} \in F$, and we conclude that $t \in \fT_{\rm Max}$. To show $\fT_{\rm Max} \subseteq \fT(\bOmega)$, take $t \in \fT_{\rm Max}$. If there exists a finite normal path $v$ and $i \leq \Vert v \Vert$ with $v_{=i} = t$ and $\mfd(v) \in \fD_{\max}$, then $\chi_v \in \Omegamax$ and $t \in \fT(\bOmega)$. Now suppose that for all finite subsets $F \subseteq \mfd(t) \fS$ with $t \ti{t} \in \fS$ for all $\ti{t} \in F$, there exists $x \in \fC$ with $x \cap \ti{t} = \emptyset$ for all $\ti{t} \in F$. It suffices to show that there exists $\chi_w \in \bOmega$ with $w \in \cW$ and $w_1 = t$. Order all finite subsets $F \subseteq \mfd(t) \fS$ with $t \ti{t} \in \fS$ for all $\ti{t} \in F$ by inclusion and find $\chi_F \in \Omegamax$ with $\chi_F(tx) = 1$. Such $\chi_F$ exist by Lemma~\ref{lem:Omegamax}. By compactness of $\Omega(\mft(t))$, we may assume that $\lim_F \chi_F = \chi_w \in \bOmega$. Then $\chi_F(t) = 1$ for all $F$ while $\chi_F(t \ti{t}) = 0$ whenever $\ti{t} \in F$. It follows that $w_1 = t$, as desired.
\eproof

The following are immediate consequences.
\bcor
\label{cor:ClInv-B}
We have $(\fT(\bOmega),\fD(\bOmega)) = (\overline{\fT_{\rm Max}},\overline{\fD_{\max}})$. Under the correspondence in Theorem~\ref{thm:ClInvSubsp}, a $\cG$-invariant, closed subspace $X$ is contained in $\bOmega$ if and only if $(\fT(X),\fD(X)) \subseteq (\overline{\fT_{\rm Max}},\overline{\fD_{\max}})$.
\ecor

\blemma
\label{lem:TMax_locfin}
If $\fS$ is locally finite, then 
\begin{align*}
 \fT_{\rm Max} = &\menge{t \in \fS}{\exists \ x \in \fC \text{ with } x \cap \ti{t} = \emptyset \ \forall \ \ti{t} \in \mfd(s) \fS \text{ with } t \ti{t} \in \fS} \\
 \cup &\menge{t \in \fS}{\exists \text{ finite normal path } v \text{ and } i \leq \Vert v \Vert \text{ with } v_{=i} = t \text{ and } \mfd(v) \in \fD_{\max}}.
\end{align*}
Moreover, $\fT_{\rm Max} = \fT(\Omegamax)$.
\elemma
\bproof
The first claim follows since $\# \, \mfd(t) \fS < \infty$ as $\fS$ is locally finite. For the second claim, it suffices to show that $\fT_{\rm Max} \subseteq \fT(\Omegamax)$. Given $t \in \fT_{\rm Max}$, take $x \in \fC$ with $x \cap \ti{t} = \emptyset$ for all $\ti{t} \in \mfd(s) \fS$ with $t \ti{t} \in \fS$. By Lemma~\ref{lem:Omegamax}, there exists $\chi_w \in \Omegamax$ (where $w \in \cW$) with $\chi_w(tx) = 1$. It then follows that $w_1 = t$, as desired.
\eproof
\setlength{\parindent}{0cm} \setlength{\parskip}{0.5cm}

\subsection{Topological freeness and local contractiveness}

Let us establish a sufficient condition for topological freeness.

\bprop
\label{prop:Gar-TopFree}
Let $(\fT_1,\fD_1)$ and $(\fT_2,\fD_2)$ be admissible, $H$-invariant and $\max_{\preceq}^{\infty}$-closed pairs with $(\fT_1,\fD_1) \subseteq (\fT_2,\fD_2)$. Set $X_* \defeq X(\fT_*,\fD_*)$ for $* = 1,2$. Assume that $\mfv \fC^* \mfv = \mfv$ for all $\mfv \in \fD_2 \setminus \fD_1$, and that for all finite normal paths $a, b$ in $\fS$ with $\mft(a) = \mft(b)$, $\mfd(a) = \mfd(b)$, $a_1 \neq b_1$ and $s \in \fT_2 \setminus \fT_1$ with $\mfd(s) = \mfd(a)$, there exists $t \in \fT_2 \setminus \fT_1$ such that $s t$ is normal and $H(at) \neq H(bt)$. Then $I_l \ltimes (X_2 \setminus X_1)$ is topologically free.
\eprop
\setlength{\parindent}{0cm} \setlength{\parskip}{0cm}

\bproof
Given $c, d \in \fC$ with $\mfd(c) = \mfd(d)$ and a basic open set $U = (X_2 \setminus X_1)(x;\mfy) \subseteq \Omega(d)$, where $x \in \fC$ and $\mfy \subseteq \fC$ is a finite set, we want to show that $[cd^{-1},U] \cap U \neq \emptyset$ or there exists $\chi \in U$ with $cd^{-1}.\chi \neq \chi$. If there exists a finite normal word $w$ with $\chi_v \in U$, then we must have $\mfd(v) \in \fD_2 \setminus \fD_1$, and $cd^{-1}.\chi_v = \chi_v$ implies that $cd^{-1}(v) = vu$ for some $u \in \mfd(w) \fC^* \mfd(w) = \mfd(w)$ (the last equality holds by assumption). Hence $c d^{-1}(v) = v$ and thus $[cd^{-1},\chi_v] = \chi_v$. 
\setlength{\parindent}{0.5cm} \setlength{\parskip}{0cm}

Now suppose that there is an infinite normal word $v$ with $\chi_v \in U$. Set $L \defeq \max \menge{\Vert y \Vert}{y \in \mfy}$. Then we claim that for every normal word $v'$ with $\Vert v' \Vert \geq L$, $v'_L = v_L$ and $\chi_{v'} \in X_2 \setminus X_1$ imply that $\chi_{v'} \in U$. Indeed, given $y \in \mfy$, if $\chi_{v'}(y) = 1$, then $y \preceq v'_n$ for some $n$, which would imply $y \preceq v'_L = v_L$ by \cite[Chapter~III, Proposition~1.14]{Deh15}, contradicting $\chi_v(y) = 0$. Now if $c v_L = d v_L$, then $cd^{-1}(v_L) = v_L$ and hence $[cd^{-1},\chi_v] = \chi_v$. If $c v_L \neq d v_L$, then we can find finite normal paths $r, a, b$ with $a_1 \neq b_1$ such that $c v_L = r a$ and $d v_L = r b$. For $s = v_{=L}$, there exists by assumption $t \in \fT_2 \setminus \fT_1$. Thus we can find a normal word $w$ with $w_1 = t$ such that $\chi_w \in X_2 \setminus X_1$. Since $s t$ is normal, we obtain that $v_L t$ is normal, so that $\chi_{v_L w} \in U$. Now write $a w = H(at) z_a$ and $b w = H(bt) z_b$ for some normal words $z_a$ and $z_b$. We conclude that
$$
 c.\chi_{v_L w} = c v_L. \chi_w = r a. \chi_w = r. \chi_{H(at)z_a} \neq r. \chi_{H(bt)z_b} = r b. \chi_w = d v_L. \chi_w = d. \chi_{v_L w}.
$$
\eproof
\setlength{\parindent}{0cm} \setlength{\parskip}{0.5cm}

Next, we present a sufficient condition for local contractiveness.
\setlength{\parindent}{0cm} \setlength{\parskip}{0cm}

\bprop
\label{prop:XminusXLocallyContractive}
Suppose that $\fC$ is left Noetherian. Let $(\fT_1,\fD_1)$ and $(\fT_2,\fD_2)$ be two admissible, $H$-invariant and $\max_{\preceq}^{\infty}$-closed pairs with $(\fT_1,\fD_1) \subseteq (\fT_2,\fD_2)$. Set $X_* \defeq X(\fT_*,\fD_*)$ for $* = 1,2$. Assume that for every admissible, $H$-invariant and $\max_{\preceq}^{\infty}$-closed pair $(\fT,\fD)$ with $(\fT_1,\fD_1) \subseteq (\fT,\fD) \subsetneq (\fT_2,\fD_2)$, there exists an infinite normal path in $\fT_2 \setminus \fT$. Further suppose that for all finite normal paths $c$ in $\fT_2 \setminus \fT_1$, there exists a maximal element $s \in \fT_2 \setminus \fT_1$ with respect to $\preceq$ together with a normal path $p$ such that $cps$ is normal and two distinct normal paths $q_1$, $q_2$ such that $s q_1 c$ and $s q_2 c$ are normal. Then $I_l \ltimes (X_2 \setminus X_1)$ is locally contractive.
\eprop

\bproof
As above, suppose that we are given a basic open set $U = (X_2 \setminus X_1)(x;\mfy) \subseteq \Omega(c^{-1}(c \cap d))$, where $x \in \fC$ and $\mfy \subseteq \fC$ is a finite set. Our first assumption implies that there exists an infinite normal word $v$ such that $\chi_v \in U$. Set $L \defeq \max \menge{\Vert y \Vert}{y \in \mfy}$ and $c \defeq v_L$. As shown above, for every normal word $v'$ with $\Vert v' \Vert \geq L$, $v'_L = c$ and $\chi_{v'} \in X_2 \setminus X_1$ imply that $\chi_{v'} \in U$. 
\setlength{\parindent}{0.5cm} \setlength{\parskip}{0cm}

Now let $\gekl{t_i}$ be the minimal elements with respect to $\preceq$ among $\menge{t \in \fS}{s t \in \fT_1}$. The elements $t_i$ exist because $\fC$ is left Noetherian. Without loss of generality assume that $t_i \neq t_j$ for all $i \neq j$. Now we claim that $\# \gekl{t_i} < \infty$. If not, then we show that by passing to a subsequence, we may arrange $\lim_i s t_i = s t$ for some $t \in \fS \cup \fC^0$. Indeed, if $s$ is not maximal among $\menge{r \in \fS}{r \preceq s t_i \text{ for almost all } i}$, by passing to a subsequence we may arrange that there exists $t' \in \fS$ with $st' \preceq s t_i$ for all $i$. If $st'$ is not maximal among $\menge{r \in \fS}{r \preceq s t_i \text{ for almost all } i}$, then we obtain, by passing to a subsequence if necessary, an element $t'' \in \fS$ with $t' \prec t''$ such that $st'' \preceq st_i$ for all $i$. Continuing this way, we obtain a sequence $t' \prec t'' \prec \dotso$, contradicting our assumption that $\fS$ is locally bounded. So there exists $t \in \fS$ such that $st$ is maximal among $\menge{r \in \fS}{r \preceq s t_i \text{ for almost all } i}$. It follows that $\lim_i s t_i = s t$. If $t = t_j$ for some $j$, then $s t = s t_j \preceq s t_i$ implies $t_j \preceq t_i$ and hence $t_j = t_i$ by minimality. But $t_i \neq t_j$ for all $i \neq j$. So $s t \prec s t_i$ for all $i$. Since $(\fT_1,\fD_1)$ is $\max_{\preceq}^{\infty}$-closed, we must have $st \in \fT_1$. This contradicts minimality of $t_i$ unless $t = \mfd(s)$, which would contradict $s \notin \fT_1$. So we conclude that $\# \gekl{t_i} < \infty$, say $\gekl{t_i} = \gekl{t_1, \dotsc, t_j}$.

Define $V \defeq \menge{\chi \in X_2(\mfd(s))}{\chi(t_i) = 0 \ \forall \ 1 \leq i \leq j}$. Given a normal word $z$ in $\fT_2$ such that $\chi_z \in V$, we claim that $cpsz$ is normal. Indeed, this follows from $H(sz) = s$, which is shown as follows: If $H(sz) = st \in \fS$, then $s t \in \fT_2$ because $(\fT_2,\fD_2)$ is $H$-invariant. If $st \neq s$, then $st \in \fT_1$ because $s$ is maximal in $\fT_2 \setminus \fT_1$. Now $H(sz) = st$ implies that $\chi_{sz}(st) = 1$ and hence $\chi_z(t) = 1$. But $\chi_z \in V$ implies that $\chi_z(t) = 0$, which is a contradiction. Moreover, $H(sz) = s$ and $s \in \fT_2 \setminus \fT_1$ imply that $z$ must be a normal word in $\fT_2 \setminus \fT_1$. It follows that $cps.V \subseteq X_2 \setminus X_1$. In addition, we have $V \neq \emptyset$. Indeed, as $s$ lies in $\fT_2 \setminus \fT_1$, there exists a normal word $w$ starting with $s$, say $w = s s_1 s_2 \dotsm$, such that $\chi_w \in X_2$. Then $s^{-1}.\chi_w = \chi_{s^{-1}w}$, where $s^{-1}w = s_1 s_2 \dotsm$, must lie in $V$ because $s^{-1}.\chi_w(t_i) = 1$ would imply $t_i \preceq s_1 \dotsm s_n$ for some $n$ and thus $st_i \preceq s s_1 \dotsm s_n$. But this, together with $st_i \in \fS$, would contradict that $s s_1 s_2 \dotsm$ is normal.

Our findings imply that $cps. V \subseteq U$. Now the bisection $[cpsq_1, cps.V]$ has source $cps.V \subseteq U$ and range $cpsq_1cps.V \subsetneq cps.V$ because $(cpsq_1cps.V) \cap (cpsq_2cps.V) = \emptyset$.
\eproof
\setlength{\parindent}{0cm} \setlength{\parskip}{0cm}

We derive the following consequences with the help of \cite{BL}.
\bcor
\label{cor:InvSub=Ideals,StronglyPurelyInf}
Suppose that $I_l \ltimes \Omega$ is Hausdorff, inner exact in the sense of \cite[Definition~3.7]{Ana} and \cite[Definition~3.5]{BL}, $\fC^{*,0} = \fC^0$, and that every admissible, $H$-invariant and $\max_{\preceq}^{\infty}$-closed pair $(\fT,\fD)$ has the property that for all finite normal paths $a, b$ in $\fS$ with $\mft(a) = \mft(b)$, $\mfd(a) = \mfd(b)$, $a_1 \neq b_1$ and $s \in \fT$ with $\mfd(s) = \mfd(a)$, there exists $t \in \fT$ such that $s t$ is normal and $H(at) \neq H(bt)$. Then the map $(\fT,\fD) \mapsto \spkl{C_0(\Omega \setminus X(\fT,\fD))}$ establishes an inclusion-preserving one-to-one correspondence between admissible, $H$-invariant, $\max_{\preceq}^{\infty}$-closed pairs $(\fT,\fD)$ with $\fT \subseteq \fS$ and $\fD \subseteq \fC^0$ and ideals of $C^*_r(I_l \ltimes \Omega)$.
\setlength{\parindent}{0.5cm} \setlength{\parskip}{0cm}

Suppose, in addition, that $\fC$ is left Noetherian, and that every admissible, $H$-invariant and $\max_{\preceq}^{\infty}$-closed pair $(\fT,\fD)$ satisfies the following: For every admissible, $H$-invariant and $\max_{\preceq}^{\infty}$-closed pair $(\fT',\fD')$ with $(\fT',\fD') \subsetneq (\fT,\fD)$, there exists an infinite normal path in $\fT \setminus \fT'$, and for all finite normal paths $c$ in $\fT$, there exists a maximal element $s \in \fT$ with respect to $\preceq$ together with a normal path $p$ such that $cps$ is normal and two distinct normal paths $q_1$, $q_2$ such that $s q_1 c$ and $s q_2 c$ are normal. Then $C^*_r(I_l \ltimes \Omega)$ is strongly purely infinite.
\ecor
\bproof
Proposition~\ref{prop:Gar-TopFree} implies that $I_l \ltimes \Omega$ is essentially principal, in the sense of \cite[\S~2.1]{BL}. Now our first claim follows from \cite[Corollary~3.12]{BL}. The second claim follows from Proposition~\ref{prop:XminusXLocallyContractive} and \cite[Theorem~4.2]{BL}.
\eproof
\setlength{\parindent}{0cm} \setlength{\parskip}{0.5cm}

\section{Examples}

We apply our findings to two concrete classes of examples, higher rank graphs and Artin-Tits monoids.

\subsection{Higher rank graphs}
\label{ss:k-graphs}

Let $P = \Zz_0^k$, where $\Zz_0 = \gekl{0, 1, 2, 3, \dotsc}$ denotes the set of non-negative integers. A higher rank graph is a small category $\fC$ equipped with a $P$-valued degree map, i.e., a functor $\bbd: \: \fC \to P$ such that the following unique factorization property holds: For all $c \in \fC$ with $\bbd(c) = pq$, there exist $a, b \in \fC$ with $c = ab$, $\bbd(a) = p$, $\bbd(b) = q$, and if we have $c = a' b'$ for some $a', b' \in \fC$ with $\bbd(a') = p$, $\bbd(b') = q$, then $a' = a$ and $b' = b$. Note that $\fC$ is automatically cancellative, and we have $\fC^* = \fC^0$.

C*-algebras attached to higher rank graphs have been introduced in \cite{KuPa}. Given a higher rank graph $\fC$, its C*-algebra in the sense of \cite{KuPa,RSY03,RSY04} is canonically isomorphic to $C^*_r(I_l \ltimes \bOmega)$ (see \cite{FMY,Sp14}). In the following, we want to apply our findings to higher rank graphs. In particular, our goal is to classify closed invariant subspaces of $\Omega$ and to compare our results with previous work.

First, we need to find a Garside family in $\fC$. This will be discussed in \cite[\S~6]{Li21b} in more detail and in a more general context. In the following, we simply summarize what we need in our specific situation. First, since $P$ is Noetherian, then so is $\fC$. Moreover, $\fC$ has disjoint mcms in the following sense: Given $\mfv \in \fC^0$ and $a, b \in \mfv \fC$, take $C \subseteq \mfv \fC$ such that the canonical projection $\fC \to \fC / { }_{\sim}$ induces a bijection $C \isom (\bbd^{-1}(\lcm(a,b)) \cap (a \fC \cap b \fC))$. Then
$
 a \fC \cap b \fC = \coprod_{c \in C} c \fC
$.
Let $S_P \defeq \menge{(0, \dotsc, 0) \neq (p_1, \dotsc, p_k) \in P}{0 \leq p_j \leq 1 \ \forall \ 1 \leq j \leq k}$. Let $\fS \defeq \bbd^{-1}(S_P)$. Then $\fS$ is a Garside family in $\fC$ which is always locally bounded.

In the following, we will always assume that $\fC$ is finitely aligned (which is not automatic in general). To give an example for a sufficient condition, if $\mfv \bbd^{-1}(p) < \infty$ for all $\mfv \in \fC^0$ and $p \in P$, then $\fC$ is finitely aligned. Actually, in that case $\fS$ will be locally finite. Note that in the literature, higher rank graphs with locally finite $\fS$ are called row-finite.

Let us now apply the classification of closed invariant subspaces of $\Omega$ in our situation. First observe that given $s, t \in \fS$, $s t$ is normal if and only if $\bbd(s) \geq \bbd(t)$. Moreover, $a \in \fS$ is an atom if and only if $\bbd(a)$ is one of the standard generators of $P$. The following are easy to see.
\blemma
$(\fT,\fD) \subseteq (\fS,\fC^0)$ is admissible if and only if the following is satisfied:
\setlength{\parindent}{0cm} \setlength{\parskip}{0cm}

\begin{enumerate}
\item[(A)] For every $t \in \fT$ there exists $t' \in \fT$ with $\bbd(t) \geq \bbd(t')$ or $\mfd(t) \in \fD$. 
\end{enumerate}
$(\fT,\fD) \subseteq (\fS,\fC^0)$ is $H$-invariant if and only if the following is satisfied:
\begin{enumerate}
\item[(I)] For every $t \in \fT \cup \fD$ and every atom $a$ with $\mfd(a) = \mft(t)$, if $\bbd(a) \not\leq \bbd(t)$, then $at \in \fT$, and if $\bbd(a) \leq \bbd(t)$ and $t = rs$ with $\bbd(s) = \bbd(a)$, then $ar \in \fT$. 
\end{enumerate}
$(\fT,\fD) \subseteq (\fS,\fC^0)$ is $\max_{\preceq}^{\infty}$-closed if and only if the following is satisfied:
\begin{enumerate}
\item[(C)] Given a sequence $a z_i \in \fT$ with $\bbd(z_i) = d \in P$, if whenever $\varepsilon \leq d$ is a standard generator of $P$ and $s_i \preceq z_i$ satisfies $\bbd(s_i) = \varepsilon$, we must have $s_i \neq s_j$ for all $i \neq j$, then $a \in \fT \cup \fD$.
\end{enumerate}
\elemma
\setlength{\parindent}{0cm} \setlength{\parskip}{0cm}

With these observations, we obtain the following applications of Theorem~\ref{thm:ClInvSubsp} and Corollary~\ref{cor:ClInvSubsp_locfin}. As before, we write $\cG \defeq I_l \ltimes \Omega$.
\bcor
\label{cor:ClInv_HigherRank}
Suppose that $\fC$ is a countable, finitely aligned higher rank graph. Then the maps $X \ma (\fT(X),\fD(X))$ and $X(\fT,\fD) \mapsfrom (\fT,\fD)$ establish an inclusion-preserving one-to-one correspondence between $\cG$-invariant, closed subspaces of $\Omega$ and pairs $(\fT,\fD)$ with $\fT \subseteq \fS$ and $\fD \subseteq \fC^0$ satisfying conditions (A), (I) and (C).
\setlength{\parindent}{0.5cm} \setlength{\parskip}{0cm}

If, in addition, $\fS$ is locally finite, then the maps $X \ma (\fT(X),\fD(X))$ and $X(\fT,\fD) \mapsfrom (\fT,\fD)$ establish an inclusion-preserving one-to-one correspondence between $\cG$-invariant, closed subspaces of $\Omega$ and pairs $(\fT,\fD)$ with $\fT \subseteq \fS$ and $\fD \subseteq \fC^0$ satisfying conditions (A) and (I).
\ecor
\setlength{\parindent}{0cm} \setlength{\parskip}{0cm}

Moreover, Proposition~\ref{prop:Gar-TopFree} yields the following sufficient condition for topological freeness, and hence also effectiveness since our groupoids are Hausdorff (as $\fC$ is cancellative and finitely aligned).
\bcor
Suppose that $\fC$ is a countable, finitely aligned higher rank graph. Let $(\fT,\fD)$ satisfy (A), (I) and (C). Set $X \defeq X(\fT,\fD)$. Assume that for all finite normal paths $a, b$ in $\fS$ with $\mft(a) = \mft(b)$, $\mfd(a) = \mfd(b)$, $a_1 \neq b_1$ and $d \in \bbd(\fT)$, there exists $t \in \fT$ such that $\bbd(t) = d$ and $H(at) \neq H(bt)$. Then $I_l \ltimes X$ is effective.
\ecor

Let us now consider the boundary. The following are consequences of Corollary~\ref{cor:ClInv-B} and Lemma~\ref{lem:TMax_locfin}.
\bcor
\label{cor:B_HigherRank}
Suppose that $\fC$ is a countable higher rank graph which is locally convex such that $\fS$ is locally finite. Then, for all $\mfv \in \fC^0$,
$$
 \mfv \fT_{\rm Max} = \mfv \fT(\Omegamax) = \menge{s \in \mfv \fS}{\bbd(s) = \max \menge{\bbd(r)}{r \in \mfv \fS}}.
$$
Under the correspondence in Corollary~\ref{cor:ClInv_HigherRank}, a $\cG$-invariant, closed subspace $X$ is contained in $\bOmega$ if and only if $(\fT(X),\fD(X)) \subseteq (\fT_{\rm Max},\fD_{\max})$.
\ecor
\setlength{\parindent}{0cm} \setlength{\parskip}{0cm}

Let us compare this last result with the classification of gauge-invariant ideals of higher rank graph C*-algebras in \cite{RSY03}. In the following, we write $\bcG \defeq I_l \ltimes \bOmega$.
\blemma
\label{lem:gauge=induced}
Let $\fC$ be a countable, finitely aligned higher rank graph. An ideal $I$ of the C*-algebra $C^*_r(\bcG)$ of $\fC$ is gauge-invariant if and only if it is induced from an open invariant subspace of $\bOmega$, in the sense that $I = \spkl{C_0(U)} = C^*_r(I_l \ltimes U)$ for some open, $\bcG$-invariant subspace $U \subseteq \bOmega$.
\elemma
\bproof
Let $\theta$ be the canonical conditional expectation on $C^*_r(\bcG)$ given by averaging over the gauge action. If $I$ is gauge-invariant, then $\theta(I) \subseteq I$. Now the image of $\theta$ is given by $C^*_r(\cH)$, where $\cH$ is the subgroupoid of $\bcG$ describing $C^*_r(\bcG)^{\theta}$ identified in \cite{RSWY}. As $\cH$ is principal, it follows from \cite{BL} that the ideal $\theta(I)$ is induced, i.e., if $\theta': \: C^*_r(\cH) \onto C_0(\bOmega)$ is the canonical conditional expectation, then $\theta'(\theta(I)) \subseteq \theta(I)$. So we conclude that $\spkl{\theta'(\theta(I))} \subseteq I$. Now $\theta' \circ \theta$ is the canonical conditional expectation $C^*_r(I_l \ltimes \bOmega) \onto C_0(\bOmega)$, and we always have $I \subseteq \spkl{\theta'(\theta(I))}$. As $\theta'(\theta(I))$ is an ideal of $C_0(\bOmega)$, it must be of the form $C_0(U)$ for some open, $\bcG$-invariant subspace $U \subseteq \bOmega$.
\eproof

In \cite{RSY03}, it was shown that gauge-invariant ideals of $C^*_r(\bcG)$ are in one-to-one correspondence to hereditary, saturated subsets $\fH \subseteq \fC^0$. The following result enables us to translate between this result in \cite{RSY03} and Corollary~\ref{cor:B_HigherRank}.
\blemma
The assignment $(\fT,\fD) \ma \fH(\fT,\fD) \defeq \menge{\mfv \in \fC^0}{\bOmega(\mfv) \subseteq \bOmega \setminus X(\fT,\fD)}$ defines a one-to-one correspondence between pairs $(\fT,\fD)$ with $(\fT(X),\fD(X)) \subseteq (\fT_{\rm Max},\fD_{\max})$ satisfying conditions (A) and (I) and hereditary, saturated subsets $\fH \subseteq \fC^0$.
\elemma
\bproof
First of all, $\fH(\fT,\fD)$ is hereditary. Indeed, given $a \in \fC$ with $\mft(a) = \mfw$ and $\mfd(a) \mfv$, where $\bOmega(\mfw) \in \bOmega \setminus X(\fT,\fD)$, take an infinite normal word $x$ with $\chi_x \in \bOmega$ and $\mft(x) = \mfv$. Then $a.\chi_x \in \bOmega \setminus X(\fT,\fD)$ as $\mft(ax) = \mfw$. As $\bOmega \setminus X(\fT,\fD)$ is $\bcG$-invariant, it follows that $\chi_x \in \bOmega \setminus X(\fT,\fD)$, as desired. Moreover, $\fH(\fT,\fD)$ is saturated because we have $\bigcup_{\mfv \in \fH(\fT,\fD)} \fC \mfv \bOmega \subseteq \bOmega \setminus X(\fT,\fD)$.
\setlength{\parindent}{0.5cm} \setlength{\parskip}{0cm}

Now given a hereditary, saturated subset $\fH \subseteq \fC^0$, define $\cO(\fH) \defeq \bigcup_{\mfv \in \fH} \fC \mfv \bOmega$. We claim that $\bOmega \setminus X(\fT,\fD) = \cO(\fH(\fT,\fD))$. \an{$\supseteq$} is clear. To show \an{$\subseteq$}, take an infinite normal word $x = s_1 s_2 \dotsm$ with $\chi_x \in \bOmega \setminus X(\fT,\fD)$. Then there exists $n$ such that $s_{n+1} \notin \fT$ and thus $s_N \notin \fT$ for all $N \geq n+1$. We then claim that $\mfv = \mfd(s_{n+1}) \in \fH(\fT,\fD)$. Indeed, if there exists $t \in \fT$ with $\mft(t) = \mfv$, then $\bbd(t) = \max \menge{\bbd(r)}{r \in \mfv \fS}$. It follows that $\bbd(t) = \bbd(s_{n+2})$. Hence $s_{n+1} t$ is normal. Since $(\fT,\fD)$ is $H$-invariant, that would imply $s_{n+1} \in \fT$, which is a contradiction. This shows that $\chi_x \in (s_1 \dotsm s_n s_{n+1}). \partial \Omega(\mfv) \subseteq \cO(\fH(\fT,\fD))$, as desired. The conclusion is that the map $(\fT,\fD) \ma \fH(\fT,\fD) \ma \cO(\fH(\fT,\fD))$ is a bijection between pairs $(\fT,\fD)$ with $(\fT(X),\fD(X)) \subseteq (\fT_{\rm Max},\fD_{\max})$ satisfying conditions (A) and (I) and $\bcG$-invariant, open subsets of $\bOmega$.

Hence it suffices to show that the map $\fH \ma \cO(\fH)$ is injective. We claim that for every hereditary, saturated subset $\fH \subseteq \fC^0$, we have $\fH = \menge{\mfv \in \fC^0}{\bOmega(\mfv) \subseteq \cO(\fH)}$. Indeed, \an{$\subseteq$} is clear, and for \an{$\supseteq$}, suppose that $\mfw \in \fC^0$ satisfies $\bOmega(\mfw) \subseteq \cO(\fH) = \bigcup_{\mfv \in \fH} \fC \bOmega(\mfv)$. Then, by compactness of $\bOmega(\mfw)$, there exist finitely many $a_i \in \fC$ and $\mfv_i \in \fH$ with $\bOmega(\mfw) \subseteq \bigcup_{i=1}^n a_i. \bOmega(\mfv_i)$. It follows that $\gekl{a_i}$ must be exhaustive in the sense of \cite{Sim}, and thus $\mfw \in \fH$ because $\fH$ is saturated.
\eproof
\setlength{\parindent}{0cm} \setlength{\parskip}{0cm}

\bremark
It would also be interesting to compare our results with the ones in \cite{Sim} for more general finitely aligned higher rank graphs.
\eremark

Finally, we specialize to one vertex higher rank graphs.
\bcor
Suppose that $\fC$ is a countable, finitely aligned higher rank graph with one vertex. In that case $\fS$ is locally finite if and only if it is finite.
\setlength{\parindent}{0.5cm} \setlength{\parskip}{0cm}

If $\fS$ is finite, then Corollary~\ref{cor:ClInv_HigherRank} yields a one-to-one correspondence between $\geq$-closed subsets $T \subseteq S_P$ and $\cG$-invariant, closed subspaces of $\Omegainf$, given by $T \ma X(\bbd^{-1}(T),\emptyset)$. Moreover, if for every standard generator $\varepsilon$ of $P$, $\# \bbd^{-1}(\varepsilon) \geq 2$, then $I_l \ltimes (X_2 \setminus X_1)$ is locally contractive for all $\cG$-invariant, closed subspaces $X_1 \subsetneq X_2$ of $\Omegainf$. Furthermore, if for all finite normal paths $a, b$ in $\fS$ with $\mft(a) = \mft(b)$, $\mfd(a) = \mfd(b)$, $a_1 \neq b_1$ and $d \in T$, there exists $t \in \fC$ such that $\bbd(t) = d$ and $H(at) \neq H(bt)$, then $I_l \ltimes X(T)$ is effective.

If $\fS$ is infinite, then Corollary~\ref{cor:ClInv_HigherRank} yields a one-to-one correspondence between $\geq$-closed subsets $T \subseteq S_P$ and $\cG$-invariant, closed subspaces of $\Omega$, given by $T \ma X(T) \defeq X(\bbd^{-1}(T),\fC^0)$. Moreover, if for every standard generator $\varepsilon$ of $P$, $\# \bbd^{-1}(\varepsilon) \geq 2$, then $I_l \ltimes (X_2 \setminus X_1)$ is locally contractive for all $\cG$-invariant, closed subspaces $X_1 \subsetneq X_2$ of $\Omega$. Furthermore, if for all finite normal paths $a, b$ in $\fS$ with $\mft(a) = \mft(b)$, $\mfd(a) = \mfd(b)$, $a_1 \neq b_1$ and $d \in T$, there exists $t \in \fC$ such that $\bbd(t) = d$ and $H(at) \neq H(bt)$, then $I_l \ltimes X(T)$ is effective.
\ecor
\setlength{\parindent}{0cm} \setlength{\parskip}{0cm}

\bremark
The general results in Theorem~\ref{thm:ClInvSubsp} and Corollary~\ref{cor:ClInvSubsp_locfin} can also be applied to monoids and categories arising from self-similar actions of groups and groupoids on graphs and higher rank graphs as in \cite{EP17, LRRW, ABRW, BKQS} (see \cite[Remark~7.10]{Li21b}).
\eremark
\setlength{\parindent}{0cm} \setlength{\parskip}{0.5cm}

\subsection{Artin-Tits monoids}

In the following, we analyse reduced C*-algebras of Artin-Tits monoids using our general approach for Garside categories. Recall that an Artin-Tits monoid $P$ is given by the following presentation: 
$$
 P = \big\langle A \ \big\vert \ (ab)^{[m_{a,b}]} = (ba)^{[m_{b,a}]} \ \forall \ a, b \in A \big\rangle ^+,
$$
where $A$ is a set (the set of atoms), $m_{a,b} \in \gekl{2, 3, \dotsc} \cup \gekl{\infty}$ with $m_{a,b} = m_{b,a}$, and $(ab)^{[m_{a,b}]}$ denotes the alternating word $abab \dotsm$ of length $m_{a,b}$. If $m_{a,b} = m_{b,a} = \infty$, then it is understood that $(ab)^{[m_{a,b}]} = (ba)^{[m_{b,a}]}$ simply means that we do not add a relation involving $a$ and $b$. For more information about Artin-Tits monoids and groups, the reader may consult for instance \cite{BS,CL02,CL07}. In the following, given $x \in P$, we write $\cL(x) \defeq \menge{a \in A}{a \preceq x}$ and $\cR(x) \defeq \menge{a \in A}{a \tipreceq x}$. It was recently shown in \cite{DH} (see also \cite{DDH}) that there exists a finite Garside family in every finitely generated Artin-Tits monoid. 

\btheo
\label{thm:AT_fingen_NotSpherical}
Let $P$ be a finitely generated, irreducible, non-spherical Artin-Tits monoid and $S \subseteq P$ a finite Garside family. Suppose $T \subseteq S$ is such that $(T,\emptyset)$ is admissible, $H$-invariant and $\max_{\preceq}^{\infty}$-closed. Then $A \subseteq T$ and hence $T = S$.
\etheo
\setlength{\parindent}{0cm} \setlength{\parskip}{0cm}

\bproof
Let $A^s \subseteq A$ be maximal such that $\Delta_{A^s} \defeq \lcm \menge{a}{a \in A^s}$ exists. By assumption, $A^s \neq A$ because otherwise, $P$ would be spherical by \cite{BS}. Take $t \in T$ arbitrary and form $x_1 \defeq \Delta_{A^s} t$. We must have $A^s = \cL(x_1)$. Indeed, if $b \in \cL(x_1)$ and $b \notin A^s$, then $\lcm(b, \Delta_{A^s})$ would exist, contradicting maximality of $A^s$. Let $A_1, \dotsc, A_n$ be the irreducible components of $A^s$. We proceed inductively on $n$. Since $A^s \neq A$, there exist $a_1 \in A_1$ and $a_2 \in A_2$ together with $b_1, \dotsc, b_N \in A \setminus A^s$ such that $m_{a_1,b_1}, m_{b_1,b_2}, \dotsc, m_{b_{N-1},b_N}, m_{b_N,a_2} > 2$. For all $1 \leq m \leq n$, set $\Delta_m \defeq \lcm \menge{a}{a \in A_m}$. We have $a_1 \not\preceq b_1 \Delta_1 \Delta_2 x_1$: If $a_1 \preceq b_1 \Delta_1 \Delta_2 x_1$, then $m_{a_1,b_1} > 2$ implies that $b_1 a_1 b_1 \preceq \lcm(b_1,a_1) \preceq b_1 \Delta_1 \Delta_2 x_1$ and thus $b_1 \preceq a_1^{-1} \Delta_1 \Delta_2 x_1$. We claim that this would imply $b_1 \preceq x_1$. Indeed, write $a_1^{-1} \Delta_1 \Delta_2 = c_1 \dotsm c_l$ for some $c_1, \dotsc, c_l \in A_1 \cup A_2$. $b_1 \preceq c_1 \dotsm c_l x_1$ implies $c_1 b_1 \preceq \lcm(c_1,b_1) \preceq c_1 \dotsm c_l x_1$ and thus $b_1 \preceq c_2 \dotsm c_l x_1$. Now proceeding inductively, we end up with $b_1 \preceq x_1$, which is a contradiction. So $a_1 \not\preceq b_1 \Delta_1 \Delta_2 x_1$. Similarly, we obtain $b_1, a_1 \not\preceq b_2 b_1 \Delta_1 \Delta_2 x_1$, ..., and finally $b_N, b_{N-1}, \dotsc, b_1, a_1 \not\preceq a_2 b_N b_{N-1} \dotsm b_2 b_1 \Delta_1 \Delta_2 x_1$. Set $x_2 \defeq a_2 b_N b_{N-1} \dotsm b_2 b_1 \Delta_1 \Delta_2 x_1$. We conclude that $\cL(x_2) \subseteq A_1 \setminus \gekl{a_1} \cup A_2 \cup \dotsm \cup A_n$. By \cite[Proposition~4.38]{LOS}, there exists a normal path $g_1 \dotsm g_k$ in $A_1$ with $\cL(g_1) = \gekl{a_1}$ and $\cL(g_k) = A_1 \setminus \gekl{a_1}$. Define $\Delta^{(1)} \defeq \Delta_2 \dotsm \Delta_n$ and $A^{(1)} \defeq A_2 \cup \dotso \cup A_n$. Then $\cL(g_j \Delta^{(1)}) = \cL(g_j) \cup A^{(1)}$ and $\cR(g_j \Delta^{(1)}) = \cR(g_j) \cup A^{(1)}$ for all $1 \leq j \leq k$. If we now set $g'_j \defeq g_j \Delta^{(1)}$, then $g'_1 \dotsm g'_k x_2$ is normal. With $x_3 \defeq g'_1 \dotsm g'_k x_2$, we obtain $\cL(x_3) = \gekl{a_1} \cup A^{(1)}$. Let $b_1, \dotsc, b_N$ be as above. We have $a_1 \not\preceq b_1 a_1 x_3$ as $b_1 \not\preceq x_3$, and proceeding inductively, we arrive at $b_N, \dotsc, b_1, a_1 \not\preceq a_2 b_N \dotsc b_1 a_1 x_3$. Hence, with $x_4 \defeq a_2 b_N \dotsc b_1 a_1 x_3$, we obtain $\cL(x_4) \subseteq A^{(1)}$. Repeating this process, we arrive at an element $x$ of the form $pt$ for some $p \in P$ with $\cL(x) = \gekl{a}$ for some $a \in A$. Now suppose that $a' \in A$ is arbitrary. Since $P$ is irreducible, there exist $d_1, \dotsc, d_M \in A$ such that $m_{a,d_1}, m_{d_1,d_2}, \dotsc, m_{d_{M-1},d_M}, m_{d_M,a'} > 2$. An analogous argument as above shows that $d_M, \dotsc, d_1, a \not\preceq a' d_M \dotsm d_1 a x$. Hence, with $y \defeq a' d_M \dotsm d_1 a x$, $y$ is of the form $qt$ for some $q \in P$, and we have $\cL(y) = \gekl{a'}$. It follows that $H(a' y) = a' \in T$, as desired.
\eproof

Recall that $P$ is called left reversible if $pP \cap qP \neq \emptyset$ for all $p, q \in P$. If our irreducible Artin-Tits monoid $P$ is not finitely generated, $P$ is left reversible if and only if $P$ is the increasing union of finitely generated, irreducible, spherical Artin-Tits submonoids. 

\btheo
\label{thm:InfGenNotRev_min}
Let $P$ be an irreducible Artin-Tits monoid which is not finitely generated. If $P$ is not left reversible, then $\Omega$ is minimal.
\etheo
\bproof
Given $F \subseteq A$, we write $P_F \defeq \spkl{F}^+ \subseteq P$. Let $U = \Omega(x; \mfy)$ be a basic open set and $F \subseteq A$ a finite subset with $x \in P_F$, $\mfy \subseteq P_F$. As $P$ is not left reversible, we may assume that $P_F$ is not left reversible and thus not spherical. Moreover, as $P$ is irreducible, there exists a finite subset $\bar{F} \subseteq A$ with $F \subsetneq \bar{F}$ such that $P_{\bar{F}}$ is irreducible. Since $P_F$ is not spherical, $P_{\bar{F}}$ is not spherical, either. This follows from the fact that Artin-Tits presentations are complete for right reversing by \cite{Deh11}, which implies that lcms in $P_F$ of elements in $P_F$ coincide with their lcms in $P_{\bar{F}}$ (see \cite[Proposition~6.10]{Deh03}). Now $F \subsetneq \bar{F}$ implies that $U \cap \Omega_{\bar{F},\infty} \neq \emptyset$. Hence, because $P_{\bar{F}}$ is irreducible and not spherical, Theorem~\ref{thm:AT_fingen_NotSpherical} implies that $U \cap \bOmega_{\bar{F}} \neq \emptyset$ and thus $U \cap \Omega_{\bar{F},\max} \neq \emptyset$. By Lemma~\ref{lem:Omegamax}, there exists $z \in P_{\bar{F}}$ with $x \preceq z$ and $z \cap y = \emptyset$ in $P_{\bar{F}}$ for all $y \in \mfy$. Now we again use that Artin-Tits presentations are complete for right reversing by \cite{Deh11}, so that lcms in $P_{\bar{F}}$ of elements in $P_{\bar{F}}$ coincide with their lcms in $P$ by \cite[Proposition~6.10]{Deh03}. It follows that $z \cap y = \emptyset$ in $P$ for all $y \in \mfy$. By Lemma~\ref{lem:Omegamax}, there exists $\chi \in \Omegamax$ with $\chi(z) = 1$. Thus we have found a character $\chi$ in $U \cap \Omegamax$. Hence $\bOmega = \Omega$. 
\eproof

Let $\kerbd$ be the boundary ideal, i.e., the kernel of the canonical quotient map $C^*_{\lambda}(P) \onto \partial C^*_{\lambda}(P)$.
\btheo
\label{thm:AT_NotFinGen_Rev}
Let $P$ be an irreducible Artin-Tits monoid which is not finitely generated. If $P$ is left reversible, then $\kerbd$ is simple. In particular, $\Omega \setminus \gekl{\infty}$ is minimal.
\etheo
\bproof
Let $F \subseteq A$ and $P_F \defeq \spkl{F}^+ \subseteq P$. Let $\lambda$ be the left regular representation of $P$. First we want to identify $C^*_{\lambda}(P_F)$ with the sub-C*-algebra $C^*(\lambda(P_F))$ of $C^*_{\lambda}(P)$ generated by the image of $P_F$ under $\lambda$. The inclusion $P_F \subseteq P$ allows us to view $\ell^2 P_F$ as a subspace of $\ell^2 P$. Restriction to this subspace induces a homomorphism $C^*(\lambda(P_F)) \onto C^*_{\lambda}(P_F)$, which restricts to an isomorphism of diagonal subalgebras $D(\lambda(P_F)) \isom C(\Omega_{P_F})$, where $D(\lambda(P_F)) = C^*(\lambda(P_F)) \cap C(\Omega)$. Using the commutative diagram 
\begin{equation*}
\begin{tikzcd}
  C^*(\lambda(P_F)) \ar["i_d"]{r} \ar{d} & C^*_{\lambda}(P_F) \ar{d}\\
  D(\lambda(P_F)) \ar{r} & C(\Omega_{P_F})
\end{tikzcd}
\end{equation*}
whose vertical arrows are the canonical faithful conditional expectations, a standard argument shows that the homomorphism $C^*(\lambda(P_F)) \onto C^*_{\lambda}(P_F)$ we just constructed is an isomorphism. Identifying $C^*(\lambda(P_F))$ with $C^*_{\lambda}(P_F)$, we obtain $C^*_{\lambda}(P) = \overline{\bigcup_F C^*_{\lambda}(P_F)}$, where $F$ runs through an increasing and exhausting family of finite subsets of $A$ such that $P_F$ is irreducible and spherical. Using for instance \cite[II.8.2.4]{Bla}, we deduce that $\kerbd = \overline{\bigcup_F {\rm Ker}_{\partial,F}}$. Now let $J$ be an ideal of $\kerbd$. Then $J = \overline{\bigcup_F J \cap {\rm Ker}_{\partial,F}}$ (again by \cite[II.8.2.4]{Bla}). By \cite[Theorem~4.39]{LOS}, we must have $J \cap {\rm Ker}_{\partial,F} \in \gekl{(0), \cK_F, {\rm Ker}_{\partial,F}}$, where $\cK_F = \cK(\ell^2 P_F)$. If $J \cap {\rm Ker}_{\partial,F} = (0)$ for all $F$, then clearly $J = (0)$. If there exists $F$ with $J \cap {\rm Ker}_{\partial,F} \neq (0)$, then $\cK_F \subseteq J \cap {\rm Ker}_{\partial,F}$, and it follows that for all $\bar{F}$ with $F \subsetneq F$, we have $J \cap {\rm Ker}_{\partial,\bar{F}} = {\rm Ker}_{\partial,\bar{F}}$ because $\spkl{\cK_F}_{{\rm Ker}_{\partial,\bar{F}}} \neq (0), \cK_{\bar{F}}$ and hence ${\rm Ker}_{\partial,\bar{F}} = \spkl{\cK_F}_{{\rm Ker}_{\partial,\bar{F}}} \subseteq J \cap {\rm Ker}_{\partial,\bar{F}}$. In that case, we conclude that $J = \kerbd$, as desired.
\eproof

\bcor
\label{cor:AT_ClInvSubsp}
Let $P$ be an irreducible Artin-Tits monoid. If $P$ is a finitely generated and spherical, then $\Omega$, $\Omegainf$ and $\bOmega = \gekl{\infty}$ are the only closed invariant subspaces of $\Omega$. If $P$ is finitely generated and not spherical, then $\Omega$ and $\Omegainf = \bOmega$ are the only closed invariant subspaces of $\Omega$. If $P$ is not finitely generated and left reversible, then $\Omega$ and $\bOmega = \gekl{\infty}$ are the only closed invariant subspaces of $\Omega$. If $P$ is not finitely generated and not left reversible, then $\Omega$ is minimal.
\ecor

\blemma
\label{lem:topfree-everywhere}
Let $P$ be an irreducible Artin-Tits monoid. If $P$ is left reversible, then $G$ acts topologically freely on every closed invariant subspace of $\Omega \setminus \gekl{\infty}$. If $P$ is not left reversible, then $G$ acts topologically freely on every closed invariant subspace of $\Omega$.
\elemma
\bproof
If $P$ is finitely generated and spherical, then the only closed invariant subspaces of $\Omega \setminus \gekl{\infty}$ are $\Omega \setminus \gekl{\infty}$ and $\Omegainf \setminus \gekl{\infty}$. On the first one, the $G$-action is topologically free because $P^* = \gekl{1}$ implies that $G \curvearrowright \Omega$ is topologically free. Topological freeness of $G \curvearrowright \Omegainf \setminus \gekl{\infty}$ follows from (the proof of) \cite[Theorem~4.39]{LOS}. If $P$ is not finitely generated and left reversible, then $\Omega \setminus \gekl{\infty}$ is minimal, and $G \curvearrowright \Omega \setminus \gekl{\infty}$ is topologically free because $P^* = \gekl{1}$ implies that $G \curvearrowright \Omega$ is topologically free. If $P$ is finitely generated and not spherical, then the only closed invariant subspaces of $\Omega$ are $\Omega$ and $\Omegainf = \bOmega$. $G \curvearrowright \Omega$ is topologically free because $P^* = \gekl{1}$. Moreover, it is shown in \cite{CL07} and also follows from Theorem~\ref{thm:Gc=1} that $G \curvearrowright \bOmega$ is topologically free if $P$ is right-angled, i.e., $m_{a,b} \in \gekl{2, \infty}$ for all $a, b \in A$. If $P$ is not right-angled, then we must have $\# A \geq 3$, so that we can find $a, b \in A$ with $2 < m_{a,b} < \infty$. Let $P_{a,b} \defeq \spkl{a,b}^+$ be the submonoid of $P$ generated by $a$ and $b$ and $G_{a,b}$ its enveloping group. $P_{a,b}$ is itself a spherical Artin-Tits monoid. Hence it follows from \cite[Remark~4.6]{LOS} that $G_{a,b} \curvearrowright (\Omega_{P_{a,b}})_{\infty} \setminus \bOmega_{P_{a,b}}$ is topologically free, where $\Omega_{P_{a,b}}$ denotes the space of characters for $P_{a,b}$ and $\bOmega_{P_{a,b}}$ its boundary. Thus there exists an infinite word $w$ in $a, b$ such that the corresponding character $\chi_{w,P_{a,b}}$ in $\Omega_{P_{a,b}}$ has trivial stabilizer group in $G_{a,b}$. Let $\chi_w \in \Omega$ be the character given by the same infinite word. If $g \in G$ satisfies $g.\chi_w = \chi_w$, then it follows that $g \in G_{a,b}$ and hence $g = 1$. This shows that $\chi_w$ also has trivial stabilizer group in $G$. Finally, if $P$ is not finitely generated and not left reversible, then $\Omega$ is minimal, and $G \curvearrowright \Omega$ is topologically free because $P^* = \gekl{1}$.
\eproof

\blemma
\label{lem:RevInfGen_pi}
Let $P$ be an irreducible Artin-Tits monoid which is not finitely generated and left reversible. Then $\kerbd$ is purely infinite.
\elemma
\bproof
For $F \subseteq A$, let $P_F \defeq \spkl{F}^+$ be the submonoid of $P$ generated by $F$, and let $\Omega_F$ be the space of characters for $P_F$ and $\bOmega_F = \gekl{\infty_F}$ its boundary. Recall from the proof of Theorem~\ref{thm:AT_NotFinGen_Rev} that $\kerbd = \overline{\bigcup_F {\rm Ker}_{\partial,F}}$, where $F$ runs through all finite subsets of $A$. Write $\ti{\Omega} \defeq \Omega \setminus \gekl{\infty}$ and $\ti{\Omega}_F \defeq \Omega_F \setminus \gekl{\infty_F}$. For a finite subset $F \subseteq A$ such that $P_F$ is irreducible and spherical, consider a basic compact open subset $U = \ti{\Omega}_F(x;\mfy)$ and $\bm{1}_U$ the corresponding characteristic function. We claim that for every finite subset $\bar{F}$ of $A$ with $F \subsetneq \bar{F}$  such that $P_{\bar{F}}$ is irreducible and spherical, $\bm{1}_U$ is infinite in ${\rm Ker}_{\partial,\bar{F}}$ and hence also in $\kerbd$. The image of $\bm{1}_U$ in ${\rm Ker}_{\partial,\bar{F}}$ is given by the characteristic function of $\ti{\Omega}_{\bar{F}}(x;\mfy)$. $F \subsetneq \bar{F}$ implies that $\ti{\Omega}_{\bar{F}}(x;\mfy) \cap (\Omega_{\bar{F}})_{\infty} \neq \emptyset$. Then the same argument for local contractiveness in \cite[Theorem~4.39]{LOS} shows that the image of $\bm{1}_U$ is infinite in ${\rm Ker}_{\partial,\bar{F}}$, as desired. Now \cite[Theorem~4.1]{BCS} (second countability is not needed, see \cite[Theorem~4.2]{BL}) implies that $\kerbd$ is purely infinite because sets of the form $\ti{\Omega}(x;\mfy)$ for $x \in P_F$, $\mfy \subseteq P_F$, where $F$ is an arbitrary finite subset of $A$ such that $P_F$ is irreducible and spherical, form a basis of compact open subsets of $\ti{\Omega}$.
\eproof

Let us summarize our analysis of left regular C*-algebras of Artin-Tits monoids. Our results on ideal structure and pure infiniteness extend the corresponding results in the right-angled case in \cite{CL02,CL07} and in the finitely generated, spherical case in \cite[\S~4.2]{LOS}. 
\bcor
Let $P$ be an irreducible Artin-Tits monoid. If $P$ is spherical, then $\kerbd = \cK(\ell^2 P)$ if $\# A = 1$ and $\cK(\ell^2 P)$ is the only non-trivial ideal of $\kerbd$ if $2 \leq \#A < \infty$. In the latter case, $\kerbd / \cK(\ell^2 P)$ is purely infinite simple. If $P$ is not finitely generated and left reversible, then $\kerbd$ is purely infinite simple. If $P$ is finitely generated and not spherical, then $\cK(\ell^2 P)$ is the only non-trivial ideal of $C^*_{\lambda}(P)$, and $C^*_{\lambda}(P) / \cK(\ell^2 P)$ is purely infinite simple. If $P$ is not finitely generated and not left reversible, then $C^*_{\lambda}(P)$ is purely infinite simple.
\setlength{\parindent}{0.5cm} \setlength{\parskip}{0cm}

If $P$ is left reversible, then $C^*_{\lambda}(P)$ is nuclear if and only if $\# A = 1$, and $\kerbd$ is nuclear if and only if $\# A \leq 2$. If $P$ is not left reversible, then $C^*_{\lambda}(P)$ is nuclear if and only if $P$ is right-angled (i.e., $m_{a,b} \in \gekl{2,\infty}$ for all $a, b \in A$).
\ecor
\setlength{\parindent}{0cm} \setlength{\parskip}{0cm}

\bproof
Our claims for finitely generated spherical $P$ follow from \cite[Remark~4.8, Theorem~4.39, Proposition~4.15]{LOS}. If $P$ is not finitely generated and left reversible, our claims follow from Lemma~\ref{lem:RevInfGen_pi} and the same argument as for \cite[Proposition~4.15]{LOS} for the failure of nuclearity. If $P$ is not left reversible, then our claims follow from Corollary~\ref{cor:AT_ClInvSubsp}, Lemma~\ref{lem:topfree-everywhere} and \cite[Corollary~5.7.17]{CELY} because $C^*_{\lambda}(P) / \cK(\ell^2 P) \cong \partial C^*_{\lambda}(P)$ if $P$ is finitely generated and $C^*_{\lambda}(P) \cong \partial C^*_{\lambda}(P)$ if $P$ is not finitely generated. Our claims about nuclearity follow from \cite[Theorem~4.2]{LL} or a similar argument as for \cite[Proposition~4.15]{LOS}.
\eproof
\setlength{\parindent}{0cm} \setlength{\parskip}{0.5cm}

\end{document}